%% file: amicable_numbers.tex
\documentclass[11pt]{article}
\usepackage{array}
\usepackage{multirow}
\usepackage{longtable}
\usepackage{tablefootnote}
\usepackage{helvet}
\usepackage{amsmath}
\usepackage{enumitem}
\usepackage{fancyvrb}
\usepackage[a4paper, portrait, margin=1.5in]{geometry}

\newcommand{\FR}{F.\,R.}
\newcommand{\JE}{J.\,E.}
\newcommand{\et}{\quad\mbox{and}\quad}
\newcommand{\smalltx}[1]{\,\,\mbox{#1}\,\,}
\newcommand{\seu}{\quad\mbox{or}\quad}
\newcommand{\tx}[1]{\quad\mbox{#1}\quad}

\newcommand{\efrac}[2]{\,\boxed{\dfrac{#1}{#2}}\,}
\newcommand{\amn}[2]{\left.\begin{cases}#1 \\ #2 \end{cases}\!\!\!\!\!\!\!\right)}
\newcommand{\nml}[1]{\multicolumn{1}{r}{#1}}
\let\intn\int
\renewcommand{\int}{\intn\!}
\newcommand{\runon}{\vspace{-\baselineskip}}
\newcounter{eulercounter}

\newenvironment{Definition}[1][]{\refstepcounter{eulercounter}\par\medskip
  \begin{center}\sc\Large Definition #1\end{center}\medskip \qquad\S\theeulercounter.\,\,\,\,}{}
\newenvironment{Hypothesis}[1][]{\refstepcounter{eulercounter}\par\medskip
  \begin{center}\sc\Large Hypothesis #1\end{center}\medskip \qquad\S\theeulercounter.\,\,\,\,}{}
\newenvironment{Lemma}[1][]{\refstepcounter{eulercounter}\par\medskip
  \begin{center}\sc\Large Lemma #1\end{center}\medskip \qquad\S\theeulercounter.\,\,\,\,}{}
\newenvironment{GeneralProblem}[1][]{\refstepcounter{eulercounter}\par\medskip
  \begin{center}\sc\Large General Problem #1\end{center}\medskip \qquad\S\theeulercounter.\,\,\,\,}{}
\newenvironment{ParticularProblem}[1][]{\refstepcounter{eulercounter}\par\medskip
  \begin{center}\sc\Large Particular Problem #1\end{center}\medskip \qquad\S\theeulercounter.\,\,\,\,}{}
\newenvironment{Problem}[1][]{\refstepcounter{eulercounter}\par\medskip
  \begin{center}\sc\Large Problem #1\end{center}\medskip \qquad\S\theeulercounter.\,\,\,\,}{}
\newenvironment{Rule}[1][]{\refstepcounter{eulercounter}\par\medskip
  \begin{center}\sc\Large Rule #1\end{center}\medskip \qquad\S\theeulercounter.\,\,\,\,}{}

\newenvironment{Scholium}[1][]{\refstepcounter{eulercounter}\par\medskip
  \begin{center}\sc\large Scholium #1\end{center}\medskip \qquad\S\theeulercounter.\,\,\,\,}{}

\newenvironment{SolutionUnnum}[1][]{\par\medskip
  \begin{center}\sc\large Solution #1\end{center}\medskip\par}{}
\newenvironment{Corollary}[1][]{\refstepcounter{eulercounter}\par\medskip
  \begin{center}\sc\large Corollary #1\end{center}\medskip \qquad\S\theeulercounter.\,\,\,\,}{}

\newenvironment{Case}[1][]{\refstepcounter{eulercounter}\par\medskip
  \begin{center}\sc\large Case #1\end{center}\medskip \qquad\S\theeulercounter.\,\,\,\,}{}
\newenvironment{RemainingCases}[1][]{\refstepcounter{eulercounter}\par\medskip
  \begin{center}\sc\large Remaining Cases #1\end{center}\medskip \qquad\S\theeulercounter.\,\,\,\,}{}

\newenvironment{Example}[1][]{\refstepcounter{eulercounter}\par\medskip
  \begin{center}\sc\large Example #1\end{center}\medskip \qquad\S\theeulercounter.\,\,\,\,}{}

\newenvironment{FurtherExamples}[1][]{\refstepcounter{eulercounter}\par\medskip
  \begin{center}\sc\large Further Examples #1\end{center}\medskip \qquad\S\theeulercounter.\,\,\,\,}{}

\newenvironment{Catalogue}[1][]{\par\medskip
  \begin{center}\sc\large Catalogue of Amicable Numbers\end{center}\medskip}{}

\begin{document}
\VerbatimFootnotes

\begin{center}
  {\Huge\sc
  On Amicable Numbers}\\

  \medskip
  \normalsize
  
  An English translation of
  \Large\medskip
  
  {\sc De Numeris Amicabilibus}\footnote{Original text
    \cite{euler152} (E152). Also available online at
    \verb|https://scholarlycommons.pacific.edu/euler/| Footnotes
    are comments by the translator (\JE) or adopted from the
    1915 Opera Omnia edition edited by Ferdinand Rudio
    (\FR).}\textsuperscript{,}\footnote{Compare with
    \cite{euler100} and \cite{euler798}. See also
    {\cite[p. 100]{krafft1749}}, and also the letters which
    Krafft wrote to Euler around 1746. \FR}
  
  \medskip

  by Leonhard Euler\\

  \bigskip

  \small
  
  Translated by Jonathan David Evans, School of Mathematical Sciences,\\
  Lancaster University \verb|j.d.evans@lancaster.ac.uk|\\

\end{center}

\paragraph{Translator's note.} This (E152) is the most
substantial of Euler's three papers entitled ``De Numeris
Amicabilibus'' (E100, E152, E798), in which he expounds at great
length the {\em ad hoc} methods he has developed to search for
pairs of amicable numbers. The concept of amicable numbers had
been known for at least two thousand years, and had been
intensively studied by the Arabic school of mathematics after
the work of Thabit Ibn Qurra in the ninth century AD, and then
by Fermat, Descartes and others in the seventeenth century.

Despite this, there were only three pairs known before Euler. In
an earlier paper he found \(26\) further pairs (and one
incorrect pair), and in this 1750 paper he expands the list of
known amicable pairs to \(62\) (and some more incorrect
pairs). According to Dickson's {\em History of the Theory of
  Numbers} {\cite[Chapter I]{DicksonHistory1}}, it wasn't until
1866 that another pair was discovered (by Paganini
\cite{Paganini}), and then until 1911 before three more were
found (by Dickson himself \cite{Dickson}).

Originally written in Latin, this paper has been translated into
Czech, French and German \cite{euler152}, and detailed summaries
can be found elsewhere, including a beautiful summary of Problem
1, Rule 1, Case 1 by Sandifer {\cite[Chapter 9/November
  2005]{Sandifer}} and a staggering four page summary of the
whole paper in Dickson {\cite[pp.42--46]{DicksonHistory1}}, who
customarily dedicates at most a paragraph to any given
paper. Moreover, the Latin original is quite readable, as
Euler's notation is completely modern. Nonetheless, I felt it
not out of place to attempt a translation into English.

There are impressively few typographical errors in the original,
and still fewer in Rudio's 1915 {\em Opera Omnia} edition. Like
Rudio, I have chosen to correct any such errors I have found,
but I have given footnotes to point out where they would have
been. I have also included more amplificatory comments in
footnotes where I found the reasoning mildly non-obvious or to
attempt to reconstruct some of the historical context. I have
also either preserved Rudio's footnotes or updated them. I
include two appendices: one which fills in some mathematical
detail about \S 81, and one which includes the Sage code I used
to check the tables of factorised divisor sums.

I would like to thank the referees for their detailed and thoughtful comments which have improved the translation greatly.

\newpage

\begin{Definition}
  {\em Two numbers are said to be amicable if they have the
    property that\footnote{{\em ita sint comparati} -- literally
      ``are so prepared''. I have chosen to translate this
      construct as ``have the property'', both here and
      elsewhere. \JE} the sum of the aliquot parts\footnote{{\em
        Aliquot parts} means the divisors of a number excluding
      the number itself. \JE}  of one is equal to the other
    number, and, in turn, the sum of the aliquot parts of the
    other equals the first number.}

  \medskip
  
  Thus the numbers \(220\) and \(284\) are amicable; indeed, the
  aliquot parts of the first, \(220\), taken together:
  \[1+2+4+5+10+11+20+22+44+55+110\] make \(284\), and the
  aliquot parts of \(284\): \[1+2+4+71+142\] make the first
  number \(220\).
\end{Definition}

\begin{Scholium} Stifel\footnote{See Stifel (1487--1567)
    {\cite[Folio 10]{stifel1544}}. \FR}, who
  first\footnote{Descartes (1596--1650) {\cite[p. 93--94 (Lettre
      CXIX de Descartes \`{a} Mersenne 31 mars
      1638)]{descartes}} and van Schooten (1615--1660)
    {\cite[Liber V, Sectio IX, p. 419--426.]{VanSchooten}}
    published these three pairs of amicable numbers:
    \(220 = 2^2\cdot 5 \cdot 11\) and \(284 = 2^2\cdot 71\),
    \(17296 = 2^4\cdot 23\cdot 47\) and \(18416=2^2\cdot 1151\),
    \(9363584 = 2^7\cdot 191\cdot 383\) and
    \(9437056 = 2^7\cdot 73727\). Of these, the first pair was
    already known to Pythagoras {\cite[p.35]{pistelli}}, the
    second Fermat (1607--1665) had shared with his friend
    Mersenne (1588--1648) and other mathematicians around 1636
    {\cite[p.136]{fermat1}}, {\cite[p. 20, 21, 71]{fermat2}},
    {\cite[p. 65, 66, 67]{fermat3}}, the third was communicated
    by Descartes around 1638 to his friend Mersenne in the
    letter commended above. \FR}\textsuperscript{,}\footnote{In
    fact, amicable numbers were heavily studied from a
    theoretical perspective by Arab mathematicians such as
    Thabit Ibn Qurra in the late ninth century BCE, who first
    discovered the rule here attributed to Descartes. For an
    overview of the Arabic work on amicable numbers from this
    period, see {\cite[Chapter I]{DicksonHistory1}}, and for a
    more detailed discussion, see Rashed {\cite[Chapter
      4]{Rashed}} \JE}  made mention of this kind of number,
  having noticed these two numbers \(220\) and \(284\) by
  chance, seems to have been led to this speculation; indeed he
  judges analysis to be unsuitable\footnote{See van Schooten,
    \cite[Liber V, Sectio IX, p. 419]{VanSchooten}, according to
    whom, ``There are those who think more arithmetical
    operations are to be found which are not subject to algebra,
    amongst them is the by-no-means-unknown arithmetician
    Michael Stifel.'' He points to a quotation {\cite[Folio
      486--7]{Rudolff}} from Stifel's commentary to the 1554
    edition of Rudolff's "Die Coss" to support this. Indeed,
    Stifel makes some cryptic remarks about ``computations which
    are not subject to algebra'' and then gives the example of
    finding amicable pairs (of which he knows only \(220\) and
    \(284\)). \JE} as a means by which more such pairs of
  numbers are found. However, Descartes tried to adapt analysis
  to this end, and discovered a rule which produced three pairs
  of such numbers; nor was van Schooten, who seems to have
  exerted himself greatly in this investigation, able to extract
  more. Since those times, hardly any Geometers are found to
  have devoted further effort working out this question. Since
  it is also without doubt in this regard that analysis would
  lead to a non-trivial development, if a method were discovered
  which allowed many more such pairs of numbers to be found, I
  judge it would be by no means out of the question if I were to
  relate methods which I have happened upon, with this end in
  view. To that end, it is necessary to introduce the following.
\end{Scholium}

\begin{Hypothesis}
  {\em If \(n\) denotes an arbitrary positive integer, which
    will always be understood hereafter, I will indicate by the
    symbol \(\int n\) the sum of all its divisors, and similarly
    the character \(\int\) prefixed to any number denotes the
    sum of all divisors of this number; so we get\footnote{I have
      chosen to translate {\em erit} (literally ``it will be'')
      into ``we get'' to reflect current idiom. I think this is
      not out of place, since Euler does use the first person
      plural elsewhere, e.g. {\em nanciscimur} ``we obtain'',
      {\em habebimus} ``we have''. \JE}
    \(\int 6=1+2+3+6=12\).}
\end{Hypothesis}

\begin{Corollary}[1]
  Since any number is considered one of its own divisors, but is
  not counted among its aliquot parts, it is clear that the sum
  of the aliquot parts\footnote{I will henceforth translate this
    as {\em aliquot sum}. \JE} of a number \(n\) is expressed as
  \(\int n-n\).
\end{Corollary}

\begin{Corollary}[2]
  Since a prime number has no other divisors than one and
  itself, if \(n\) is a prime number then we get \(\int
  n=1+n\). However, in the case \(n=1\) we get \(\int 1=1\),
  which shows that it would not be correct to count \(1\) amongst the
  prime numbers.
\end{Corollary}

\begin{Lemma}[1]
  If \(m\) and \(n\) are numbers which are relatively prime,
  that is they have no common divisors other than one, then
  \(\int mn=\int m\cdot\int n\); said another way, the sum of
  divisors of the product \(mn\) is equal to the product of the
  sums of divisors of each of the numbers \(m\) and \(n\).

  Indeed, the product \(mn\) has firstly one divisor for each
  factor \(m\) and \(n\), and moreover is divisible by the
  product of each divisor of \(m\) with each divisor of
  \(n\). Verily, all these divisors of \(mn\) appear in
  combination if \(\int m\) is multiplied by \(\int n\).
\end{Lemma}

\begin{Corollary}[1]
  If the numbers \(m\) and \(n\) are each prime, and so
  \(\int m=1+m\) and \(\int n=1+n\), the divisor sum of the
  product will be \[\int mn=(1+m)(1+n)=1+m+n+mn.\] If, moreover,
  \(p\) is a prime number different from \(m\) and \(n\), we get
  \[\int mnp=\int mn\cdot\int p=\int m\cdot\int n\cdot\int
  p=(1+m)(1+n)(1+p).\] And hence the divisor sum of all numbers
  which are products of distinct primes will easily be assigned.
\end{Corollary}

\begin{Corollary}[2]
  If \(m\), \(n\), and \(p\) are not themselves prime numbers,
  but are nonetheless such that they have no common factors
  other than \(1\), then \(mn\) and \(p\) will be relatively
  prime, and therefore
  \(\int mnp=\int mn\cdot\int p\). But since
  \(\int mn=\int m\cdot\int n\), we get
  \(\int mnp=\int m\cdot \int n\cdot\int p\).
\end{Corollary}

\begin{Scholium}
  Unless the factors \(m\), \(n\), \(p\) are relatively prime,
  the divisor sum of the product, as stated in the lemma, is
  incorrect. Indeed, since according to the lemma each divisor
  of the factors \(m\), \(n\), \(p\) is considered as a divisor
  of the product, if they were to have a divisor in common, it
  would be counted twice as a divisor of the product; however,
  when calculating the divisor sum of an arbitrary number, no
  divisor should be counted twice. Hence, if \(m\) and \(n\) are
  prime numbers and \(m=n\), we do not get
  \(\int nn=\int n\cdot \int n=(1+n)^2=1+2n+nn\), but will have
  \(\int nn=1+n+nn\), as the divisor \(n\) should not be counted
  twice. Therefore, since divisor sums of numbers which are
  products of distinct primes are assigned correctly by this
  lemma, all that is left is to find the rule for equal factors,
  by means of which the divisor sum of a product may be
  determined.
\end{Scholium}

\begin{Lemma}[2]
  If \(n\) is a prime number, we get \[\int n^2=1+n+n^2,\quad
  \int n^3=1+n+n^2+n^3,\quad \int n^4=1+n+n^2+n^3+n^4,\] and
  in general we get
  \(\int n^k=1+n+n^2+\cdots+n^k=\dfrac{n^{k+1}-1}{n-1}\).
\end{Lemma}

\begin{Corollary}[1]
  Since \(\int n=1+n\), we get \(\int n^2=\int n+n^2\), or
  equivalently \(\int n^2=1+n\int n\). In a similar manner, we
  get \(\int n^3=\int n^2+n^3\) or \(\int n^3=1+n\int n^2\);
  moreover, \(\int n^4=\int n^3+n^4\) or\footnote{The original
    has a typo here: \(\int n^4=1+\int n^3\). This is tacitly
    corrected in the Opera Omnia edition and I have also
    corrected it here. \JE} \(\int n^4=1+n\int n^3\), and so
  on. And thus, by knowing the divisor sum of each power
  \(n^k\), the divisor sum of the subsequent power \(n^{k+1}\)
  can easily be assigned, since
  \(\int n^{k+1}=\int n^k+n^{k+1}\), or
  \(\int n^{k+1}=1+n\int n^k\).
\end{Corollary}

\begin{Corollary}[2]
  So that divisor sums may easily be factorised, note that
  \begin{gather*}
    \int n^3 = (1+n)(1+n^2)=(1+n^2)\!\!\int n\\
    \int n^5 =(1+n^2+n^4)\!\!\int n,\,\,
    \int n^7=(1+n^2+n^4+n^6)\!\!\int n=(1+n^4)(1+n^2)\!\!\int n;
  \end{gather*}
  and thus the divisor sums of odd powers always factorise, but
  divisor sums of even powers will sometimes be prime.
\end{Corollary}

\begin{Corollary}[3]
  Hence therefore it will be easy to compile a table, which
  shows the divisor sums not only of prime numbers, but also of
  their powers. Such a Table is seen to be attached here, in
  which divisor sums (in factorised form) are given for all
  prime numbers not bigger than a thousand, as well as their
  powers up to the third, or higher for smaller
  numbers.\footnote{In the first edition and, according to
    Rudio, also in the edition of Fuss \cite{eulerfuss}, there were
    some errors, which were corrected in the 1915 Opera Omnia
    edition. I have given the corrected versions here and also
    corrected an error from the 1915 edition: namely \(7^{10}\)
    is stated as \(329554457\) there, instead of its prime
    factorisation \(1123\cdot 293459\). The other, earlier,
    errors were that: the powers of \(79\) were omitted, and the
    following entries were given: \newline\scriptsize
    \vspace{-0.3cm}
    \begin{center}
      \begin{tabular}{p{0.7cm}|p{2.5cm}||p{0.7cm}|p{3cm}} \(5^5\) & \(2\cdot 3^3\cdot 7\cdot 31\)  &  \(523^3\) & \(2^3\cdot 5\cdot 7\cdot 131\cdot 1609\) \\
        \(37^3\) & \(2^2\cdot 5\cdot 2603\)  & \(563^2\) & \(2^3\cdot 35\cdot 29\cdot 47\cdot 1093\)\\
        \(41^3\) & \(2^2\cdot 3\cdot 7\cdot 29_2\)& \(571^3\) & \(2^3\cdot 11\cdot 13\cdot 163041\) \\
        \(149^3\) & \(2^2\cdot 3\cdot 5^2\cdot 11\cdot 101\) & \(613^2\) & \(3\cdot 125461\) \\
        \(173^2\) & \(67\cdot 449\) & \(769^3\) & \(2^2\cdot 5\cdot 7\cdot 11\cdot 71\cdot 17393\)\\
        \(283^3\) & \(2^2\cdot 5\cdot 71\cdot 8009\) & \(811\) & \(2\cdot 7\cdot 29\) \\
        \(461^3\) & \(2^2\cdot 3\cdot 7\cdot 11106261\)&\(827\) & \(2^2\cdot 3^3\cdot 23\)\end{tabular}\end{center}

    \footnotesize Instead of copying these tables out by hand, I
    used some Sage code to produce the LaTeX source. See
    Appendix \ref{app:code} for the code. \JE}
\end{Corollary}

\bgroup
\let\clearpage\relax
\def\arraystretch{1.25}
\input{bigtab}
\egroup

\begin{Scholium}
  The use of this table is primarily in resolving questions
  revolving around divisors and aliquot parts.  Indeed with its
  help, the divisor sum of any given number can be an easy
  matter to find; if that same given number were taken away from
  the result, what remains is its aliquot sum. From this it is
  immediately clear that, with the help of this table, one can
  easily check whether the amicable numbers, which I am about to
  relate, are correct or not. In the following lemma, I will
  explain how one can know the divisor sum of an arbitrary
  number by means of this table.
\end{Scholium}

\newpage

\begin{Lemma}[3]
  {\em For any given number whatsoever, its divisor sum is
    produced in the following manner.}

  \medskip
  
  Since every number is either prime or a product of primes, the
  given number is resolved into its prime factors, and any
  amongst them which are equal are grouped together. In this
  manner the given number will always be put into the form
  \(m^\alpha\cdot n^\beta\cdot p^\gamma\cdot
  q^\delta\cdot\,\mbox{etc.}\) where \(m\), \(n\), \(p\), \(q\),
  etc. are prime numbers. Therefore, calling the given number
  \(=N\), since
  \(N=m^\alpha\cdot n^\beta\cdot p^\gamma\cdot
  q^\delta\cdot\,\mbox{etc.}\) and the factors \(m^\alpha\),
  \(n^\beta\), \(p^\gamma\), \(q^\delta\), etc. are relatively
  prime, we get
  \(\int N = \int m^\alpha\cdot \int n^\beta\cdot\int
  p^\gamma\cdot \int q^\delta\cdot\,\mbox{etc.}\) and the values
  of \(\int m^\alpha\), \(\int n^\beta\), \(\int p^\gamma\),
  \(\int q^\delta\) etc. will be evident from the annexed table.

  \medskip
  
  {\sc 1. Example.} \emph{Let the given number be \(N=360\).}

  \medskip
  
  By resolving this number into its prime factors, we get
  \(N = 2^3\cdot 3^2\cdot 5\), and so
  \(\int 360 = \int 2^3\cdot \int 3^2\cdot\int 5 =
    3\cdot 5 \cdot 13\cdot 2 \cdot 3\), because
  \(\int 2^3=3\cdot 5\), \(\int 3^2=13\), \(\int 5=2\cdot
  3\).

  Whence by arranging these factors we get
  \(\int 360 = 2\cdot 3^2\cdot 5 \cdot 13 = 1170\).

  \medskip

  {\sc 2. Example.} {\em It is checked whether or not the
    numbers \(2620\) and \(2924\) are amicable.}

  \medskip
  
  Since we have \(2620=2^2\cdot 5 \cdot 131\) and
  \(2924=2^2\cdot 17\cdot 43\), the calculation will be organised as follows.

  \begin{center}
    \bgroup
    \def\arraystretch{1.2}
    \begin{tabular}{p{5.2cm}|c|c|}
      \qquad Given numbers & \(2620\) & \(2924\) \\
      expressed in factors & \(2^2\cdot 5 \cdot 131\) & \(2^2\cdot 17\cdot 43\) \\
      divisor sum & \(7\cdot 6 \cdot 132\) & \(7\cdot 18\cdot 44\)\\
      \centering or & \(5544\) & \(5544\)\\
      Aliquot sum & \(2924\) & \(2620\)
    \end{tabular}
    \egroup
  \end{center}

  Therefore, since the aliquot sums are equal to the alternate
  numbers, this shows the given numbers to be amicable.
\end{Lemma}

\begin{Scholium}
  Therefore, with this said in advance about finding the
  divisors of all numbers, I will proceed to the problem of
  looking for amicable numbers, and I will examine in what ratio
  such numbers must be with their divisor sum, so that in turn
  they can more easily be discovered following the rules related
  below.
\end{Scholium}

\newpage

\begin{GeneralProblem}
  {\em To find amicable numbers, that is two numbers of the
    character that each is equal to the aliquot sum of the
    other.}
\end{GeneralProblem}

\begin{SolutionUnnum}
  Let \(m\) and \(n\) be two such amicable numbers, and by
  hypothesis \(\int m\) and \(\int n\) their divisor sums. The
  aliquot sum of \(m\) will be \(=\int m-m\) and the aliquot sum
  of \(n\) will be \(=\int n-n\). Hence by the nature of
  amicable numbers, these two equations will follow:
  \begin{gather*}
    \int m-m=n\quad \mbox{and}\quad \int n-n=m,\\
    \mbox{or}\qquad\int m=\int n=m+n.
  \end{gather*}
  \par Therefore the amicable numbers \(m\)
  and \(n\) must in the first place have the same divisor sum,
  and moreover it must be that this common divisor sum is equal
  to the sum \(m+n\) of these same numbers.
\end{SolutionUnnum}

\begin{Corollary}[1]
  The problem is therefore reduced to that of finding two
  numbers which have the same divisor sum, and this is equal to
  the sum of those same numbers.
\end{Corollary}

\begin{Corollary}[2]
  Indeed the motivation for the problem demands that the two
  numbers we seek are not equal to one another. If however they
  were desired to be equal, so that \(m=n\), then we get
  \(\int n=2n\) and \(\int n-n=n\); that is to say, the aliquot
  sum of this repeated number is equal to itself, which is the
  property of a perfect number. Therefore any perfect number
  repeated counts as a pair of amicable numbers.
\end{Corollary}

\begin{Corollary}[3]
  If however, the amicable numbers \(m\) and \(n\) are unequal,
  as the nature of the question postulates, it is clear that one
  is abundant and the other deficient; that is to say, the
  aliquot sum of one is bigger than itself, and of the other is
  smaller than itself.
\end{Corollary}

\newpage

\begin{Scholium}
  Indeed, from this general property we obtain very little
  assistance in finding amicable numbers, because this kind of
  analysis, which can be derived by means of the equation
  \(\int m=\int n=m+n\) is, even now, heavily
  underdeveloped. Because of this deficiency we are forced to
  contemplate more specific formulas, from which rules of a
  special nature may be derived for the discovery of amicable
  numbers; to which class also belongs the rule of Descartes
  related by van Schooten.\footnote{See \cite{euler100}. {\FR}
    See also {\cite[Liber V, Sectio IX]{VanSchooten}}. \JE} And
  indeed first, even if it is not settled whether or not there
  exist amicable numbers which are relatively prime, I will
  restrict general formulas so that the amicable numbers possess
  common factors.
\end{Scholium}

\begin{ParticularProblem}
  To find the kind of amicable numbers which have common factors.
\end{ParticularProblem}

\begin{SolutionUnnum}
  Let \(a\) be the common factor of the amicable numbers, of
  which we put one \(=am\) and the other \(=an\); let \(m\) and
  \(a\), and also \(n\) and \(a\), be relatively prime, so that
  each of their divisor sums may be found using the rule
  given. First, therefore, since each divisor sum must be the
  same, we get \(\int a\cdot \int m=\int a\cdot n\), and
  so \(\int m=\int n\). Then indeed it is necessary that
  \(\int a\cdot \int m\) or \(\int a\cdot \int n\) are equal to
  the sum of the same numbers \(am+an\), whence we have
  \[\frac{a}{\int a}=\frac{\int m}{m+n}=\frac{\int n}{m+n}.\]
  Therefore, by supposing \(am\) and \(an\) to be amicable
  numbers, it must hold first that \(\int m=\int n\), and then
  indeed we need that \(a(m+n)=\int a\cdot m\).
\end{SolutionUnnum}

\begin{Corollary}[1]
  Therefore, if \(m\) and \(n\) are taken to be numbers such
  that \(\int m=\int n\), then a number \(a\) must be sought
  such that \(\dfrac{a}{\int a}=\dfrac{\int m}{m+n}\);
  equivalently this same number \(a\) must be sought from the
  ratio it must hold with its divisor sum.
\end{Corollary}

\begin{Corollary}[2]
  If the common factor \(a\) is given, the question is reduced
  to finding numbers \(m\) and \(n\) which are,
  case-by-case\footnote{My interpretation of {\em prouti} here,
    based on what Euler actually does in the sequel. \JE},
  assumed to be either prime or composed of two or more primes;
  since the divisor sums can then be computed explicitly,
  special rules for their discovery will be formulated.
\end{Corollary}

\begin{Corollary}[3]
  However, it is seen at once that both numbers \(m\) and \(n\)
  cannot be prime\footnote{i.e. at the same time, because then
    \(m+1=\int m=\int n=n+1\) so \(m=n\). \JE}; wherefore the
  simplest case is if one is assumed prime and the other a
  product of two prime numbers. One can then take each of the
  numbers to be the product of two or more primes, and whence
  countless special rules for finding amicable numbers can be
  derived.
\end{Corollary}

\begin{Scholium}
  Therefore the many forms of amicable numbers which can be
  generated from this can be represented in the following
  manner. Let \(a\) be the common factor of each, and \(p\),
  \(q\), \(r\), \(s\), etc. primes of which none are divisors of
  the common factor \(a\); and the forms of amicable numbers
  will be:

  \begin{center}
    \begin{tabular}{p{2.5cm}cccp{2cm}}
      First form & - & - & - & \(\begin{cases}apq\\ ar\end{cases}\)\\
      Second form & - & - & - & \(\begin{cases}apq\\ ars\end{cases}\)\\
      Third form & - & - & - & \(\begin{cases}apqr\\ as\end{cases}\)\\
      Fourth form & - & - & - & \(\begin{cases}apqr\\ ast\end{cases}\)\\
      Fifth form & - & - & - & \(\begin{cases}apqr\\ astu\end{cases}\)\\
      \qquad etc. & & &
    \end{tabular}
  \end{center}
  
  \par Although the number of these forms can increase to
  infinity, one can by no means conclude that all amicable
  numbers are in one of these forms. Indeed, firstly, whilst
  these letters \(p\), \(q\), \(r\), \(s\), \(t\), etc. signify
  distinct primes, it is unlikely that no amicable numbers can
  be given in which no powers of the same prime
  occur\footnote{Except as part of the common factor
    \(a\). \JE}. Next, equally, it is not certain whether there
  exist amicable numbers which have either no common factor, or
  for which this factor does not appear in the same way: that is
  if there exist amicable numbers of the form \(m^\alpha P\) and
  \(m^\beta Q\), in which the exponents \(\alpha\) and \(\beta\)
  are different; on account of this such forms are not amongst
  those listed above, even if \(P\) and \(Q\) were products of
  pure\footnote{{\em meris} = pure (ablative plural). I believe Euler means ``not raised to any power''. \JE} prime numbers, different from one another\footnote{But see
    pairs LX and LXI in the Calatogue at the end of the
    paper. {\FR}}. From this, the question
  of amicable numbers is seen to be most extensive and, because
  of this, to be difficult to the extent that a complete
  solution should scarcely be expected. Therefore, for my part,
  I will devote myself merely to particular solutions, and
  reveal various methods by means of which it was permitted to
  me to elicit more amicable numbers from the traditional
  formulas.  Moreover, each form supplied me with two methods,
  according to whether the common factor \(a\) is assumed given
  or is sought; and I will explain these methods in the
  following problems.
\end{Scholium}

\begin{Problem}[1]
  {\em To find amicable numbers of the first form \(apq\) and
    \(ar\) if the common factor \(a\) is given.}
\end{Problem}

\begin{SolutionUnnum}
  Since \(p\), \(q\), and \(r\) should be prime numbers, and
  \(\int r=\int p\cdot \int q\) or \[r+1=(p+1)(q+1),\] putting
  \(p+1=x\) and \(q+1=y\) makes \(r=xy-1\). And so \(x\) and
  \(y\) must be numbers such that not only \(x-1\) and \(y-1\)
  but also \(xy-1\) are prime numbers. Then, in order that
  \(a(x-1)(y-1)\) and \(a(xy-1)\) are amicable numbers, it is
  necessary that their sum \(a(2xy-x-y)\) be equal to the
  divisor sum \(xy\int a\) of either\footnote{i.e. the divisor sum of either of the
    amicable numbers (they have the same divisor sum). \JE};
  whence we obtain this equation
  \[xy\int a=2axy-ax-ay\tx{or}y=\frac{ax}{\left(2a-\int
        a\right)x-a}.\] For sake of brevity, let
  \(\dfrac{a}{2a-\int a}=\dfrac{b}{c}\), where \(\dfrac{b}{c}\) is the
  value of the fraction \(\dfrac{a}{2a-\int a}\) written in
  lowest terms; we get
  \[y=\frac{bx}{cx-b}\tx{or}cy =
    \frac{bcx}{cx-b}=b+\frac{bb}{cx-b},\] whence we
  have \[(cx-b)(cy-b)=bb.\] Therefore since \(cx-b\) and
  \(cy-b\) are factors of \(bb\), the known square number \(bb\)
  must be resolved into two factors, each of which, when
  increased by \(b\) becomes divisible by \(c\), and the
  quotients \(x\) and \(y\) arising from this have the property
  that \(x-1\), \(y-1\), and \(xy-1\) result in prime numbers.
  Whenever this condition can be satisfied, which will
  immediately be discerned for any assumed value of \(a\), then
  amicable numbers will be obtained, which will be
  \(a(x-1)(y-1)\) and \(a(xy-1)\).  Q.E.J.\footnote{This seems
    to be a form of {\em Q.E.I. (Quod Erat Inveniendum)}: which was to be
    found. \JE}
\end{SolutionUnnum}

\newpage

\begin{Corollary}
  Therefore, as one number or another is taken for \(a\), which
  determines the values of \(b\) and \(c\), particular rules
  emerge, by means of which amicable numbers can easily be
  extracted, if any exist in that form.
\end{Corollary}

\begin{Rule}[1]
  Let the common factor \(a\) be any power of two, say \(2^n\);
  we get \(\int a=2^{n+1}-1\), and so \(2a-\int a=1\), whence we
  get \(\dfrac{a}{2a-\int a}=2^n\), and therefore
  \(b=2^n\) and \(c=1\). Hence it arises that
  \[(x-2^n)(y-2^n)=2^{2n}.\]

  \par Since \(2^{2n}\) has no other factors than powers of two,
  we get

  \bgroup
  \setlength{\tabcolsep}{1pt}
  \begin{center}
    \begin{tabular}{p{1cm}p{0.5cm}p{1.2cm}p{1.5cm}p{0.2cm}p{0.5cm}p{1.6cm}}
      \(x-2^n\)&\centering\(=\)&\(2^{n+k}\)&\centering\multirow{2}{1.5cm}{\quad or\quad}&\(x\)&\centering\(=\)&\(2^{n+k}+2^n\)\\
      \(y-2^n\)&\centering\(=\)&\(2^{n-k}\)&&\(y\)&\centering\(=\)&\(2^{n-k}+2^n.\)
    \end{tabular}
  \end{center}
  \egroup
  
  \par Wherefore it is to be considered whether a value of \(k\)
  exists which makes the following three numbers prime
  \begin{align*}
    x-1&=2^{n+k}+2^n-1,\\
    y-1&=2^{n-k}+2^n-1,\\
    xy-1&=2^{2n+1}+2^{2n+k}+2^{2n-k}-1.
  \end{align*}
  If this succeeds, the amicable numbers will be:
  \begin{gather*}
    2^n\left(2^{n+k}+2^n-1\right)\left(2^{n-k}+2^n-1\right),\\
    2^n\left(2^{2n+1}+2^{2n+k}+2^{2n-k}-1\right).
  \end{gather*}

  \par Alternatively, let \(n-k=m\), or \(n=m+k\), making
  \begin{alignat*}{3}
    x-1&=2^m\left(2^{2k}+2^k\right)-1&\,\,=q,\\
    y-1&=2^m\left(1+2^k\right)-1&\,\,=p,\\
    xy-1&=2^{2m}\left(2^{2k+1}+2^{3k}+2^k\right)-1&\,\,=r,
  \end{alignat*}
  which numbers, whenever they become prime, will furnish us
  with amicable numbers.
\end{Rule}

\newpage

\begin{Case}[1]
  Let \(k=1\); amicable numbers will be obtained whenever the
  following three numbers become prime:
  \[3\cdot 2^m-1,\quad 6\cdot 2^m-1\quad\et\quad 18\cdot
    2^{2m}-1.\] Indeed by then putting
  \[p=3\cdot 2^m - 1,\quad q=6\cdot 2^m - 1\et
    r=18\cdot 2^{2m}-1\] the amicable numbers will be
  \(2^{m+1}pq\) and \(2^{m+1}r\) because \(n=m+k=m+1\). And this is
  the rule of Descartes related by van Schooten.\footnote{See
    \cite{euler100}. \FR See also {\cite[Liber V, Sectio
      IX]{VanSchooten}}. \JE}
\end{Case}

\begin{Example}[1]
  Let \(m=1\) and we get
  \begin{alignat*}{3}
    p&=&\,\,3\cdot 2-1 \,&=&\,\, 5&\quad\mbox{prime number},\\
    q&=&\,\,6\cdot 2-1 \,&=&\,\, 11&\quad\mbox{prime number},\\
    r&=&\,\,18\cdot 2-1 \,\,&=&\,\, 71&\quad\mbox{prime number}.
  \end{alignat*}
  Hence, therefore, the amicable numbers arising are
  \[2^2\cdot 5\cdot 11\et 2^2\cdot 71\] or \(220\) and \(284\)
  which are the smallest that can be produced.
\end{Example}

\begin{Example}[2]
  Let \(m=2\) and we get \(2^m=4\) and \(2^{2m}=16\) and also
  \begin{alignat*}{3}
    p&=&\,\,3\cdot 4-1 \,&=&\,\, 11&\quad\mbox{prime number},\\
    q&=&\,\,6\cdot 4-1 \,&=&\,\, 23&\quad\mbox{prime number},\\
    r&=&\,\,18\cdot 16-1 \,&=&\,\, 287&\quad\mbox{non-prime number};
  \end{alignat*}
  and hence we find no amicable numbers.
\end{Example}

\begin{Example}[3]
  Let \(m=3\) and we get \(2^m=8\) and \(2^{2m}=64\) and also
  \begin{alignat*}{3}
    p&=&\,\,3\cdot 8-1\, &=&\,\, 23&\quad\mbox{prime},\\
    q&=&\,\,6\cdot 8-1\, &=&\,\, 47&\quad\mbox{prime},\\
    r&=&\,\,18\cdot 64-1\,&=&\,\, 1151&\quad\mbox{prime}.
  \end{alignat*}
  Therefore the amicable numbers will be
  \[2^4\cdot 23\cdot 47\et 2^4\cdot 1151\qquad
    \mbox{or}\qquad 17296\quad\mbox{and}\quad 18416.\]
\end{Example}

\begin{FurtherExamples}
  These examples, like the following in which the exponent \(m\)
  is assigned bigger values, can be displayed more conveniently
  in one view thus

  \bgroup
  \def\arraystretch{1.1}
  \begin{center}
    \begin{tabular}{>{\raggedleft\arraybackslash}p{1.9cm}r|
      >{\raggedleft\arraybackslash}p{0.6cm}@{}p{0.01cm}|
      >{\raggedleft\arraybackslash}p{0.7cm}@{}p{0.01cm}|
      >{\raggedleft\arraybackslash}p{0.8cm}@{}p{0.01cm}|
      >{\raggedleft\arraybackslash}p{1cm}@{}p{0.01cm}|
      >{\raggedleft\arraybackslash}p{1cm}@{}p{0.01cm}|
      >{\raggedleft\arraybackslash}p{1.2cm}@{}p{0.01cm}|
      >{\raggedleft\arraybackslash}p{1.4cm}@{}p{0.01cm}}
      \hline\hline
      Let \(m=\)&\(1\) & \(2\) && \(3\) && \(4\) && \(5\) && \(6\) && \(7\) && \(8\)&\\
      \hline\hline
      we get \(p=\)&\(5\)  & \(11\) && \(23\) && \(47\) && \(95\)&\(^*\) & \(191\) && \(383\) && \(767\) &\(^*\)\\
       \(q=\)&\(11\)      & \(23\) && \(47\) && \(95\) &\(^*\) & \(191\) && \(383\) && \(767\)&\(^*\) & \(1535\)&\(^*\)\\
      \(r=\)&\(71\)      & \(287\) &\(^*\) & \(1151\) && \(4607\)&\(^*\) & \(18431\)&\(^*\) & \(73727\) && \(294911\) && \(1179647\)&\(^\dagger\)\\
    \end{tabular}
  \end{center}
  \egroup
  
  \noindent where non-prime numbers are denoted with an
  asterisk\footnote{As Rudio notes, some of the asterisks are
    missing on the tables in the original. In most cases this is
    likely to be because the case has already been dismissed and
    the remaining number is too big to factorise. Instead of
    adding the asterisks back in, we put \(\dagger\)s in for the
    missing asterisks. \JE \label{fn:daggers}}; whence only three amicable
  numbers\footnote{These numbers constitute the same three pairs
    that were known before Euler. \FR} are obtained in this way,
  namely:
  \[\mbox{I.}\amn{2^2\cdot 5\cdot 11}{2^2\cdot
      71}\qquad \mbox{II.}\amn{2^4\cdot 23\cdot 47}{2^4\cdot
      1151}\qquad \mbox{III.}\amn{2^7\cdot 191\cdot
      383}{2^7\cdot 73727.}\] But further we may not progress,
  because the values of \(r\) become exceedingly large, so that
  one cannot tell whether or not they are prime. For indeed the
  tables of prime numbers constructed to date\footnote{See
    Glaisher {\cite[p.1, especially
      p.34--40]{Glaisher}}. Amongst the books enumerated in this
    report, the following three, being published before around
    1750 are worth mentioning here:\\-- 1) Rahn \cite{Rahn}:
    this contains tables of divisors of odd numbers as far as
    \(24000\).\\ -- 2) The translation \cite{Brancker} of this
    notable book (see \cite[p.113]{Wertheim}). For in the
    translation, J. Pell has continued the table of divisors
    constructed by Rahn as far as 10000. J. Wallis gave several
    (30) corrections to this Pellian table in the work
    \cite[Additional Treatise IV, p.136]{Wallis}. \\ -- 3)
    Kr\"{u}ger \cite{Kruger}: this contains a table of prime
    numbers as far as \(100999\), which table Euler was without
    doubt in the habit of using. This work is found indeed in
    Euler's index {\em Catalogus librorum meorum}, of which
    G. Enestr\"{o}m made mention in the report
    {\cite[p.197]{Enestrom}}, ``[Sixth notebook]...S.363--402 is
    a directory of Leonhard Euler's library (539 titles)''. \FR}
  hardly extend beyond \(100000\).
\end{FurtherExamples}

\begin{Case}[2]
  Let \(k=2\) and the values of the letters \(p\), \(q\), \(r\)
  which must be prime will be
  \begin{alignat*}{2}
    p&=&\,\,5\cdot 2^m-1,\\
    q&=&\,\,20\cdot 2^m-1,\\
    r&=&\,\,100\cdot 2^{2m}-1.
  \end{alignat*}
  Since the last of these is always divisible by three,
  because \(2^{2m}=3\alpha+1\) and \(r=300\alpha+99\), no new
  amicable numbers are obtained from this.
\end{Case}

\begin{Case}[3]
  Suppose \(k=3\) and so we get
  \begin{alignat*}{2}
    p&=&\,\,9\cdot 2^m-1,\\
    q&=&\,\,72\cdot 2^m-1,\\
    r&=&\,\,648\cdot 2^{2m}-1.
  \end{alignat*}
  Since we see that none of these necessarily
  admit divisors, I will represent the values of \(p\), \(q\),
  \(r\) arising from the simpler value of \(m\) together here:

  \bgroup
  \def\arraystretch{1.1}
  \begin{center}
    \begin{tabular}{>{\raggedleft\arraybackslash}p{0.7cm}
      >{\raggedleft\arraybackslash}p{0.8cm}@{}p{0.01cm}|
      >{\raggedleft\arraybackslash}p{1.1cm}@{}p{0.01cm}|
      >{\raggedleft\arraybackslash}p{1cm}@{}p{0.01cm}|
      >{\raggedleft\arraybackslash}p{1.2cm}@{}p{0.01cm}|
      >{\raggedleft\arraybackslash}p{1.2cm}@{}p{0.01cm}}
      \hline\hline
      \(m=\) & \(1\) && \(2\) && \(3\) && \(4\) && \(5\)& \\
      \hline\hline
      \(p=\) & \(17\) && \(35\)&\(^*\) & \(71\) && \(143\)&\(^*\) & \(287\)&\(^*\)\\
      \(q=\) & \(143\)&\(^*\) & \(287\)&\(^*\) & \(575\)&\(^*\) & \(1151\) && \(2303\)&\(^*\)\\
      \(r=\) & \(2591\) && \(10367\)&\(^*\) & \(41471\)&\(^*\) & \(165887\) && \(663551\)&\(^\dagger\)
    \end{tabular}
  \end{center}
  \egroup
  
  \par Hence therefore, since further progress may not be
  made,\footnote{Because \(r\) is already much bigger than the
    limit of \(100000\) imposed by the tables of primes
    available to Euler. \JE} no amicable numbers are found.
\end{Case}

\begin{Case}[4]
  Suppose \(k=4\); the following three numbers must become prime
  \begin{alignat*}{2}
    p&=&\,\,17\cdot 2^m-1,\\
    q&=&\,\,272\cdot 2^m-1,\\
    r&=&\,\,4624\cdot 2^{2m}-1,
  \end{alignat*}
  where, since \(r\) is always a multiple of three, it is clear
  that no amicable numbers appear from this.
\end{Case}

\begin{Case}[5]
  Suppose \(k=5\); the following three numbers must become prime
  \begin{alignat*}{2}
    p&=&\,\,33\cdot 2^m-1,\\
    q&=&\,\,1056\cdot 2^m-1,\\
    r&=&\,\,34848\cdot 2^{2m}-1,
  \end{alignat*}
  where it is immediately clear that the case \(m=1\) is of no
  use, since it would give \(p=65\). Therefore let \(m=2\)
  making \[p=131,\quad q=4223^*,\quad r=557567;\] where since
  \(q\) is not prime, and bigger values of \(m\) do not submit to
  examination because of insufficient tables of prime numbers,
  no new amicable numbers are extracted in this case once again.
  But indeed for the same reason, \(k\) cannot be allowed to
  take on bigger values.
\end{Case}

\begin{Scholium}
  Since putting powers of two for \(a\) will yield unity as the
  value of \(c\) in the fraction
  \(\dfrac{b}{c}=\dfrac{a}{2a-\int a}\) and hence allow us to
  obtain solutions, I will use other values of \(a\) which also
  give \(c\) the value \(=1\). Especially to be noted
  amongst them are those arising in the form
  \(a=2^n\left(2^{n+1}+e\right)\), where \(2^{n+1}+e\) is a
  prime number; indeed then we get
  \[2a-\int
    a=e+1\et\frac{b}{c}=\frac{2^n\left(2^{n+1}+e\right)}{e+1};\]
  if therefore \(e+1\) is a divisor of the number
  \(2^n\left(2^{n+1}+e\right)\), the value of \(c\) will
  likewise be \(=1\).
\end{Scholium}

\begin{Rule}[2]
  Let the common factor be \(a=2^n\left(2^{n+1}+2^k-1\right)\),
  with \(2^{n+1}+2^k-1\) a prime number; because \(e+1=2^k\),
  the fraction will be\footnote{There is a typo in the original:
    \(\tfrac{c}{b}\). \JE}
  \[\dfrac{b}{c}=\dfrac{2^n\left(2^{n+1}+2^k-1\right)}{2^k} =
    2^{n-k}\left(2^{n+1}+2^k-1\right),\] as long as it is not
  the case that \(k>n\). Therefore by hypothesis we will have
  \[b = 2^{n-k}\left(2^{n+1}+2^k-1\right)\et c=1.\] Therefore
  the square \(bb\) is to be resolved into two factors
  \((x-b)(y-b)\) from which not only the values of the numbers
  \(x-1=p\) and \(y-1=q\), but also \(xy-1=r\) become prime
  numbers. If it is possible to find such an occurrence, we will
  get amicable numbers \(apq\) and \(ar\).  However, it is
  important here for those cases to be rejected in which any of
  the prime numbers \(p\), \(q\), \(r\) turn out to be a divisor
  of \(a\), that is to say, equal to \(2^{n+1}+2^k-1\) (because
  \(a\) is divisible by no other prime numbers).

  Let \(n-k=m\) or \(n=m+k\); we get
  \[a=2^{m+k}\left(2^{m+k+1}+2^k-1\right)\et
    b=2^m\left(2^{m+k+1}+2^k-1\right).\] Now because
  \(2^{m+k+1}+2^k-1\) must be a prime number, put
  \[2^{m+k+1}+2^k-1=f\seu
    f=2^k\left(2^{m+1}+1\right)-1,\] so that\runon
  \[a=2^{m+k}f\et b=2^mf;\] we get\runon \[bb=2^{2m}ff= (x-b)(y-b).\]
  Now because \(f\) is a prime number, the number \(2^{2m}ff\)
  can, in general, be resolved into two factors in two ways.

  Resolving in the first way gives
  \[(x-b)(y-b) = 2^{m-\alpha}f\cdot 2^{m+\alpha}f\]
  and so
  \begin{align*}
    x&=2^{m-\alpha}f+2^mf,& p&=\left(2^{m-\alpha}+2^m\right)f-1, \\
    y&=2^{m+\alpha}f+2^mf,&q&=\left(2^{m+\alpha}+2^m\right)f-1 \\
    &\qquad\qquad\mbox{and}& r&=\left(2^{2m+1}+2^{2m+\alpha}+2^{2m-\alpha}\right)ff-1,
  \end{align*}
  which three numbers \(p\), \(q\), \(r\) must be prime.

  Resolving in the other way gives
  \[(x-b)(y-b) = 2^{m\pm \alpha}\cdot 2^{m\mp\alpha}ff,\]
  whence we get
  \begin{align*}
    x&=2^{m\pm\alpha}+2^mf,& p&=2^{m\pm\alpha}+2^mf-1, \\
    y&=2^{m\mp\alpha}ff+2^mf,& q&=\left(2^{m\mp\alpha}f+2^m\right)f-1\\
     &\qquad\qquad\mbox{and}&r&=\left(2^{2m+1}f+2^{2m\pm\alpha}+2^{2m\mp\alpha}ff\right)f-1,
  \end{align*}
  and whenever prime numbers \(p\), \(q\), \(r\) appear in this
  way, there arise thereby the amicable numbers \(apq\) and
  \(ar\).
\end{Rule}

\begin{Case}[1]
  Let \(k=1\); we get \(a=2^{m+1}\left(2^{m+2}+1\right)\),
  \(b=2^m\left(2^{m+2}+1\right)\), and \(f=2^{m+2}+1\), which
  number must be prime. Therefore, since \((x-b)(y-b)=2^{2m}ff\),
  we get\vspace{-0.8cm}
  \begin{center}
    \begin{tabular}{p{6.1cm}|p{6.7cm}}
      \centering either & \centering or \arraybackslash\\
      \(\!\!\!\!p=\left(2^{m-\alpha}+2^m\right)f-1\), & \(p=2^{m\pm\alpha}+2^mf-1\),\\
      \(\!\!\!\!q=\left(2^{m+\alpha}+2^m\right)f-1\), & \(q=\left(2^{m\mp\alpha}f+2^m\right)f-1\),\\
      \(\!\!\!\!r=\left(2^{2m+1}+2^{2m+\alpha}+2^{2m-\alpha}\right)ff-1\), &
                                                                     \(r=\left(2^{2m+1}f+2^{2m\pm\alpha}+2^{2m\mp\alpha}ff\right)f-1\).\\
    \end{tabular}
  \end{center}

  \par However, note that, in order for \(2^{m+2}+1\) to be
  a prime number, the exponent \(m+2\) must be a power of two;
  therefore values of \(m\) will be \(0\), \(2\), \(6\), \(14\)
  etc. But the case \(m=0\) must be rejected, because no value of
  \(\alpha\) can be assigned.
\end{Case}

\newpage

\begin{Example}[1]
  Therefore let \(m=2\), so that \(a=8\cdot 17\) and
  \(b=4\cdot 17=68\) and also \(f=17\). Therefore, since it must
  be that\footnote{Typo \((x-b)(y-b)=4^2\cdot 17\) in original,
    but the error does not propagate. \JE}
  \((x-b)(y-b)=4^2\cdot 17^2\), it will be established by resolution into
  factors:

  \bgroup
  \def\arraystretch{1.1}
  \begin{center}
    \begin{tabular}{>{\raggedleft\arraybackslash}p{1.5cm}|
      >{\raggedleft\arraybackslash}p{1.2cm}@{}p{0.01cm}|
      >{\raggedleft\arraybackslash}p{1.2cm}@{}p{0.01cm}|
      >{\raggedleft\arraybackslash}p{1.2cm}@{}p{0.01cm}|
      >{\raggedleft\arraybackslash}p{1.2cm}@{}p{0.01cm}|}
      \(x-68=\) & \(2\) && \(4\) && \(8\) && \(34\)&\\
      \(y-68=\) & \(8\cdot 17\)&\(^2\) & \(1156\) && \(578\) &&
                                                         \(136\)&\\
      \(x=\) & \(70\) && \(72\) && \(76\) && \(102\)&\\
      \(y=\) & \(2380\) && \(1224\) && \(646\) && \(204\)&\\
      \(p=\) & \(69\)&\(^*\) & \(71\) && \(75\)&\(^*\) & \(101\)&\\
      \(q=\) & \(2379\)&\(^*\) & \(1223\) && \(645\)&\(^*\) & \(203\)&\(^*\)\\
      \(r=\) & \(166599\)&\(^*\) & \(88127\)&\(^*\) & \(49095\)&\(^*\) & \(20807\)&
    \end{tabular}
  \end{center}
  \egroup
  
  \par Hence, therefore, no amicable numbers will be obtained.
\end{Example}

\begin{Example}[2]
  Let \(m=6\), so that \(a=2^7\cdot 257\), \(b=2^6\cdot 257\) and
  \(f=257\). Therefore, since we get \[(x-b)(y-b)=2^{12}\cdot 257^2,\]
  the resolution must be established thus:

  \bgroup
  \def\arraystretch{1.1}
  \begin{center}
    \begin{tabular}{>{\raggedleft\arraybackslash}p{2cm}|
      >{\raggedleft\arraybackslash}p{1.6cm}@{}p{0.01cm}}
      \(x-16448=\) & \(32\cdot 257\)&\\
      \(y-16448=\) & \(128\cdot 257\)&\\
      \(x=\) & \(24672\)&\\
      \(y=\) & \(49344\)&\\
      \(p=\) & \(24671\)&\\
      \(q=\) & \(49343\)&\(^*\)\\
      \(r=\) & \(\cdots\)&\\
    \end{tabular}
  \end{center}
  \egroup

  \par The values arising from the remaining factors become
  still larger, so that it is difficult to judge whether they
  are prime or not.
\end{Example}

\begin{RemainingCases}
  Since \(f=2^{m+k+1}+2^k-1\) must be a prime number, we look
  first for the simpler cases in which this happens, since it is
  not possible to develop cases that are too complex. Therefore,
  let \(k=2\), and because \(f=2^{m+3}+3\), suitable values for
  \(m\) will be: \(1\), \(3\), \(4\). Let \(k=3\); we get
  \(f=2^{m+4}+7\) and suitable values for \(m\) will be \(2\),
  \(4\), \(6\). In the case \(k=4\) we have \(f=2^{m+5}+15\) and
  \(m\) will be either \(1\) or \(3\); and no further progress
  is possible.
\end{RemainingCases}

\begin{Example}[1]
  Suppose therefore \(k=2\) and \(m=1\); we get \(f=19\) and \(a=8\cdot
  19\) and also \(b=2\cdot 19=38\), whence we get
  \[(x-38)(y-38)=2^2\cdot 19^2=1444,\] and the resolutions will
  be given as:

  \bgroup
  \def\arraystretch{1.1}
  \begin{center}
    \begin{tabular}{>{\raggedleft\arraybackslash}p{1.5cm}|
      >{\raggedleft\arraybackslash}p{0.6cm}@{}p{0.01cm}|
      >{\raggedleft\arraybackslash}p{0.8cm}@{}p{0.01cm}|
      p{4.3cm}}
      \(x-38=\) & \(2\) && \(4\)&&\multirow{4}{4.5cm}{\begin{tabular}{p{4.3cm}@{}}Clearly neither factor may be assumed to be odd.\end{tabular}}\\
      \(y-38=\) & \(722\) && \(361\)&&\\
      \(x=\) & \(40\) && &&\\
      \(y=\) & \(760\) && odd\footnotemark &&\\
      \(p=\) & \(39\)&\(^*\) &&
    \end{tabular}
    \footnotetext{In the original this is ``imp:'' which could be an abbreviation for {\em impar} (odd) or {\em impossibile} (impossible). \JE}
  \end{center}
  \egroup

  \par Because here already \(p\) is not prime, it is clear
  that no amicable numbers result from this.
\end{Example}

\begin{Example}[2]
  Suppose \(k=2\) and \(m=3\), so that \(f=67\); we get
  \(a=32\cdot 67\) and \(b=8\cdot 67=536\), whence we get
  \[(x-536)(y-536)=2^6\cdot 67^2.\]
  
  \bgroup
  \def\arraystretch{1.1}
  \begin{center}
    \begin{tabular}{>{\raggedleft\arraybackslash}p{2cm}|
      >{\raggedleft\arraybackslash}p{1cm}@{}p{0.01cm}|
      >{\raggedleft\arraybackslash}p{1cm}@{}p{0.01cm}|
      p{6.5cm}}
      \(x-536=\) & \(268\) &&
                             \(16\)&&\multirow{6}{5cm}{\begin{tabular}{p{6.5cm}@{}}Remaining
                                                        values of \(p\)
                                                        yield numbers divisible by \(3\)
                                                        and are omitted on account
                                                        of that.
                                                        Further
                                                        examples
                                                        lead to
                                                        exceedingly
                                                        large numbers.\end{tabular}}\\
      \(y-536=\) & \(1072\) && \(17956\)&&\\
      \(x=\) & \(804\) && \(552\) &&\\
      \(y=\) & \(1608\) && \(\cdots\)&&\\
      \(p=\) & \(803\)&\(^*\) & \(551\)&\(^*\)\\
      \(q=\) & \(1607\) && \(\cdots\)&      
    \end{tabular}
  \end{center}
  \egroup
\end{Example}

\begin{Rule}[3]
  As before, let \(a=2^n\left(2^{n+1}+2^k-1\right)\) and
  \(2^{n+1}+2^k-1=f\) be a prime number, but in the fraction
  \(\dfrac{b}{c}=\dfrac{2^n\left(2^{n+1}+2^k-1\right)}{2^k}\)
  let \(k>n\); then we get
  \[b=2^{n+1}+2^k-1\et c=2^{k-n}.\]
  Suppose \(k-n=m\), so that \(k=m+n\); we get
  \[a=2^n\left(2^{n+1}+2^{m+n}-1\right),\quad
    b=2^{n+1}+2^{m+n}-1=f\et c=2^m,\]
  whence we will have this equation
  \[\left(2^mx-b\right)\left(2^my-b\right)=bb.\] But since
  \(b=f\) is a prime number, no resolution takes place other
  than \(1\cdot bb\), from which we get
  \begin{alignat*}{4}
    x&=\frac{1+b}{2^m}&\qquad\mbox{and}&\qquad&
                                        y&=\frac{b(1+b)}{2^m}&&\mbox{or}\\
    x&=2^n+2^{n+1-m}&\qquad\mbox{and}&\qquad&
                         y&=\left(2^{n+1}+2^{m+n}-1\right)&&\left(2^n+2^{n+1-m}\right).
  \end{alignat*}
  \par Now note that these four numbers must be prime
  \begin{gather*}
    f=2^{n+1}+2^{m+n}-1,\\
    p=x-1,\quad q=y-1,\et r=xy-1
  \end{gather*}
  and it is also necessary that \(m<n+1\). If these conditions
  are satisfied, we will get the amicable numbers \(apq\) and
  \(ar\).
\end{Rule}

\begin{Case}[1]
  Let \(m=1\); we get \(f=2^{n+2}-1\), \(x=2^{n+1}\), and
  \(p=2^{n+1}-1\); but it is not possible that both \(f\) and
  \(p\) are simultaneously prime numbers except in the case
  \(n=1\), which however gives \(q=27\). Therefore from the
  hypothesis \(m=1\) no amicable numbers arise.
\end{Case}

\begin{Case}[2]
  Therefore let \(m=2\), so that
  \[f=3\cdot 2^{n+1}-1,\quad x=3\cdot 2^{n-1}\et y=3\cdot
    2^{n-1}\left(3\cdot 2^{n+1}-1\right)\tx{and}
    a=2^n\cdot f.\] Therefore the following four numbers must be
  prime:
  \begin{align*}
    f&=3\cdot 2^{n+1}-1, &q&=3\cdot 2^{n-1}\left(3\cdot 2^{n+1}-1\right)-1,\\
    p&=3\cdot 2^{n-1}-1,
                         &r&=9\cdot 2^{2n-2}\left(3\cdot 2^{n+1}-1\right)-1,
  \end{align*}
  whence these examples are laid out as follows:

  \bgroup
  \def\arraystretch{1.1}
  \begin{center}
    \begin{tabular}{>{\raggedleft\arraybackslash}p{0.7cm}|
      >{\raggedleft\arraybackslash}p{0.8cm}@{}p{0.01cm}|
      >{\raggedleft\arraybackslash}p{0.8cm}@{}p{0.01cm}|
      >{\raggedleft\arraybackslash}p{0.8cm}@{}p{0.01cm}|
      >{\raggedleft\arraybackslash}p{0.8cm}@{}p{0.01cm}|
      >{\raggedleft\arraybackslash}p{0.8cm}@{}p{0.01cm}|}
      \hline\hline
      \(n=\) & \(1\) && \(2\) && \(3\) && \(4\) && \(5\)&\\
      \hline\hline
      \(f=\) & \(11\) && \(23\) && \(47\) && \(95\)&\(^*\) & \(191\)&\\
      \(p=\) & \(2\) && \(5\) && \(11\) && \(\cdots\) && \(47\)&\\
      \(q=\) & \(32\)&\(^*\) & \(137\) && \(563\) && \(\cdots\) && \(9167\)&\(^*\)\\
      \(r=\) & \(98\)&\(^*\) & \(827\) && \(6767\)&\(^*\) & \(\cdots\) && \(\cdots\)&\\
      \multicolumn{1}{c}{} &\multicolumn{1}{c}{}&\multicolumn{1}{c}{}&\multicolumn{1}{c}{valid}&\multicolumn{1}{c}{}&\multicolumn{1}{c}{}&\multicolumn{1}{c}{}&\multicolumn{1}{c}{}&\multicolumn{1}{c}{}&\multicolumn{1}{c}{}&\multicolumn{1}{c}{}
    \end{tabular}
  \end{center}
  \egroup

  \noindent and hence therefore from \(n=2\) and \(a=4\cdot 23\)
  we obtain the amicable numbers
  \[\amn{4\cdot 23\cdot 5\cdot 137}{4\cdot 43\cdot 827.}\]
\end{Case}

\begin{RemainingCases}
  Let \(m=3\), again either \(f\) or \(p\) becomes divisible by
  \(3\), and the same happens if \(m=5\), or \(7\),
  etc. Therefore let \(m=4\); we get
  \[f=9\cdot 2^{n+1}-1,\quad x=9\cdot 2^{n-3}\et y=9\cdot
    2^{n-3}\left(9\cdot 2^{n+1}-1\right)\et a=2^n\cdot f,\]
  whence these examples are laid out as follows:

  \bgroup
  \def\arraystretch{1.1}
  \begin{center}
    \begin{tabular}{>{\raggedleft\arraybackslash}p{0.7cm}|
      >{\raggedleft\arraybackslash}p{1cm}@{}p{0.01cm}|
      >{\raggedleft\arraybackslash}p{1cm}@{}p{0.01cm}|
      >{\raggedleft\arraybackslash}p{1cm}@{}p{0.01cm}|
      >{\raggedleft\arraybackslash}p{1cm}@{}p{0.01cm}|}
      \hline\hline
      \(n=\) & \(1\) && \(4\) && \(5\) && \(6\) &\\
      \hline\hline
      \(f=\) & \(35\)&\(^*\) & \(287\)&\(^*\) & \(575\)&\(^*\) & \(1151\)& \\
      \(x=\) & \(\cdots\) && \(\cdots\) && \(\cdots\) && \(72\)& \\
      \(y=\) & \(\cdots\) && \(\cdots\) && \(\cdots\) && \(82871\)& \\
      \(p=\) & \(\cdots\) && \(\cdots\) && \(\cdots\) && \(71\)&\\
      \(q=\) & \(\cdots\) && \(\cdots\) && \(\cdots\) && \(82871\)&\(^*\)\\
      \(r=\) & \(\cdots\) && \(\cdots\) && \(\cdots\) && \(\cdots\)&\\
    \end{tabular}
  \end{center}
  \egroup

  \par Neither from this nor from taking larger values of
  \(m\) is it possible to elicit amicable numbers.
\end{RemainingCases}

\begin{Rule}[4]
  Still more expressions for the common factor can be found, from
  which the denominator \(c\) of the fraction \(\dfrac{b}{c}\)
  becomes equal to either unity or a power of two. For indeed, let
  us imagine that \(a=2^n(g-1)(h-1)\), where \(g-1\)
  and \(h-1\) are prime numbers; we get
  \[\int a=\left(2^{n+1}-1\right)gh=2^{n+1}gh-gh;\] but
  \(2a=2^{n+1}gh-2^{n+1}g-2^{n+1}h+2^{n+1}\), whence we get
  \[2a-\int a=gh-2^{n+1}g-2^{n+1}h+2^{n+1}.\] Suppose
  \(2a-\int a=d\); we get \(gh-2^{n+1}(g+h)+2^{n+1}=d\)
  and\footnote{Typo
    \(\left(g-2^{n+1}\right)\left(h=2^{n+1}\right)=d-2^{n+1}+2^{2n+2}\)
    in original. \JE}
  \[\left(g-2^{n+1}\right)\left(h-2^{n+1}\right)=d-2^{n+1}+2^{2n+2};\]
  whence, by resolving into factors, values for \(g\) and \(h\)
  must be found so that \(g-1\) and \(h-1\) become prime
  numbers, and then we get
  \[a=2^n(g-1)(h-1)\et\frac{b}{c}=\frac{a}{d}.\]
  \begin{itemize}[align=right,itemindent=2em,labelsep=2pt,labelwidth=1em,leftmargin=0pt,nosep]
  \item[I.] Suppose \(n=1\); we get
    \[(g-4)(h-4)=d+12;\] now splitting \(d+12\) into two
    even factors, the following values will appear:

    \medskip
    
    Let \(d=4\); we get
    \begin{gather*}
      (g-4)(h-4) = 16=2\cdot 8,\tx{whence}g=6,\quad h=12,\\
      a=2\cdot 5\cdot 11\tx{and}\frac{b}{c}=\frac{2\cdot
        5\cdot 11}{4};\tx{therefore}b=5\cdot 11\et c=2.
    \end{gather*}
    Let \(d=8\); we get
    \begin{gather*}
      (g-4)(h-4) = 20=2\cdot 10,\tx{whence}g=6,\quad h=14,\\
      a=2\cdot 5\cdot 13\tx{and}\frac{b}{c}=\frac{2\cdot
        5\cdot 13}{8};\tx{therefore}b=5\cdot 13\et c=4.
    \end{gather*}
    Let \(d=16\); we get
    \begin{gather*}
      (g-4)(h-4) = 28=2\cdot 14,\tx{whence}g=6,\quad h=18,\\
      a=2\cdot 5\cdot 17\tx{and}\frac{b}{c}=\frac{2\cdot
        5\cdot 17}{16};\tx{therefore}b=5\cdot 17\et c=8.
    \end{gather*}
  \item[II.] Suppose \(n=2\); we get
    \[(g-8)(h-8)=d+56\] and \(a=4(g-1)(h-1)\), whence the
    following cases result:

    \medskip
    
    Let \(d=4\); we get
    \begin{gather*}
      (g-8)(h-8) = 60=6\cdot 10,\tx{whence}g=14,\quad h=18,\\
      a=4\cdot 13\cdot 17\tx{and}\frac{b}{c}=\frac{4\cdot
        13\cdot 17}{4};\tx{therefore}b=13\cdot 17\et c=1.
    \end{gather*}
    Let \(d=8\); we get
    \begin{gather*}
      (g-8)(h-8) = 64=4\cdot 16,\tx{whence}g=12,\quad h=24,\\
      a=4\cdot 11\cdot 23\tx{and}\frac{b}{c}=\frac{4\cdot
        11\cdot 23}{8};\tx{therefore}b=11\cdot 23\et c=2.
    \end{gather*}
    Let \(d=16\); we get
    \begin{gather*}
      (g-8)(h-8) = 72=6\cdot 12,\tx{whence}g=14,\quad h=20,\\
      a=4\cdot 13\cdot 19\tx{and}\frac{b}{c}=\frac{4\cdot
        13\cdot 19}{16};\tx{therefore}b=13\cdot 19\et c=4.
    \end{gather*}
  \item[III.] Suppose \(n=3\), so that \(a=8(g-1)(h-1)\),
    and it will have to be
    \[(g-16)(h-16)=d+240.\]
    Let \(d=4\); we get
    \begin{gather*}
      (g-16)(h-16) = 244=2\cdot 122,\tx{whence}g=18,\quad h=138,\\
      a=8\cdot 17\cdot 137\smalltx{and}\frac{b}{c}=\frac{8\cdot
        17\cdot 137}{4};\smalltx{therefore}b=2\cdot 17\cdot 137\smalltx{and} c=1.
    \end{gather*}
    Let \(d=8\); we get
    \begin{gather*}
      (g-16)(h-16) = 248=2\cdot 124,\tx{whence}g=18,\quad h=140,\\
      a=8\cdot 17\cdot 139\smalltx{and}\frac{b}{c}=\frac{8\cdot
        17\cdot 139}{8};\smalltx{therefore}b=17\cdot 139\smalltx{and} c=1.
    \end{gather*}
    Let \(d=16\); we get
    \begin{gather*}
      (g-16)(h-16) = 256=4\cdot 64,\tx{whence}g=20,\quad h=80,\\
      a=8\cdot 19\cdot 79\smalltx{and}\frac{b}{c}=\frac{8\cdot
        19\cdot 79}{16};\smalltx{therefore}b=19\cdot 79\smalltx{and} c=2.
    \end{gather*}
    Again, let \(d=16\) and
    \begin{gather*}
      (g-16)(h-16) = 8\cdot 32,\tx{whence}g=24,\quad h=48,\\
      a=8\cdot 23\cdot 47\smalltx{and}\frac{b}{c}=\frac{8\cdot
        23\cdot 47}{16};\smalltx{therefore}b=23\cdot 47\smalltx{and} c=2.
    \end{gather*}
  \end{itemize}
  \par Indeed, by taking values for \(a\) in this way, if
  \(a(x-1)(y-1)\) and \(a(xy-1)\) were established as amicable
  numbers, so that \(x-1\), \(y-1\), and \(xy-1\) are prime
  numbers, it is necessary that \((cx-b)(cy-b)=bb\).
\end{Rule}

\begin{Example}[1]
  Let \(a=2\cdot 5\cdot 11\); we get \(b=5\cdot 11=55\) and
  \(c=2\), whence we get
  \[(2x-55)(2y-55)=5^2\cdot 11^2.\]
  \bgroup
  \def\arraystretch{1.1}
  \begin{center}
    \begin{tabular}{>{\raggedleft\arraybackslash}p{1.2cm}|
      >{\raggedleft\arraybackslash}p{1cm}@{}p{0.01cm}|
      >{\raggedleft\arraybackslash}p{1cm}@{}p{0.01cm}|
      >{\raggedleft\arraybackslash}p{1cm}@{}p{0.25cm}|
      p{5cm}}
      \(2x-55\) & \(1\) && \(5\)&&\(25\)&
      &\multirow{7}{5cm}{\begin{tabular}{p{5cm}@{}}Hence
                           therefore no amicable numbers are obtained.\end{tabular}}\\
      \(2y-55\) & \(3025\) && \(605\)&& \(121\)&)\footnotemark&\\
      \(x\) & \(28\) && \(30\) && \(40\) &&\\
      \(y\) & \(1540\) && \(330\)&& \(88\) &&\\
      \(x-1\) & \(27\)&\(^*\) & \(29\) && \(39^*\) &&\\
      \(y-1\) & \(\cdots\) && \(329\)&\(^*\) & \(\cdots\) &&\\
      \(xy-1\) & \(\cdots\) && \(\cdots\) && \(\cdots\) &&
    \end{tabular}
  \end{center}
  \footnotetext{In the original table, the following typos
    appear: \(125\) instead of \(121\), \(90\) instead of
    \(88\). \JE} \egroup
\end{Example}

\newpage

\begin{Example}[2]
  Let \(a=2\cdot 5\cdot 13\); we get \(b=5\cdot 13=65\) and
  \(c=4\), whence we get
  \[(4x-65)(4y-65) = 5^2\cdot 13^2.\]

  \par But this number
  \(5^2\cdot 13^2\) cannot be resolved in two factors which,
  when increased by \(65\) become divisible by \(4\); the same
  applies to the value \(a=2\cdot 5\cdot 17\).
\end{Example}

\begin{Example}[3]
  Let \(a=4\cdot 13\cdot 17\); we get \(b=13\cdot 17=221\) and
  \(c=1\) and it must be that \((x-221)(y-221)=13^2\cdot 17^2\),
  whence

  \bgroup
  \def\arraystretch{1.1}
  \begin{center}
    \begin{tabular}{>{\raggedleft\arraybackslash}p{1.2cm}|
      >{\raggedleft\arraybackslash}p{1cm}@{}p{0.01cm}|
      >{\raggedleft\arraybackslash}p{1cm}@{}p{0.01cm}|
      >{\raggedleft\arraybackslash}p{1.3cm}@{}p{0.01cm}|}
      \(x-221\) & \(13\) && \(17\) && \(169\)&\\
      \(y-221\) & \(3757\) && \(\cdots\) && \(289\)&\\
      \(x-1\) & \(233\) && \(237\)&\(^*\) & \(389\)&\\
      \(y-1\) & \(3977\)&\(^*\) & \(\cdots\) && \(509\)&\\
      \(xy-1\) & \(\cdots\) && \(\cdots\) && \(198899\)&\\
    \end{tabular}
  \end{center}
  \egroup
      
  In the final factorisation, \(x-1\) and \(y-1\) are prime
  numbers so the question reduces to whether \(xy-1=198899\) is
  a prime number or not. But even though this number exceeds the
  limit of \(100000\), I can nevertheless show it to be prime,
  whence the amicable numbers will be
  \[\amn{4\cdot 13\cdot 17\cdot 389\cdot 509}{4\cdot
      13\cdot 17\cdot 198899.}\]
\end{Example}

\begin{Scholium}
  I infer this number \(198899\) to be prime because I have
  observed that \(198899=2\cdot 47^2+441^2\), so that \(198899\)
  is a number of the form \(2aa+bb\). But it is certain that if
  a number can be put into the form \(2aa+bb\) in a unique way,
  then it is prime; otherwise, if it may be reduced to the form
  \(2aa+bb\) in two or more ways, then it is
  composite.\footnote{See {\cite[Theorem 10]{euler256}}. \FR}
  Therefore I have looked for whether from this number
  \(198899\) any doubled square other than \(47^2\) may be
  subtracted to leave a square residue, and, by a drawn-out
  calculation, I found none; from which I have safely concluded
  this number to be prime, and so the discovered numbers to be
  amicable. From the remaining values of \(a\) which I listed,
  no amicable numbers can be found.
\end{Scholium}

\newpage

\begin{Rule}[5]
  Yet more numbers may be suitable to take as values for \(a\),
  from which amicable numbers can be extracted. But since
  general rules cannot be formulated for them, I will only work
  out some of them, by imitation of which it is not difficult to
  contrive others.

  \begin{itemize}
  \item[I.] Let therefore \(a=3^2\cdot 5\cdot 13\); we get
    \(\int a=13\cdot 6\cdot 14\) and because \(2a=90\cdot 13\)
    and \(\int a=84\cdot 13\) we get \(2a-\int a=6\cdot 13\) and
    also\footnote{Typo \(\frac{a}{-2a\int a}\) in original. \JE}
    \[\dfrac{b}{c}=\dfrac{a}{2a-\int a}=\dfrac{3^2\cdot 5\cdot
      13}{6\cdot 13}=\dfrac{15}{2}\tx{and so}b=15\tx{and}c=2.\]
  \item[II.] Let \(a=3^2\cdot 7\cdot 13\); we get \(\int a=13\cdot
    8\cdot 14=16\cdot 7 \cdot 13\), whence because \(2a=18\cdot 7\cdot
    13\) we get \(2a-\int a=2\cdot 7\cdot 13\) and so
    \[\dfrac{b}{c}=\dfrac{3^2\cdot 7\cdot 13}{2\cdot 7\cdot
        13}=\dfrac{9}{2},\tx{whence}b=9\et c=2.\]
  \item[III.] Let \(a=3^2\cdot 7^2\cdot 13\); we get
    \(\int a=13\cdot 3\cdot 19\cdot 14=2\cdot 3\cdot 7\cdot
    13\cdot 19\) and\footnote{Typo: \(4\cdot 2\) instead of \(42\)
      in original. \JE} \(2a=42\cdot 3\cdot 7\cdot 13\), whence
    \(2a-\int a=4\cdot 3\cdot 7\cdot 13\) and so
    \[\dfrac{b}{c}=\dfrac{3^2\cdot 7^2\cdot 13}{4\cdot 3\cdot
        7\cdot 13}=\dfrac{21}{4},\tx{therefore}b=21\et c=4.\]
  \item[IV.] Let \(a=3^3\cdot 5\); we get
    \(\int a=5\cdot 8\cdot 6=16\cdot 3\cdot 5\). Therefore
    because \(2a=18\cdot 3\cdot 5\) we get
    \(2a-\int a=2\cdot 3\cdot 5\) and hence
    \[\dfrac{b}{c}=\dfrac{3^3\cdot 5}{2\cdot 3\cdot
      5}=\dfrac{9}{2}\et b=9\et c=2.\]
  \item[V.] Let \(a=3^2\cdot 5\cdot 13\cdot 19\); we get
    \(\int a=13\cdot 6\cdot 14\cdot 20=16\cdot 3\cdot 5\cdot
    7\cdot 13\) and because \(2a=114\cdot 3\cdot 5\cdot 13\) and
    \(\int a=112\cdot 3\cdot 5 \cdot 13\) we get
    \[\dfrac{b}{c}=\dfrac{3^2\cdot 5\cdot 13\cdot 19}{2\cdot
      3\cdot 5\cdot 13}=\dfrac{3\cdot 19}{2}\et
    b=3\cdot 19=57\et c=2.\]
\item[VI.] Let \(a=3^2\cdot 7^2\cdot 13\cdot 19\); we get
  \(\int a=13\cdot 3\cdot 19\cdot 14\cdot 20=8\cdot 3\cdot
  5\cdot 7 \cdot 13\cdot 19\) and because
  \(2a=42\cdot 3\cdot 7\cdot 13\cdot 19\) we
  get
  \[\dfrac{b}{c}=\dfrac{3^2\cdot 7^2\cdot 13\cdot 19}{2\cdot
      3\cdot 7\cdot 13\cdot 19}=\dfrac{21}{2},\tx{whence we
      get}b=21\et c=2.\]
  \end{itemize}

  \par However, supposing \(a(x-1)(y-1)\) and \(a(xy-1)\) to
  be amicable numbers, it must be that \((cx-b)(cy-b)=bb\).
\end{Rule}

\begin{Example}[1]
  Let \(b=15\), \(c=2\); we get \(a=3^2\cdot 5\cdot 13\) and
  this equation must be satisfied \((2x-15)(2y-15)=225\).

  \bgroup
  \def\arraystretch{1.1}
  \begin{center}
    \begin{tabular}{>{\raggedleft\arraybackslash}p{1.2cm}|
      >{\raggedleft\arraybackslash}p{1cm}@{}p{0.01cm}|
      >{\raggedleft\arraybackslash}p{1cm}@{}p{0.01cm}|
      >{\raggedleft\arraybackslash}p{1.3cm}@{}p{0.01cm}|
      p{6cm}}
      \(2x-15\) & \(1\) && \(5\) && \(9\) &\\
      \(2y-15\) & \(225\) && \(45\) && \(25\) &\\
      \(x\) & \(8\) && \(10\) && \(12\) & \\
      \(y\) & \(120\) && \(30\) && \(20\) & \\
      \(x-1\) & \(7\) && \(9\)&\(^*\) & \(11\) & \\
      \(y-1\) & \(119\)&\(^*\) & \(\cdots\) && \(19\) & \\
      \(xy-1\) & \(\cdots\) && \(\cdots\) && \(239\) &
    \end{tabular}
  \end{center}
  \egroup
  
  Therefore the amicable numbers will be
  \(\amn{3^2\cdot 5\cdot 13\cdot 11\cdot
      19}{3^2\cdot 5\cdot 13\cdot 239.}\)
\end{Example}

\begin{Example}[2]
  Let \(b=9\), \(c=2\); we get \(a=3^2\cdot 7\cdot 13\) or
  \(a=3^3\cdot 5\) and the equation to be resolved is
  \((2x-9)(2y-9)=81\).

  \bgroup
  \def\arraystretch{1.1}
  \begin{center}
    \begin{tabular}{>{\raggedleft\arraybackslash}p{1.2cm}|
      >{\raggedleft\arraybackslash}p{1cm}@{}p{0.01cm}|
      p{6cm}}
      \(2x-9\) & \(3\) && 
                          \multirow{7}{6cm}{\begin{tabular}{p{6cm}@{}}
                                              So, since \(x-1=5\),
                                              this value cannot
                                              be combined with
                                              \(a=3^3\cdot 5\).
                                              Therefore the amicable numbers will be \[\amn{3^2\cdot
                                              7\cdot
                                              13\cdot
                                              5\cdot
                                              17}{3^2\cdot
                                              7\cdot
                                              13\cdot
                                              107.}\]
                                             \end{tabular}}\\
      \(2y-9\) & \(27\) &&\\
      \(x\) & \(6\) &&  \\
      \(y\) & \(18\) &&\\
      \(x-1\) & \(5\) &&\\
      \(y-1\) & \(17\) &&\\
      \(xy-1\) & \(107\) &&
    \end{tabular}
  \end{center}
  \egroup
\end{Example}

\begin{Example}[3]
  Let \(b=21\) and \(c=4\); we get \(a=3^2\cdot 7^2\cdot 13\) and
  the equation to be resolved is \((4x-21)(4y-21)=441\).

  \bgroup
  \def\arraystretch{1.1}
  \begin{center}
    \begin{tabular}{>{\raggedleft\arraybackslash}p{1.2cm}|
      >{\raggedleft\arraybackslash}p{1cm}@{}p{0.01cm}|
      p{6.5cm}}
      \(4x-21\) & \(3\) && 
                          \multirow{7}{6.5cm}{\begin{tabular}{p{6.5cm}@{}}\qquad 
                                                Because
                                                \(x\) and \(y\)
                                                must be even
                                                numbers, other
                                                resolutions do
                                                not take place.\\  
                                                \qquad Hence, therefore, these amicable numbers appear: \\\centering\(\amn{3^2\cdot
                                                 7^2\cdot
                                                 13\cdot
                                                 5\cdot
                                                 41}{3^2\cdot
                                                 7^2\cdot
                                                 13\cdot
                                                 251.}\)\end{tabular}}\\
      \(4y-21\) & \(147\) &&\\
      \(x\) & \(6\) &&  \\
      \(y\) & \(42\) &&\\
      \(x-1\) & \(5\) &&\\
      \(y-1\) & \(41\) &&\\
      \(xy-1\) & \(251\) &&
    \end{tabular}
  \end{center}
  \egroup
\end{Example}

\begin{Example}[4]
  Let \(b=21\) and \(c=2\); we get \(a=3^2\cdot 7^2\cdot 13\cdot 19\) and
  the equation to be resolved is \((2x-21)(2y-21)=441\).

  \bgroup
  \def\arraystretch{1.1}
  \begin{center}
    \begin{tabular}{>{\raggedleft\arraybackslash}p{1.2cm}|
      >{\raggedleft\arraybackslash}p{1cm}@{}p{0.01cm}|
      >{\raggedleft\arraybackslash}p{1cm}@{}p{0.01cm}|
      p{6.5cm}}
      \(2x-21\) & \(3\) && \(7\) && 
                                    \multirow{7}{6.5cm}{\begin{tabular}{p{6.5cm}@{}}\qquad 
                                                          But
                                                          because
                                                          the value 
                                                          \(x-1=13\) is
                                                          already
                                                          contained
                                                          in \(a\),
                                                          we obtain hence
                                                          no
                                                          amicable
                                                          numbers.\end{tabular}}\\
      \(2y-21\) & \(147\) && \(63\) &&\\
      \(x\) & \(12\) &&  \(14\) &&\\
      \(y\) & \(84\) && \(42\) &&\\
      \(x-1\) & \(11\) && \(13\) &&\\
      \(y-1\) & \(83\) && \(41\) && \\
      \(xy-1\) & \(1007\)&\(^*\) & \(587\) &&
    \end{tabular}
  \end{center}
  \egroup
\end{Example}

\begin{Example}[5]
  Let \(b=57\) and \(c=2\); we get \(a=3^2\cdot 5\cdot 13\cdot 19\)
  and the equation to be resolved is \((2x-57)(2y-57)=3249\).

  \bgroup
  \def\arraystretch{1.1}
  \begin{center}
    \begin{tabular}{>{\raggedleft\arraybackslash}p{1.2cm}|
      >{\raggedleft\arraybackslash}p{1cm}@{}p{0.01cm}|
      >{\raggedleft\arraybackslash}p{1cm}@{}p{0.3cm}|
      p{6.5cm}}
      \(2x-57\) & \(3\) && \(19\) && \multirow{7}{6cm}{\begin{tabular}{p{5cm}@{}}\qquad Hence, therefore, these amicable numbers arise: \[\amn{3^2\cdot 5\cdot 13\cdot 19\cdot 29\cdot 569}{3^2\cdot 5\cdot 13\cdot 19\cdot 17099.}\]\end{tabular}}\\
      \(2y-57\) & \(1083\) && \(171\) && \\
      \(x\) & \(30\) && \(38\) &&  \\
      \(y\) & \(570\) && \(114\) &&\\
      \(x-1\) & \(29\) && \(37\)&)\footnotemark &\\
      \(y-1\) & \(569\) && \(113\) &&\\
      \(xy-1\) & \(17099\) && \(4331\)&\(^*\) &
    \end{tabular}
  \end{center}
  \footnotetext{Typo: \(34\) appears instead of \(37\) in the original. \JE}
  \egroup
\end{Example}

\begin{Example}[6]
  Let \(b=45\) and \(c=2\); we get \(a=3^4\cdot 5\cdot 11\)
  and the equation to be resolved is \((2x-45)(2y-45)=2025\).

  \bgroup
  \def\arraystretch{1.1}
  \begin{center}
    \begin{tabular}{>{\raggedleft\arraybackslash}p{1.2cm}|
      >{\raggedleft\arraybackslash}p{1cm}@{}p{0.01cm}|
      >{\raggedleft\arraybackslash}p{1cm}@{}p{0.01cm}|
      p{6cm}}
      \(2x-45\) & \(3\) && \(15\) &&  
                                     \multirow{7}{6cm}{\begin{tabular}{p{6cm}@{}}\qquad 
                                                         Hence, therefore,
                                                         arise the
                                                         amicable numbers: \[\amn{3^4\cdot
                                                         5\cdot
                                                         11\cdot
                                                         29\cdot
                                                         89}{3^4\cdot
                                                         5\cdot
                                                         11\cdot
                                                         2699.}\]\end{tabular}}\\
      \(2y-45\) & \(675\) && \(135\) && \\
      \(x\) & \(24\) && \(30\) &&  \\
      \(y\) & \(360\) && \(90\) &&\\
      \(x-1\) & \(23\) && \(29\) &&\\
      \(y-1\) & \(359\) && \(89\) &&\\
      \(xy-1\) & \(8639\)&\(^*\) & \(2699\) &&
    \end{tabular}
  \end{center}
  \egroup
\end{Example}

\begin{Example}[7]
  Let \(b=77\) and \(c=2\); we get
  \(a=3^2\cdot 7^2\cdot 11\cdot 13\) and the equation to be
  resolved is
  \((2x-77)(2y-77)=49\cdot 121\).

  \bgroup
  \def\arraystretch{1.1}
  \begin{center}
    \begin{tabular}{>{\raggedleft\arraybackslash}p{1.2cm}|
      >{\raggedleft\arraybackslash}p{1cm}@{}p{0.01cm}|
      >{\raggedleft\arraybackslash}p{1cm}@{}p{0.01cm}|
      p{6cm}}
      \(2x-77\) & \(7\) && \(11\) &&  
                                   \multirow{7}{6cm}{\begin{tabular}{p{5cm}@{}}\qquad 
                                                       Hence, therefore,
                                                       arise the
                                                       amicable numbers: \[\amn{3^2\cdot
                                                           7^2\cdot
                                                           11\cdot
                                                           13\cdot
                                                           41\cdot
                                                           461}{3^2\cdot
                                                           7^2\cdot
                                                           11\cdot
                                                           13\cdot
                                                           19403.}\]\end{tabular}}\\
      \(2y-77\) & \(847\) && \(539\) && \\
      \(x\) & \(42\) && \(44\) &&  \\
      \(y\) & \(462\) && \(308\) &&\\
      \(x-1\) & \(41\) && \(43\) &&\\
      \(y-1\) & \(461\) && \(307\) &&\\
      \(xy-1\) & \(19403\) && \(13551\)&\(^*\) &
    \end{tabular}
  \end{center}
  \egroup
\end{Example}

\begin{Example}[8]
  Let \(b=105\), \(c=2\); we get \(a=3^2\cdot 5\cdot 7\) and the
  equation to be resolved is \((2x-105)(2y-105)=105^2\).

  \bgroup
  \def\arraystretch{1.1}
  \begin{center}
    \begin{tabular}{>{\raggedleft\arraybackslash}p{1.5cm}|
      >{\raggedleft\arraybackslash}p{1cm}@{}p{0.01cm}|
      >{\raggedleft\arraybackslash}p{1cm}@{}p{0.01cm}|
      >{\raggedleft\arraybackslash}p{1cm}@{}p{0.01cm}|
      >{\raggedleft\arraybackslash}p{1cm}@{}p{0.01cm}|
      p{5cm}}
      \(2x-105\) & \(3\) && \(7\) && \(15\) && \(35\) &&  
                                                     \multirow{7}{5cm}{\begin{tabular}{p{5cm}@{}}
                                                                         \qquad Since
                                                                         \(102059\)
                                                                         is a prime number,
                                                                         because
                                                                         it
                                                                         can
                                                                         be
                                                                         put
                                                                         in
                                                                         the
                                                                         form
                                                                         \(8a+3\)
                                                                         and
                                                                         in
                                                                         a
                                                                         unique
                                                                         way
                                                                         in
                                                                         the
                                                                         form
                                                                         \(2aa+bb\),
                                                                         the
                                                                         amicable
                                                                         numbers
                                                                         arising
                                                                         hence
                                                                         will
                                                                         be \\
                                                                         \centering\(\amn{3^2\cdot
                                                                             5\cdot
                                                                             7\cdot
                                                                             53\cdot
                                                                             1889}{3^2\cdot
                                                                             5\cdot
                                                                             7\cdot
                                                                             102059.}\)\end{tabular}}\\
      \(2y-105\) & \(3675\) && \(\cdots\) && \(735\) && \(\cdots\) && \\
      \(x\) & \(54\) && \(56\) && \(60\) && \(70\) &&  \\
      \(y\) & \(1890\) && \(\cdots\) && \(420\) && \(\cdots\) &&\\
      \(x-1\) & \(53\) && \(55\)&\(^*\) & \(59\) && \(69\)&\(^*\) &\\
      \(y-1\) & \(1889\) && \(\cdots\) && \(419\) && \(\cdots\) &&\\
      \(xy-1\) & \(102059\) && \(\cdots\) && \(25199\)&\(^*\) & \(\cdots\) &&
    \end{tabular}
  \end{center}
  \egroup
\end{Example}

\begin{Scholium}
  Therefore the amicable numbers of the form \(apq\), \(ar\)
  which we have found thus far are:

  \begin{center} 
    \begingroup
    \renewcommand{\arraystretch}{1.5} % Default value: 1
    \setlength\tabcolsep{0pt}
    \begin{tabular}{p{4cm}p{5cm}p{4.5cm}}
      \centering \(\mbox{I.}\,\amn{2^2\cdot 5\cdot 11}{2^2\cdot 71}\) & \centering \(\mbox{II.}\,\amn{2^4\cdot 23\cdot 47}{2^4\cdot 1151}\) & \centering \(\mbox{III.}\,\amn{2^7\cdot 191\cdot 383}{2^7\cdot 73727}\)\arraybackslash \\
      \smallskip\centering \(\mbox{IV.}\,\amn{4\cdot 23\cdot 137}{4\cdot 23\cdot 827}\) & \smallskip\centering \(\mbox{V.}\,\amn{4\cdot 13\cdot 17\cdot 389\cdot 509}{4\cdot 13\cdot 17\cdot 198899}\) & \smallskip\centering \(\mbox{VI.}\,\amn{3^2\cdot 5\cdot 13\cdot 11\cdot 19}{3^2\cdot 5\cdot 13\cdot 239}\)\arraybackslash \\
      \smallskip\centering \(\mbox{VII.}\,\amn{3^2\cdot 7\cdot 13\cdot 5\cdot 17}{3^2\cdot 7\cdot 13\cdot 107}\) & \smallskip\centering \(\!\!\!\!\mbox{VIII.}\,\amn{3^2\cdot 7^2\cdot 13\cdot 5\cdot 41}{3^2\cdot 7^2\cdot 13\cdot 251}\) & \smallskip\centering \(\!\!\!\!\mbox{IX.}\,\amn{3^2\cdot 5\cdot 13\cdot 19\cdot 29\cdot 569}{3^2\cdot 5\cdot 13\cdot 19\cdot 17099}\)\arraybackslash \\
      \smallskip\centering \(\mbox{X.}\,\amn{3^4\cdot 5\cdot 11\cdot 29\cdot 89}{3^4\cdot 5\cdot 11\cdot 2699}\) & \smallskip\centering \(\mbox{XI.}\,\amn{3^2\cdot 7^2\cdot 11\cdot 13\cdot 41\cdot 461}{3^2\cdot 7^2\cdot 11\cdot 13\cdot 19403}\) & \smallskip\centering \(\mbox{XII.}\,\amn{3^2\cdot 5\cdot 7\cdot 53\cdot 1889}{3^2\cdot 5\cdot 7\cdot 102059.}\)      
    \end{tabular}

    \endgroup
  \end{center}
\end{Scholium}

\begin{Problem}[2]
  {\em To find amicable numbers of the second form \(apq\),
    \(ars\) by supposing \(p\), \(q\), \(r\), \(s\) to be prime
    numbers and the common factor \(a\) given.}
\end{Problem}

\begin{SolutionUnnum}
  Since the common factor \(a\) is given, we seek from this the
  value of the fraction \(\dfrac{b}{c}=\dfrac{a}{2a-\int a}\) in
  lowest terms and hence we get \(a:\int a=b:2b-c\). Then, since
  it must be that \(\int p\cdot \int q=\int r\cdot \int s\), or
  \((p+1)(q+1)=(r+1)(s+1)\), suppose that each side is
  \(=\alpha\beta xy\) and take
  \[p=\alpha x-1,\quad q=\beta y-1,\quad r=\beta x-1,\quad
    s=\alpha y-1,\] where it is clear that these numbers \(\alpha\),
  \(\beta\), \(x\), \(y\) must be such that \(p\), \(q\),
  \(r\), \(s\) become prime numbers, and the amicable numbers will be
  \[a(\alpha x-1)(\beta y-1)\et a(\beta x-1)(\alpha y-1).\]
  Indeed, moreover, by the nature of amicable numbers we must
  have
  \[\alpha\beta xy\int a = a(\alpha x-1)(\beta y-1)+a(\beta
    x-1)(\alpha y-1)\] or, because \(\int a:a=2b-c:b\),
  \[\left.\begin{array}{c}2b\alpha\beta xy\\
            -c\alpha\beta xy\end{array}\right\} =
        \left\{\begin{array}{c}2b\alpha\beta xy-b\alpha x-b\beta
                 y+2b\\ -b\beta x-b\alpha y\end{array}\right.\] or\runon
  \[c\alpha\beta xy=b(\alpha+\beta)(x+y)-2b.\] Whence we
  get\footnote{Multiply by \(c\alpha\beta\) and add
    \(bb(\alpha+\beta)^2\) to both sides. \JE}
  \[\begin{array}{c}cc\alpha^2\beta^2 xy -
      bc\alpha\beta(\alpha+\beta)x + bb(\alpha+\beta)^2\\
      -bc\alpha\beta(\alpha+\beta)y\end{array}
    \begin{array}{c}=-2bc\alpha\beta +
      bb(\alpha+\beta)^2.\\\,\end{array}\] Wherefore this
  equation must be satisfied
  \[(c\alpha\beta x-b(\alpha+\beta))(c\alpha\beta
    y-b(\alpha+\beta))=bb(\alpha+\beta)^2-2bc\alpha\beta.\]
  Therefore, in all cases, the number
  \(bb(\alpha+\beta)^2-2bc\alpha\beta\) must be resolved into
  two factors, say \(PQ\), such that, by putting
  \[x=\frac{P+b(\alpha+\beta)}{c\alpha\beta}\et y
    =\frac{Q+b(\alpha+\beta)}{c\alpha\beta}\] these numbers
  \(x\) and \(y\) not only become integers, but also
  \(\alpha x-1\), \(\beta y-1\), \(\beta x-1\), and
  \(\alpha y-1\) become prime numbers. Therefore we
  get\footnote{Typo: the denominators of \(r\) and \(s\) are
    switched in the original. \JE}
  \begin{align*}
    p&=\frac{P+b\alpha+(b-c)\beta}{c\beta},
    & q= \frac{Q+b\beta+(b-c)\alpha}{c\alpha},\\
    r&=\frac{P+b\beta+(b-c)\alpha}{c\alpha},
    & s= \frac{Q+b\alpha+(b-c)\beta}{c\beta}.
  \end{align*}

  \par Therefore for any given value of \(a\), from which we
  obtain \(\dfrac{b}{c}=\dfrac{a}{2a-\int a}\), it is to be
  considered, whether, with the numbers \(\alpha\) and \(\beta\)
  assumed, a factorisation
  \[bb(\alpha+\beta)^2-2bc\alpha\beta = PQ\] may be made so that
  the values as given above for \(p\), \(q\), \(r\), and \(s\)
  become prime numbers and indeed such that the common factor
  involves none of them. And whenever this condition can be
  satisfied, the amicable numbers will be \(apq\) and \(ars\).
\end{SolutionUnnum}

\begin{Corollary}
  Simpler numbers are put for \(\alpha\) and \(\beta\), and
  since it cannot be that \(\alpha=\beta\), the following cases
  arise from this:
  \begin{itemize}
  \item[I.] Let \(\alpha=1\), \(\beta=2\); we get \(PQ=9bb-4bc\) and
    \begin{alignat*}{3}
      p&=\frac{P+3b-2c}{2c},&\quad& & q&= \frac{Q+3b-c}{c},\\
      r&=\frac{P+3b-c}{c}, &\quad& & s&= \frac{Q+3b-2c}{2c}.
    \end{alignat*}
    
  \item[II.] Let \(\alpha=1\), \(\beta=3\); we get \(PQ=16bb-6bc\) and
    \begin{alignat*}{3}
      p&=\frac{P+4b-3c}{3c},&\quad& & q&= \frac{Q+4b-c}{c},\\
      r&=\frac{P+4b-c}{c},&\quad& & s&= \frac{Q+4b-3c}{3c}.
    \end{alignat*}
    
  \item[III.] Let \(\alpha=2\), \(\beta=3\); we get \(PQ=25bb-12bc\) and
    \begin{alignat*}{3}
      p&=\frac{P+5b}{3c}-1,&\quad& & q&= \frac{Q+5b}{2c}-1,\\
      r&=\frac{P+5b}{2c}-1,&\quad& & s&= \frac{Q+5b}{3c}-1.
    \end{alignat*}
    
  \item[IV.] Let \(\alpha=1\), \(\beta=4\); we get \(PQ=25bb-8bc\) and
    \begin{alignat*}{3}
      p&=\frac{P+5b}{4c}-1,&\quad& & q&= \frac{Q+5b}{c}-1,\\
      r&=\frac{P+5b}{c}-1,&\quad& & s&= \frac{Q+5b}{4c}-1.
    \end{alignat*}

  \item[V.] Let \(\alpha=3\), \(\beta=4\); we get \(PQ=49bb-24bc\) and
    \begin{alignat*}{3}
      p&=\frac{P+7b}{4c}-1,&\quad& & q&= \frac{Q+7b}{3c}-1,\\
      r&=\frac{P+7b}{3c}-1,&\quad& & s&= \frac{Q+7b}{4c}-1.
    \end{alignat*}

  \item[VI.] Let \(\alpha=1\), \(\beta=5\); we get \(PQ=36bb-10bc\) and
    \begin{alignat*}{3}
      p&=\frac{P+6b}{5c}-1,&\quad& & q&= \frac{Q+6b}{c}-1,\\
      r&=\frac{P+6b}{c}-1,&\quad& & s&= \frac{Q+6b}{5c}-1.
    \end{alignat*}

  \item[VII.] Let \(\alpha=2\), \(\beta=5\); we get \(PQ=49bb-20bc\) and
    \begin{alignat*}{3}
      p&=\frac{P+7b}{5c}-1,&\quad& & q&= \frac{Q+7b}{2c}-1,\\
      r&=\frac{P+7b}{2c}-1,&\quad& & s&= \frac{Q+7b}{5c}-1.
    \end{alignat*}

  \item[VIII.] Let \(\alpha=3\), \(\beta=5\); we get \(PQ=64bb-30bc\) and
    \begin{alignat*}{3}
      p&=\frac{P+8b}{5c}-1,&\quad& & q&= \frac{Q+8b}{3c}-1,\\
      r&=\frac{P+8b}{3c}-1,&\quad& & s&= \frac{Q+8b}{5c}-1.
    \end{alignat*}

  \item[IX.] Let \(\alpha=4\), \(\beta=5\); we get \(PQ=81bb-40bc\) and
    \begin{alignat*}{3}
      p&=\frac{P+9b}{5c}-1,&\quad& & q&= \frac{Q+9b}{4c}-1,\\
      r&=\frac{P+9b}{4c}-1,&\quad& & s&= \frac{Q+9b}{5c}-1.
    \end{alignat*}

  \item[X.] Let \(\alpha=1\), \(\beta=6\); we get \(PQ=49bb-12bc\) and
    \begin{alignat*}{3}
      p&=\frac{P+7b}{6c}-1,&\quad& & q&= \frac{Q+7b}{c}-1,\\
      r&=\frac{P+7b}{c}-1,&\quad& & s&= \frac{Q+7b}{6c}-1.
    \end{alignat*}

  \item[XI.] Let \(\alpha=5\), \(\beta=6\); we get \(PQ=121bb-60bc\) and
    \begin{alignat*}{3}
      p&=\frac{P+11b}{6c}-1,&\quad& & q&= \frac{Q+11b}{5c}-1,\\
      r&=\frac{P+11b}{5c}-1,&\quad& & s&= \frac{Q+11b}{6c}-1.
    \end{alignat*}
  \end{itemize}
  Therefore in line with these cases, I will work through the
  values of \(a\) already used earlier, because in comparison
  with others they seemed suitable for finding amicable numbers;
  but from these I will choose chiefly the ones which, having
  been worked out, lead to amicable numbers.
\end{Corollary}

\newpage

\begin{Example}[1]
  Let \(a=2^2\); we get \(b=4\) and \(c=1\). We take the second
  case, where \(\alpha=1\), \(\beta=3\), so that the amicable
  numbers would be \(2^2pq\) and \(2^2rs\), and we must get
  \begin{gather*}
    PQ=16\cdot 16-6\cdot 4=232\tx{and also}\\
    p=\frac{P+16}{3}-1,\quad q=Q+16-1,\quad r=P+16-1,\quad
    s=\frac{Q+16}{3}-1.
  \end{gather*}
  Therefore the factors of the number \(232\) must have the
  property that, when they are increased by \(16\) they become
  divisible by \(3\).

  \bgroup
  \def\arraystretch{1.1}
  \begin{center}
    \begin{tabular}{r@{}rp{8cm}}
      \(P=\)&\(2\) & \multirow{7}{8cm}{\begin{tabular}{p{8cm}@{}}
                                    \qquad 
                                    No other resolution
                                    succeeds; indeed if we were
                                    to put\footnotemark \(P=8\), \(Q\) would
                                    become an odd number, and
                                    therefore neither \(q\) nor
                                    \(s\) could be prime
                                    numbers. Hence, therefore, we obtain
                                    the amicable numbers
                                    \[\amn{2^2\cdot 5 \cdot 131}{2^2\cdot 17\cdot 43.}\]\end{tabular}}\\
      \(Q=\)&\(\,116\) & \\
      \(P+16=\)&\(\,18\) & \\
      \(Q+16=\)&\(\,132\) &\\
      \cline{1-2}
      \(p=\)&\(\,5\) & \\
      \(q=\)&\(\,131\) & \\
      \(r=\)&\(\,17\) & \\
      \(s=\)&\(\,43\) &
    \end{tabular}
  \end{center}
  \egroup
  \footnotetext{There is a typo in the Opera Omnia edition here,
    where \(p=8\) is written instead of \(P=8\). \JE}
\end{Example}

\begin{Example}[2]
  If \(\alpha=1\) and \(\beta=3\) and \(a\) is a higher power of
  two, we do not succeed in finding amicable numbers until we
  reach \(a=2^8\). But then we get \(b=2^8\) and \(c=1\)
  and also
  \begin{gather*}
    PQ=16\cdot 2^{16}-6\cdot
    2^8=2^9\left(2^{11}-3\right)=512\cdot 2045=512\cdot 5\cdot
    409,\\
    p=\frac{P+1024}{3}-1,\quad q=Q+1024-1,\quad r=P+1024-1,\quad s=\frac{Q+1024}{3}-1,
  \end{gather*}
  whence the factors \(P\) and \(Q\) must have the
  property\footnote{Indeed, \(1020\) is divisible by \(3\), so
    we only need \(P+4\) and \(Q+4\) (or indeed \(P+1\) and
    \(Q+1\)) to be divisible by \(3\). \JE} that, when increased
  by \(4\) they become divisible by \(3\), or, if the quotient
  becomes even, divisible by \(6\).

  \bgroup
  \def\arraystretch{1.1}
  \begin{center}
    \begin{tabular}{>{\raggedleft\arraybackslash}p{1.9cm}|
      >{\raggedleft\arraybackslash}p{0.8cm}@{}p{0.01cm}|
      >{\raggedleft\arraybackslash}p{0.8cm}@{}p{0.01cm}|
      >{\raggedleft\arraybackslash}p{0.8cm}@{}p{0.01cm}|
      >{\raggedleft\arraybackslash}p{0.8cm}@{}p{0.01cm}|
      >{\raggedleft\arraybackslash}p{1cm}@{}p{0.01cm}|
      >{\raggedleft\arraybackslash}p{1cm}@{}p{0.01cm}|
      >{\raggedleft\arraybackslash}p{1cm}@{}p{0.01cm}|
      >{\raggedleft\arraybackslash}p{1cm}@{}p{0.01cm}|}
      \(P=\) & \(2\) && \(8\) && \(20\) && \(32\) && \(80\) && \(128\) && \(320\) && \(1280\)&\\
      \(Q=\) & \(\cdots\) && \(\cdots\) && \(\cdots\) && \(\cdots\) && \(13088\) && \(8180\) && \(\cdots\) && \(\cdots\)& \\
      \(P+1024=\) & \(1026\) && \(1032\) && \(1044\) && \(1056\) && \(1104\) && \(1152\) && \(1344\) && \(2304\)& \\
      \(Q+1024=\) & \(\cdots\) && \(\cdots\) && \(\cdots\) && \(\cdots\) && \(14112\) && \(9204\) && \(\cdots\) && \(\cdots\) & \\
      \(p=\) & \(341\)&\(^*\) & \(343\)&\(^*\) & \(347\) && \(\cdots\) && \(367\) && \(383\) && \(447\)&\(^*\) & \(767\)&\(^*\)\\
      \(q=\) & \(\cdots\) && \(\cdots\) && \(\cdots\) && \(\cdots\) && \(14111\)&\(^*\) & \(9203\) && \(\cdots\) && \(\cdots\)& \\
      \(r=\) & \(1025\)&\(^*\) & \(\cdots\) && \(1043\)&\(^*\) & \(1055\)&\(^*\) & \(1103\) && \(1151\) && \(1343\)&\(^*\) & \(2303\)&\(^\dagger\)\\
      \(s=\) & \(\cdots\) &&  \(\cdots\) &&  \(\cdots\) && \(\cdots\) && \(4703\) && \(3067\) && \(\cdots\) && \(\cdots\) &
    \end{tabular}
  \end{center}
  \egroup

  \medskip
  
  \par Therefore the amicable numbers will be \(\amn{2^8\cdot 383\cdot 9203}{2^8\cdot 1151\cdot 3067.}\)
\end{Example}
  
\begin{Example}[3]
  Let \(\alpha=2\) and \(\beta=3\) and take \(a=3^2\cdot 5\cdot
  13\), so that \(b=15\) and \(c=2\); we get
  \begin{gather*}
    PQ=25\cdot 225-12\cdot 30=3^4\cdot 5\cdot 13,\\
    p=\frac{P+75}{6}-1,\quad q=\frac{Q+75}{4}-1,\quad r=\frac{P+75}{4}-1,\quad
    s=\frac{Q+75}{6}-1,
  \end{gather*}
  whence the factors \(PQ\) must be such that, when increased by
  \(3\), they become divisible by\footnote{At first sight, it
    appears that it is enough for \(P+3\) and \(Q+3\) to be
    divisible by \(12\), since that makes \(p\), \(q\), \(r\),
    and \(s\) into integers. However, if for example
    \(P=24k+12\) then \(\frac{P+75}{4}-1=6k+20\), which is even,
    and hence not prime. So we may as well assume \(P+3\) and
    \(Q+3\) to be divisible by \(24\). \JE} \(24\).

  \bgroup
  \def\arraystretch{1.1}
  \begin{center}
    \begin{tabular}{>{\raggedleft\arraybackslash}p{1.5cm}
      >{\raggedleft\arraybackslash}p{0.8cm}@{}p{0.01cm}|
      p{7cm}}
      \(P=\)&\(45\) && \multirow{8}{6cm}{\begin{tabular}{p{6cm}@{}}
                                     \qquad Other
                                     resolutions do not take
                                     place here; whence appear the amicable numbers:
                                     \[\amn{3^2\cdot 5 \cdot 13\cdot 19\cdot 47}{3^2\cdot 5\cdot 13\cdot 29\cdot 31.}\]\end{tabular}}\\
      \(Q=\)&\(117\) && \\
      \(P+75=\)&\(120\) && \\
      \(Q+75=\)&\(192\) &&\\
      \(p=\)&\(19\) && \\
      \(q=\)&\(47\) && \\
      \(r=\)&\(29\) && \\
      \(s=\)&\(31\) &&
    \end{tabular}
  \end{center}
  \egroup
\end{Example}

\begin{Example}[4\footnotemark]
  \footnotetext{Mislabelled Example 3 in original. \JE} Let
  \(\alpha=1\) and \(\beta=4\) and take \(a=3^3\cdot 5\), so
  that\footnote{Typo: \(p=9\) in original. \JE} \(b=9\) and
  \(c=2\); we get
  \begin{gather*}
    PQ=25\cdot 81-8\cdot 18=9\cdot 11\cdot 19,\et\\
    p=\frac{P+45}{8}-1,\quad q=\frac{Q+45}{2}-1,\quad
    r=\frac{P+45}{2}-1,\quad s=\frac{Q+45}{8}-1,
  \end{gather*}
  whence \(P\) and \(Q\) must be numbers which\footnote{Again,
    \(40\) is already divisible by \(8\). \JE}, when increased
  by \(5\) become divisible by \(8\).

  \bgroup
  \def\arraystretch{1.1}
  \begin{center}
    \begin{tabular}{>{\raggedleft\arraybackslash}p{1.9cm}
      >{\raggedleft\arraybackslash}p{0.8cm}@{}p{0.3cm}|
      >{\raggedleft\arraybackslash}p{0.8cm}@{}p{0.01cm}|
      p{7cm}}
      \(P=\) & \(3\) && \(19\) && \multirow{8}{7cm}{\begin{tabular}{p{7cm}@{}}
                                    \qquad Therefore the
                                                    amicable
                                                    numbers
                                                    arising in
                                                    this way are
                                    \[\amn{3^3\cdot 5\cdot 7\cdot 71}{3^3\cdot 5\cdot 31\cdot 17.}\]\end{tabular}}\\
      \(Q=\) & \(627\) && \(99\) && \\
      \(P+45=\) & \(48\)&)\footnotemark & \(64\) && \\
      \(Q+45=\) & \(672\) && \(144\) && \\
      \(p=\) & \(5\) && \(7\) && \\
      \(q=\) & \(335\)&\(^*\) & \(71\) && \\
      \(r=\) & \(23\) && \(31\) &&\\
      \(s=\) & \(83\) && \(17\) &&
    \end{tabular}
  \end{center}
  \egroup
  \footnotetext{\(48\) is given as \(28\) in the original.}
\end{Example}

\begin{Scholium}
  However, these operations are too undependable and
  usually more are undertaken in vain before amicable numbers
  show up. Moreover, the work would be extremely lengthy if, for
  each value of \(a\) which I showed above, we wanted to run
  through each case of the letters \(\alpha\) and \(\beta\), as
  it happens exceedingly rarely that the four resulting numbers
  for \(p\), \(q\), \(r\), and \(s\) become simultaneously prime.
  Indeed even then the discovery of amicable numbers by
  determination of the ratio of \(\alpha\) and \(\beta\) is
  too constrained, and cases exist of amicable numbers
  where the ratio \(\alpha:\beta\) is so complicated that it
  could not have been chosen by probable reason; of this sort
  are the amicable numbers\footnote{That these numbers
    \(2^4\cdot 19\cdot 8563\) and \(2^4\cdot 83\cdot 2039\) are
    not amicable was observed by K. Hunrath
    \cite{EnestromKM}. \FR} \(2^4\cdot 19\cdot 8563\) and
  \(2^4\cdot 83\cdot 2039\), for the discovery of which, by this
  method, one would need to assume the ratio to be \(5:21\) or
  \(1:102\). Because of this, I will not linger longer on this
  exceedingly sterile and laborious method, but instead I will
  explain another way by which one can more easily and freely
  investigate amicable numbers, both of this second form and of
  other more composite forms; a way which is similar to what
  went before in that it is solved by finding only three prime
  numbers.
\end{Scholium}

\begin{Problem}[3\footnotemark]
  \footnotetext{Mislabelled ``Problem 2'' in original. \JE} {\em
    To find amicable numbers of the form \(apq\) and \(afr\),
    where \(p\), \(q\), and \(r\) are prime numbers, \(f\) is
    either prime or composite, and as before the common factor
    \(a\) is given.}
\end{Problem}

\begin{SolutionUnnum}
  Again, from knowing the common factor \(a\), values of \(b\)
  and \(c\) are sought so that
  \(\dfrac{b}{c}=\dfrac{a}{2a-\int a}\); let the divisor sum of
  the number \(f\) be \(\int f = gh\). Therefore, since it is
  first required that \(\int p\cdot \int q=\int f\cdot \int r\),
  we get \((p+1)(q+1)=gh(r+1)\). We put \(r+1=xy\), \(p+1=hx\),
  and \(q+1=gy\) and it will be necessary that these three
  numbers be prime, namely: \(p=hx-1\), \(q=gy-1\) and
  \(r=xy-1\). Then it is necessary that
  \begin{gather*}
    \int apq = ghxy\int a = a(hx-1)(gy-1)+af(xy-1)\\
    =a((gh+f)xy-hx-gy+1-f)
  \end{gather*}
  or\runon
  \[2bghxy-cghxy=b(gh+f)xy-bhx-bgy+b(1-f)\]
  or\runon
  \[(bf-bgh+cgh)xy-bhx-bgy=b(f-1).\]
  \par For the sake of brevity, we put
  \[bf-bgh+cgh=e;\]
  we get \(eexy-ebhx-ebgy=eb(f-1)\) or
  \[(ex-bg)(ey-bh) = bbgh+be(f-1).\]
  \par Therefore the number \(bbgh+be(f-1)\) must be resolved
  into two factors, \(P\) and \(Q\), such that
  \[x=\frac{P+bg}{e}\et y=\frac{Q+bh}{e}\] become integers, and
  moreover \(hx-1\), \(gy-1\), and \(xy-1\) become prime
  numbers. Whenever this condition can be satisfied, the
  amicable numbers will be \[a(hx-1)(gy-1)\et af(xy-1).\] It is
  important to note that none of these prime
  numbers: \[hx-1,\quad gy-1,\quad xy-1,\] nor any factor of
  \(f\) is allowed to be a divisor of \(a\), and that \(f\) and
  \(xy-1\) must be relatively prime.
\end{SolutionUnnum}

\begin{Corollary}[1]
  Let \(f\) be a prime number, as the second form for amicable
  numbers postulates; we get \(f+1=gh\), and therefore
  \(f=gh-1\). Therefore in this case we get \(e=cgh-b\) and
  \(PQ=bbgh+be(gh-2)\) or
  \[PQ=bcgghh-2bcgh+2bb.\] Whence numbers \(x\) and \(y\) must
  be sought, possessing the properties mentioned above, such
  that \[x=\frac{P+bg}{e}\et y=\frac{Q+bh}{e}.\]
\end{Corollary}

\begin{Corollary}[2]
  Therefore it will be convenient to use these formulas in the
  following way: for \(a\) is substituted successively one or
  another value from those which I have shown above, and for each
  letter \(f\) various numbers, either prime or composite, are
  substituted, indeed those which seem suitable for finding
  amicable numbers.
\end{Corollary}

\begin{Case}[1]
  Let \(a=4\) (indeed, from the value \(a=2\) I have observed no
  amicable numbers to be obtained) and we get \(b=4\) and
  \(c=1\). Then by supposing \(4pq\) and \(4fr\) to be amicable
  numbers let \(\int f=gh\) and \(e=4f-3gh\). Then by resolving,
  factors \(P\) and \(Q\) are sought such that
  \[PQ=16gh+4e(f-1).\]
  \par And hence we extract integers \(x\) and
  \(y\), such that
  \[x=\frac{P+4g}{e}\et y=\frac{Q+4h}{e},\] and from this are
  derived the values of the letters \(p=hx-1\), \(q=gy-1\), and
  \(r=xy-1\); if these are prime numbers, \(4pq\) and \(4fr\)
  will be amicable numbers.
\end{Case}

\begin{Example}[1]
  Let \(f=3\); we get \(\int f=gh=4\) and hence \(e=12-12=0\),
  whence it is clear that nothing is obtained from this
  hypothesis.
\end{Example}

\begin{Example}[2]
  Let \(f=5\); we get \(\int f=gh=6\), \(e=20-18=2\) and
  \[PQ=16\cdot 6+8\cdot 4=128.\] Now from \(gh=6\) suppose first
  that
  \(g=2\) and \(h=3\), and we get
  \[x=\frac{P+8}{2}\et y=\frac{Q+12}{2}.\] Wherefore we will
  have the following resolutions:

  \bgroup
  \def\arraystretch{1.1}
  \begin{center}
    \begin{tabular}{>{\raggedleft\arraybackslash}p{2.5cm}|
      >{\raggedleft\arraybackslash}p{0.5cm}@{}p{0.4cm}|
      >{\raggedleft\arraybackslash}p{0.5cm}@{}p{0.01cm}|
      >{\raggedleft\arraybackslash}p{0.6cm}@{}p{0.01cm}|
      >{\raggedleft\arraybackslash}p{0.6cm}@{}p{0.01cm}|
      >{\raggedleft\arraybackslash}p{0.6cm}@{}p{0.01cm}|
      >{\raggedleft\arraybackslash}p{0.6cm}@{}p{0.01cm}|
      p{3cm}}
      \(P=\) & \(2\) && \(4\) && \(8\) && \(16\) && \(32\) && \(64\)
      && \multirow{7}{3cm}{\begin{tabular}{p{3cm}@{}}
                            \qquad Therefore these amicable numbers appear:\\
                            \(\amn{4\cdot 17 \cdot 43}{4\cdot 5\cdot 131}\) and
                            \(\amn{4\cdot 13 \cdot 107}{4\cdot 5\cdot 251}\)\end{tabular}}\\
      \(Q=\) & \(64\) && \(32\) && \(16\) && \(8\) && \(4\) && \(2\)
      && \\
      \(x=\) & \(5\) && \(6\) && \(8\) && \(12\) && \(20\) && \(36\)
      && \\
      \(y=\) & \(38\) && \(22\) && \(14\) && \(10\) && \(8\) && \(7\)
      && \\
      \(p=3x-1=\) & \(14\)&\(^*\))\footnotemark & \(17\) && \(23\) && \(35\)&\(^*\) & \(59\) && \(107\)
      && \\
      \(q=2y-1=\) & \(\cdots\) && \(43\) && \(27\)&\(^*\) & \(19\) && \(15\)&\(^*\) & \(13\)
      && \\
      \(r=xy-1=\) & \(\cdots\) && \(131\) && \(111\)&\(^*\) & \(119\)&\(^*\) & \(159\)&\(^\dagger\) & \(251\)
      && \\
    \end{tabular}
  \footnotetext{Given as \(19^*\) in original. \JE}
  \end{center}
  \egroup

  Secondly, suppose \(g=1\), \(h=6\) and we get
  \[x=\frac{P+4}{2}\et y=\frac{Q+24}{2}.\]

  \bgroup
  \def\arraystretch{1.1}
  \begin{center}
    \begin{tabular}{>{\raggedleft\arraybackslash}p{2.5cm}|
      >{\raggedleft\arraybackslash}p{0.5cm}@{}p{0.3cm}|
      >{\raggedleft\arraybackslash}p{0.6cm}@{}p{0.01cm}|
      >{\raggedleft\arraybackslash}p{0.6cm}@{}p{0.01cm}|
      >{\raggedleft\arraybackslash}p{0.6cm}@{}p{0.01cm}|
      >{\raggedleft\arraybackslash}p{0.6cm}@{}p{0.01cm}|
      >{\raggedleft\arraybackslash}p{0.6cm}@{}p{0.01cm}|
      p{3cm}}
      \(P=\) & \(2\) && \(4\) && \(8\) && \(16\) && \(32\) && \(64\)
      && \multirow{7}{3cm}{\begin{tabular}{p{3cm}@{}}
                            \qquad Therefore the same
                            two amicable numbers appear as
                            before.\end{tabular}}\\
      \(Q=\) & \(64\) && \(32\) && \(16\) && \(8\) && \(4\) && \(2\)
      && \\
      \(x=\) & \(3\) && \(4\) && \(6\) && \(10\) && \(18\) && \(34\)
      && \\
      \(y=\) & \(44\) && \(28\) && \(20\) && \(16\) && \(14\) && \(13\)
      && \\
      \(p=6x-1=\) & \(17\) &)\footnotemark& \(23\) && \(35\)&\(^*\) & \(59\) && \(107\) && \(203\)&\(^*\)
      & \\
      \(q=1y-1=\) & \(43\) && \(27\)&\(^*\) & \(19\) && \(15\)&\(^*\) & \(13\) && \(12\)&\(^*\)
      & \\
      \(r=xy-1=\) & \(131\) && \(111\)&\(^*\) & \(119\)&\(^*\) & \(159\)&\(^\dagger\) & \(251\) && \(441\)&\(^*\)
      & \\
    \end{tabular}
    \footnotetext{Erroneously asterisked in the original. \JE}
  \end{center}
  \egroup

  Therefore the amicable numbers are
  \[\amn{4\cdot 17 \cdot 43}{4\cdot 5\cdot 131}\et
    \amn{4\cdot 13 \cdot 107}{4\cdot 5\cdot 251.}\]
\end{Example}

\begin{Example}[3]
  Let \(f=7\); we get \(\int f=gh=8\), \(e=28-24=4\) and
  \[PQ=16\cdot 8+16\cdot 6=224.\]
  Therefore first let \(g=2\), \(h=4\); we get
  \[x=\frac{P+8}{4},\quad y=\frac{Q+16}{4},\quad p=4x-1,\quad
    q=2y-1,\quad r=xy-1.\]

  \bgroup
  \def\arraystretch{1.1}
  \begin{center}
    \begin{tabular}{>{\raggedleft\arraybackslash}p{1cm}|
      >{\raggedleft\arraybackslash}p{0.6cm}@{}p{0.01cm}|
      >{\raggedleft\arraybackslash}p{0.5cm}@{}p{0.3cm}|
      >{\raggedleft\arraybackslash}p{0.6cm}@{}p{0.01cm}|
      >{\raggedleft\arraybackslash}p{0.6cm}@{}p{0.01cm}|}
      \(P\) & \(4\) && \(8\) && \(28\) && \(56\)& \\
      \(Q\) & \(56\) && \(28\) && \(8\) && \(4\)& \\
      \(x\) & \(3\) && \(4\) && \(9\) && \(16\) &\\
      \(y\) & \(18\) && \(11\) && \(6\) && \(5\)&\\
      \(4x-1\) & \(11\) && \(15\)&\(^*\) & \(35\)&\(^*\) & \(63\)&\(^*\)\\
      \(2y-1\) & \(35\)&\(^*\) & \(21\)&\(^*\) & \(11\) && \(9\)&\(^*\)\\
      \(xy-1\) & \(53\) && \(43\) &)\footnotemark& \(53\) && \(79\)&\\
    \end{tabular}
    \footnotetext{Typo: \(42\) in original. \JE}
  \end{center}
  \egroup
  
  Second, let \(g=1\), \(h=8\); we get
  \[x=\frac{P+4}{4},\quad y=\frac{Q+32}{4},\quad p=8x-1,\quad
    q=y-1,\quad r=xy-1.\]

  \bgroup
  \def\arraystretch{1.1}
  \begin{center}
    \begin{tabular}{>{\raggedleft\arraybackslash}p{1cm}|
      >{\raggedleft\arraybackslash}p{0.6cm}@{}p{0.01cm}|
      >{\raggedleft\arraybackslash}p{0.5cm}@{}p{0.3cm}|
      >{\raggedleft\arraybackslash}p{0.6cm}@{}p{0.01cm}|
      >{\raggedleft\arraybackslash}p{0.6cm}@{}p{0.01cm}|}
      \(P\) & \(4\) && \(8\) && \(28\) && \(56\) &\\
      \(Q\) & \(56\) && \(28\) && \(8\) && \(4\) &\\
      \(x\) & \(2\) && \(3\) && \(8\) && \(15\) &\\
      \(y\) & \(22\) && \(15\) && \(10\) && \(9\)&\\
      \(8x-1\) & \(15\)&\(^*\) & \(23\) && \(63\)&\(^*\) & \(119\)&\(^*\)\\
      \(y-1\) & \(21\)&\(^\dagger\) & \(14\)&\(^*\) & \(9\)&\(^\dagger\) & \(8\)&\(^\dagger\)\\
      \(xy-1\) & \(43\) && \(44\)&\(^\dagger\) & \(79\) && \(134\)&\(^*\)\\
    \end{tabular}
  \end{center}
  \egroup
  
  Hence, therefore, no amicable numbers appear.\footnote{Recall
    from Footnote \ref{fn:daggers} that \(\dagger\) in these
    tables stands for an asterisk that Euler omitted. \JE}
\end{Example}

\begin{Example}[4]
  Let \(f=11\); we get \(gh=12\), \(e=8\),
  \(PQ=16\cdot 12+32\cdot 10=512\), in other words we get
  \((8x-4g)(8y-4h)=512\), which equation reduces to
  \((2x-g)(2y-h)=32\); by this resolution we get \(p=hx-1\),
  \(q=gy-1\), and \(r=xy-1\). On the other hand, if we suppose \(g=1\),
  \(h=12\), or \(g=2\), \(h=6\), or \(g=3\), \(h=4\), no prime
  numbers appear for \(p\), \(q\), and \(r\).
\end{Example}

\begin{Example}[5]\label{exm:commentworthy}
  Let \(f=13\); we get \(gh=14\), \(e=10\),
  \(PQ=224+40\cdot 12=704\) and \((10x-4g)(10y-4h)=704\), which
  reduces to \((5x-2g)(5y-2h)=176\). And hence no amicable
  numbers are obtained other than
  \[\amn{4\cdot 5\cdot 251}{4\cdot 13 \cdot 107,}\]
  which was already found earlier (\S 78).  At the same time it
  is already clear that even if bigger prime numbers are used
  for \(f\) no new amicable numbers appear, because either \(p\)
  or \(q\) will end up being smaller than the value which could
  have been assumed for \(f\).\footnote{This comment deserves
    amplification. See Appendix \ref{app:discussion}. \JE}
\end{Example}

\begin{Example}[6]
  Let \(f=5\cdot 13\); we get \(gh=6\cdot 14=84\), \(e=8\),
  \(PQ=16\cdot 84+32\cdot 64=64\cdot 53\) and
  \((8x-4g)(8y-4h)=64\cdot 53\) or \((2x-g)(2y-h)=4\cdot
  53\). And hence a solution is found in the prime numbers
  \(p=43\), \(q=2267\), and \(r=1187\); whence the amicable
  numbers will be
  \[\amn{4\cdot 43\cdot 2267}{4\cdot 5\cdot 13 \cdot 1187.}\]
\end{Example}

\newpage

\begin{Case}[2]
  Let \(a=2^3=8\); we get \(b=8\), \(c=1\); then by supposing
  \(8pq\) and \(8fr\) to be amicable numbers and \(\int f=gh\)
  we get \(e=8f-7gh\) and also
  \[(ex-8g)(ey-8h)=64gh+8e(f-1),\] whence the cases are to be
  determined in which the following numbers become prime
  \[p=hx-1,\quad q=gy-1\et r=xy-1.\]
\end{Case}

\begin{Example}[1]
  Let \(f=11\); we get \(gh=12\), \(e=4\), and also
  \[(4x-8g)(4y-8h)=64\cdot 12+32\cdot 10=64\cdot 17\]
  or \[(x-2g)(y-2h)=4\cdot 17=68.\] And hence no amicable
  numbers are discovered.
\end{Example}

\begin{Example}[2]
  Let \(f=13\); we get \(gh=14\), \(e=6\), and also
  \[(6x-8g)(6y-8h) = 64\cdot 14+48\cdot 12=64\cdot 23\]
  or\footnote{Typo: \((3x-4g)(hy-4h)\) in
    original. \JE}\runon \[(3x-4g)(3y-4h)=16\cdot 23;\] again,
  indeed, this hypothesis is of no use.
\end{Example}

\begin{Example}[3]
  Let \(f=17\); we get \(gh=18\), \(e=10\) and
  also\footnote{Typo: \(64\cdot 17+40\cdot 16\) in
    original. \JE}
  \[(10x-8g)(10y-8h)=64\cdot 18+80\cdot 16=64\cdot 38\]
  or
  \[(5x-4g)(5y-4h) = 32\cdot 19;\]
  and hence appear the amicable numbers
  \[\amn{8\cdot 23\cdot 59}{8\cdot 17\cdot 79.}\]
\end{Example}

\begin{Example}[4]
  More productive is the hypothesis \(f=11\cdot 23\); indeed
  a smaller composite value of \(f\) cannot be substituted; we
  get \(gh=12\cdot 24\), \(e=8\), whence
  \[(8x-8g)(8y-8h)=64\cdot 12\cdot 24+64\cdot 252\]
  or\runon
  \[(x-g)(y-h)=540.\] And hence the following amicable numbers
  are discovered
  \[\amn{8\cdot 383\cdot 1907}{8\cdot 11\cdot 23\cdot 2543}
    \quad
    \amn{8\cdot 467\cdot 1151}{8\cdot 11\cdot 23\cdot 1871}
    \quad
    \amn{8\cdot 647\cdot 719}{8\cdot 11\cdot 23\cdot 1619.}\]
  \par By taking such composite numbers for \(f\) many other amicable
  numbers are found besides.
\end{Example}

\begin{Scholium}
  The enormous number of combinations which appeared in this
  example gave me the leverage to reduce the solution to another
  more convenient form. Namely, since
  \[e=bf-(b-c)gh,\quad PQ = bbgh+be(f-1)=(ex-bg)(ey-bh),\] from
  the formula
  \[x=\frac{P+bg}{e}\et y=\frac{Q+bh}{e}\] we extract the values
  \[p=\frac{hP+bgh}{e}-1,\quad q=\frac{gQ+bgh}{e}-1,\quad
    r=\frac{PQ+b(hP+gQ)+bbgh}{ee}-1.\] Therefore because
  \(gh=\int f\), let
  \[e=bf-(b-c)\int f,\quad L=bb\int f+be(f-1)\et MN=L\int f;\]
  we get
  \[p=\frac{M+b\int f}{e}-1,\quad q=\frac{N+b\int f}{e}-1,\quad
    r=\frac{L+b(M+N)+bb\int f}{ee}-1,\] and now the question
  reduces to that of resolving the number \(L\int f\) into two
  factors \(M\) and \(N\), each of which, when increased by the
  quantity \(b\int f\) becomes divisible by \(e\), and such that
  when one is subtracted from the resulting quotients, they
  become prime. Lastly it is necessary that
  \(r+1=\dfrac{(p+1)(q+1)}{\int f}\) and \(r\) be a prime
  number. I will therefore illustrate this calculation in
  several cases.
\end{Scholium}

\begin{Case}[3]
  Let \(a=2^4=16\); we get \(b=16\), \(c=1\), and also
  \[e=16f-15\int f,\quad L=256\int f+16e(f-1)\et MN=L\int f.\]
  \par Therefore these numbers must be prime
  \[p=\frac{M+16\int f}{e}-1,\quad q=\frac{N+16\int
      f}{e}-1,\quad r=\frac{L+256\int f+16(M+N)}{ee}-1,\] and by
  finding them we will get the amicable numbers \(16pq\) and
  \(16fr\).
\end{Case}

\begin{Example}[1]
  Let \(f=17\); we get
  \begin{gather*}
    \int f=18,\quad e=2,\quad L=1024\cdot 5\et MN=1024\cdot
    5\cdot 18=2^{11}\cdot 3^2\cdot 5,\\
    p=\frac{M+288}{2}-1,\quad q=\frac{N+288}{2}-1,\quad
    r=\frac{512\cdot 19+16(M+N)}{4}-1;
  \end{gather*}
  or let \(M=2m\), \(N=2n\), so that \(mn=2^9\cdot 3^2\cdot 5\);
  we get
  \[p=m+143,\quad q=n+143,\et r=8(m+n)+2431,\] which three
  numbers must be prime so that \(16pq\) and
  \(16\cdot 17r\) are amicable numbers.

  \par But this succeeds in two
  ways, first if 
  \(m=24\), \(n=960\), and second, if \(m=96\) and \(n=240\);
  whence appear the amicable numbers:
  \[\amn{16\cdot 167\cdot 1103}{16\cdot 17\cdot 10303}
    \quad\amn{16\cdot 383\cdot 239}{16\cdot 17\cdot 5119.}\]
\end{Example}

\begin{Example}[2]
  Let \(f=19\); we get
  \[\int f = 20,\quad e=4,\quad L=128\cdot 49\et MN=512\cdot
    5\cdot 49=2^9\cdot 5\cdot 7^2.\] Therefore\footnote{Typo:
    \(r=\frac{128\cdot 59+16(N+N)}{16}-1\) in original. \JE}
  \[p=\frac{M+320}{4}-1,\quad q=\frac{N+320}{4}-1,\quad
    r=\frac{128\cdot 89+16(M+N)}{16}-1;\] or let \(M=4m\) and
  \(N=4n\), so that \(mn=32\cdot 5\cdot 49=2^5\cdot 5\cdot 7^2\);
  we get \[p=m+79,\quad q=n+79\et r=4(m+n)+711.\]

  \par Hence, if \(m=70\), \(n=112\), there appear the amicable
  numbers:
  \[\amn{16\cdot 149\cdot 191}{16\cdot 19\cdot 1439.}\]
\end{Example}

\begin{Example}[3]
  Let \(f=23\); we get
  \begin{gather*}
    \int f= 24, \quad e=8,\quad L=256\cdot 5\cdot 7\et
    MN=2048\cdot 3\cdot 5\cdot 7=2^{11}\cdot 3\cdot 5\cdot 7,\\
    p=\frac{M+16\cdot 24}{8}-1,\quad q=\frac{N+16\cdot
      24}{8}-1,\quad r=\frac{256\cdot 59+16(M+N)}{64}-1;
  \end{gather*}
  or let \(M=8m\), \(N=8n\), and \(mn=2^5\cdot 3\cdot 5\cdot 7\);
  we get
  \[p=m+47,\quad q=n+47\et r=2(m+n)+235.\]

  \par Hence three cases
  arise: \qquad
  \(\begin{cases}m=56\\n=60\end{cases}\quad
    \begin{cases}m=42\\n=80\end{cases}\quad
    \begin{cases}m=6\\n=560\end{cases}\)
  \newline and the amicable numbers are:
  \[\amn{16\cdot 103\cdot 107}{16\cdot 23\cdot 467}\quad
    \amn{16\cdot 89\cdot 127}{16\cdot 23\cdot 479}\quad
    \amn{16\cdot 53\cdot 607}{16\cdot 23\cdot 1367.}\]
\end{Example}

\begin{Example}[4]
  Let \(f=31\); we get\footnote{For completeness, note that
    \(e=16\) in this example. \JE}
  \begin{gather*}
    \int f= 32, \quad L=512\cdot 31\et
    MN=2^{14}\cdot 31,\\
    p=\frac{M+16\cdot 32}{16}-1,\quad q=\frac{N+16\cdot
      32}{16}-1,\quad r=\frac{16(M+N)+512\cdot 47}{256}-1.
  \end{gather*}

  \par Let therefore \(M=16m\), \(N=16n\), so that
  \(mn=2^6\cdot 31\); we get
  \[p=m+31,\quad q=n+31,\quad r=m+n+93.\]

  \par And hence no amicable numbers appear.
\end{Example}

\begin{Example}[5]
  Let \(f=47\), \(\int f=48\); we get
  \[e=32\et L=1024\cdot 5\cdot 7\et MN=2^{14}\cdot 3\cdot 5\cdot
    7,\] whence
  \[p=\frac{M+16\cdot 48}{32}-1,\quad q=\frac{N+16\cdot
      48}{32}-1,\quad r=\frac{16(M+N)+1024\cdot 47}{1024}-1.\]
  Let \(M=32m\) and \(N=32n\), so that \(mn=2^4\cdot 3\cdot 5\cdot 7\); we get
  \[p=m+23,\quad q=n+23,\quad r=\frac{1}{2}(m+n)+46.\] Therefore
  \(m+n\) must be an odd-times-even number\footnote{That is even
    but not divisible by \(4\). \JE}, so that
  \(\frac{1}{2}(m+n)\) becomes odd, which happens if either
  \(m\) or \(n\) is odd-times-even. Let \(m=30\), \(n=56\);
  the amicable numbers will be:
  \[\amn{16\cdot 53\cdot 79}{16\cdot 47\cdot 89.}\]
\end{Example}

\setcounter{eulercounter}{93}
\begingroup
\renewenvironment{Example}[1][]{\refstepcounter{eulercounter}\par\medskip
  \begin{center}\sc\large Exemplum #1\end{center}\medskip
  \theeulercounter(a)\footnote{In the first edition the number
    94 is used again by mistake. {\FR} So, like Rudio, we keep
    this mistake for backward compatibility of referencing. To
    be fair to the original, it does use different numerals
    (LXXXXIV and CIV)... \JE}.\,\,\,\,}{}

\begin{Example}[6]
  Let \(f=17\cdot 137\); we get
  \begin{gather*}
    \int f = 18\cdot 138=4\cdot 27\cdot 23=2484,\quad e=4,\\
    L = 256\cdot 2484+64\cdot 2328=512\cdot 3\cdot 7 \cdot 73\et MN=2048\cdot 81\cdot 7\cdot 23\cdot 73,\\
    p = \frac{M+16\cdot 2484}{4}-1,\quad q=\frac{N+16\cdot 2484}{4}-1,\\
    r=\frac{512\cdot 2775+16(M+N)}{16}-1.
  \end{gather*}
  Let \(M=4m\), \(N=4n\); we get \(mn=128\cdot 81\cdot 7\cdot 23\cdot 73\) and
  \[p=m+9935,\quad q=n+9935,\quad r=4(m+n)+88799.\]

  \par But this always produces a value of \(r\) bigger than
  \(100000\), so it is difficult to tell whether or not it is
  prime.
\end{Example}
\endgroup

\newpage

\begin{Example}[7]
  Let \(f=17\cdot 151\); we get\footnote{There is a typo in the
    Opera Omnia edition here: \(\int f\) is given as \(276\). In
    the original, \(1024\) is given as \(1084\) the second time
    it appears. \JE}
  \begin{gather*}
    \int f= 18\cdot 152=16\cdot 9\cdot 19=2736,\quad e=32\\
    L=1024\cdot 1967=1024\cdot 7\cdot 281\tx{and also} MN=2^{14}\cdot 9\cdot 7\cdot 19\cdot 281.
  \end{gather*}
  Let \(M=32m\), \(N=32n\); we get \(mn=16\cdot 9\cdot 7\cdot 19\cdot 281\) and
  \[p=m+1367,\quad q=n+1367,\quad r=\frac{1}{2}(m+n)+2650.\] Let
  \(m=2\mu\), \(n=8\nu\); we get
  \(\mu\nu=9\cdot 7\cdot 19\cdot 281\) and
  \[p=2\mu+1367,\quad q=8\nu+1367,\quad r=\mu+4\nu+2650.\]

  \par Hence, first, it is clear neither \(\mu\) nor \(\nu\) can
  be a number of the form \(3\alpha+2\); and then \(\mu\) cannot
  end in a \(9\) nor \(\nu\) in a \(1\); by these observations,
  only the following resolutions take place:

  \bgroup
  \def\arraystretch{1.2}
  \begin{center}
    \begin{tabular}{c|
      >{\centering\arraybackslash}p{1.2cm}|
      >{\centering\arraybackslash}p{1.2cm}|
      >{\centering\arraybackslash}p{1.2cm}|
      >{\centering\arraybackslash}p{1.2cm}|
      >{\centering\arraybackslash}p{1.2cm}|
      >{\centering\arraybackslash}p{1.3cm}|
      >{\centering\arraybackslash}p{1.5cm}|}
      & \(*\) & & & & \(*\) & \(*\) & \\
      \(\mu\) & \(3\cdot 281\) & \(7\cdot 19\) & \(21\cdot 281\)
              & \(21\) & \(63\cdot 281\) & \(3\) & \(1\)\\
      \(\nu\) & \(21\cdot 19\) & \(9\cdot 281\) & \(57\) &
                                                           \(57\cdot
                                                           281\)
                & \(19\) & \(399\cdot 281\) & \(1197\cdot 281\)
    \end{tabular}
  \end{center}
  \egroup
  
  \noindent of which those which are denoted with an asterisk
  are excluded for the reason that none of \(p\), \(q\) or \(r\)
  can be divisible by \(7\). The fourth resolution will give
  these amicable numbers
  \[\amn{16\cdot 1409\cdot 129503}{16\cdot 17\cdot 151\cdot
      66739,}\] only if this number
  \(129503\) is prime.\footnote{However
    \(129503=11\cdot 61\cdot 193\) so the corresponding numbers
    are not amicable. \FR}
\end{Example}

\begin{Example}[8]
  Let \(f=17\cdot 167\); we get
  \begin{gather*}
    \int f=18\cdot 168=16\cdot 27\cdot 7=3024,\quad e=64,\\
    L=2048\cdot 1797=2048\cdot 3\cdot 599\et MN=2^{15}\cdot
    3^4\cdot 7\cdot 599.
  \end{gather*}
  Let \(M=64m\), \(N=64n\); we get \(mn=2^3\cdot 3^4\cdot 7\cdot
  599\) and
  \[p=m+755,\quad q=n+755,\quad
    r=\frac{1}{4}(m+n)+\frac{2173}{2}.\]
  Let \(m=2\mu\), \(n=4\nu\); we get \(\mu\nu=3^4\cdot 7\cdot
  599\) and
  \[p=2\mu+755,\quad q=4\nu+755,\quad
    r=\nu+\frac{\mu+1}{2}+1086,\] where it is clear that we must
  have \(\mu=4\alpha-1\), lest \(r\) becomes an even
  number\footnote{Note that \(\mu\) and \(\nu\) are odd because
    their product is odd, so \((\mu+1)/2\) needs to be even to
    make \(r\) odd. \JE} nor \(\mu=3\alpha+2\) nor
  \(\nu=3\alpha+1\). From here appear the amicable numbers 
  \[\amn{16\cdot 809\cdot 51071}{16\cdot 17\cdot 167\cdot 13679.}\]
\end{Example}

\begin{Case}[4]
  Let either \(a=3^3\cdot 5\) or \(a=3^2\cdot 7\cdot 13\), so
  that \(b=9\), \(c=2\); we get\footnote{There is a typo in the
    Opera Omnia here: \(MM\) is written instead of \(MN\). \JE}
  \begin{gather*}
    e=9f-7\int f,\quad L=81\int f+9e(f-1)\et MN=L\int f,\\
    p=\frac{M+9\int f}{e}-1,\quad q=\frac{N+9\int f}{e}-1,\\
    r=\frac{9(M+N)+L+81\int f}{ee}-1;
  \end{gather*}
  if these \(p\), \(q\), \(r\) become prime, the amicable
  numbers will be \(\amn{apq}{afr.}\)
\end{Case}

\begin{Example}
  Let \(f=7\), \(\int f=8\); we get
  \begin{gather*}
    e=7,\quad L=2\cdot 27\cdot 19,\quad MN=16\cdot 27\cdot 19,\\
    p=\frac{M+72}{7}-1,\quad q=\frac{N+72}{7}-1,\quad
    r=\frac{9(M+N)+2\cdot 27\cdot 31}{49}-1.
  \end{gather*}
  Whence by putting \(M=54\), \(N=152\) the following amicable
  numbers arise:\footnote{Note that the other value of \(a\),
    \(3^2\cdot 7\cdot 13\), does not appear here because it
    shares a factor of \(7\) with \(f\). \JE}
  \[\amn{a\cdot 17\cdot 31}{a\cdot 7\cdot 71}\tx{or}
    \amn{3^3\cdot 5\cdot 17\cdot 31}{3^3\cdot 5\cdot 7\cdot 71.}\]
\end{Example}

\newpage

\begin{Problem}[4]
  {\em To find amicable numbers of the form \(agpq\) and
    \(ahr\), where \(p\), \(q\), \(r\) are prime numbers, but
    \(g\) and \(h\), whether prime or composite, are given along
    with one common factor \(a\).}
\end{Problem}

\begin{SolutionUnnum}
  From the common factor \(a\) we seek the fraction
  \(\dfrac{b}{c}=\dfrac{a}{2a-\int a}\) in lowest terms; then
  let \(\dfrac{\int g}{\int h}=\dfrac{m}{n}\) and from the first
  property of amicable numbers we get
  \[(p+1)(q+1)\int g=(r+1)\int h\tx{or}
    r+1=\frac{m}{n}(p+1)(q+1).\] But the other property yields
  \[(r+1)\int a\cdot \int h= a(gpq+hr);\]
  or, because \(\dfrac{\int a}{a}=\dfrac{2b-c}{b}\), we get
  \[(r+1)(2b-c)\int h=b(gpq+hr)\] and by substituting the value for \(r\)
  \[m(2b-c)(p+1)(q+1)\int h=b(ngpq+mh(p+1)(q+1)-nh).\] For the
  sake of brevity, let \(p+1=x\), \(q+1=y\); we
  get\footnote{Typo: \(\cdots+ag-nh)\) in the original. \JE}
  \[m(2b-c)xy\int h=b(mhxy+ngxy-ngx-ngy+ng-nh)\]
  or\runon
  \[\left(mbh+nbg-2mb\int h+mc\int
      h\right)xy-nbgx-nbgy=nb(h-g).\] For the sake of brevity, set
  \[e=b(mh+ng)-(2b-c)m\int h\] and we
  get \[eexy-nbgex-nbgey+nnbbgg=nnbbgg+nb(h-g)e\]
  or\runon \[(ex-nbg)(ey-nbg)=nnbbgg+nb(h-g)e.\] Therefore suppose
  \(nnbbgg+nb(h-g)e=MN\) and we get
  \begin{alignat*}{4}
    x&=\frac{M+nbg}{e}&\mbox{and}\quad& &y&=\frac{N+nbg}{e}&\mbox{or}\,\,\quad&\\
    p&=\frac{M+nbg}{e}-1,&& &q&=\frac{N+nbg}{e}-1,& r&=\frac{m}{n}xy-1.
  \end{alignat*}
  If these three numbers \(p\), \(q\), and \(r\) are prime, the
  amicable numbers will be \(agpq\) and \(ahr\), provided that
  the factors of either one are relatively prime.
\end{SolutionUnnum}

\begin{Corollary}
  If \(g\) and \(h\) are prime numbers, we get
  \(\dfrac{m}{n}=\dfrac{g+1}{h+1}\); therefore let \(g=km-1\) and
  \(h=kn-1\); we get \(\int h=kn\), whence
  \begin{align*}
    e&=b(2kmn-m-n)-(2b-c)kmn\\
    &=ckmn-b(m+n),\\
    MN&=nb\bigl(nb(km-1)^2+k(n-m)e\bigr)\\
    &=\bigl(ex-bn(km-1)\bigr)\bigl(ey-bn(km-1)\bigr)
  \end{align*}
  and\runon
  \[p=x-1,\quad q=y-1\tx{and also} r=\frac{m}{n}xy-1.\]
\end{Corollary}
\runon
\begin{Case}[1]
  Let \(m=1\), \(n=3\), therefore \(g=k-1\), \(h=3k-1\) and we get
  \[e=3ck-4b\et MN=3b\bigl(3b(k-1)^2+2ke\bigr)\] and so
  \[x=\frac{M+3b(k-1)}{e}\et y=\frac{N+3b(k-1)}{e}\] and lastly
  \(p=x-1\), \(q=y-1\), and \(r=\frac{1}{3}xy-1\).
\end{Case}

\begin{Example}[1]
  Let \(a=4\), \(b=4\), \(c=1\); we get
  \[e=3k-16\et MN=12\bigl(12(k-1)^2+2ke\bigr)\] and
  \[x=\frac{M+12(k-1)}{e}\et y=\frac{N+12(k-1)}{e}.\]
  Here one can put
  \begin{itemize}
  \item[I.] \(k=6\) making \(g=5\), \(h=17\), and \(e=2\), but
    nothing is produced from this.
  \item[II.] \(k=8\) making \(g=7\), \(h=23\), and \(e=8\),
    \(MN=12(12\cdot 49+128)\) or \(MN=16\cdot 3\cdot
    179=(8x-84)(8y-84)\) and so \(3\cdot 179=(2x-21)(2y-21)\),
    whence similarly nothing follows.
  \end{itemize}
\end{Example}

\begin{Example}[2]
  Let \(a=8\), \(b=8\), \(c=1\); we get
  \[e=3k-32,\quad MN=24\bigl(24(k-1)^2+2ke\bigr)\] or
  \[MN=48\bigl(15kk-56k+12\bigr) =
    \bigl(ex-24(k-1)\bigr)\bigl(ey-24(k-1)\bigr).\] But,
  again, one cannot conclude anything from this.
\end{Example}

\begin{Case}[2]
  Let \(m=3\), \(n=1\); we get
  \begin{gather*}
    e=3ck-4b\et g=3k-1,\quad h=k-1,\\
    MN=b\bigl(b(3k-1)^2-2ke\bigr) = \bigl(ex-b(3k-1)\bigr)\bigl(ey-b(3k-1)\bigr)
  \end{gather*}
  and also \(p=x-1\), \(q=y-1\), and \(r=3xy-1\).
\end{Case}

\begin{Example}[1]
  Let \(a=10\), \(b=5\), \(c=1\); we get
  \[e=3k-20\et 5\bigl(5(3k-1)^2-2ke\bigr) =
    \bigl(ex-5(3k-1)\bigr)\bigl(ey-5(3k-1)\bigr).\]
  If one puts here \(k=8\), we get \(5\cdot 29\cdot
  89=(4x-115)(4y-115)\).

  \par Whence appear \(x=30\), \(y=674\),
  \(3xy=60660\) and the amicable numbers will be:
  \[\amn{10\cdot 23\cdot 29\cdot 673}{10\cdot 7\cdot 60659.}\]
\end{Example}

\begin{Example}[2]
  Let \(a=3^3\cdot 5\), \(b=9\), \(c=2\); we get\footnote{Typo:
    \((93k-1)^2-3ke\) in original. \JE}
  \[e=6k-36\et 9(3k-1)^2-2ke =
    \left(\frac{1}{3}ex-3(3k-1)\right)
    \left(\frac{1}{3}ey-3(3k-1)\right).\]

  \par Now put \(k=8\); we get \(e=12\) and
  \(3\cdot 1523=(4x-69)(4y-69)\) and hence arise the values
  \(x=18\), \(y=398\), \(3xy=21492\), and the numbers \(g=23\),
  \(h=7\), \(p=17\), \(q=397\), \(r=21491\) will be prime and
  the amicable numbers will be:
  \[\amn{3^3\cdot 5\cdot 23\cdot 17\cdot 397}{3^3\cdot 5\cdot
      7\cdot 21491.}\]
\end{Example}

\begin{Scholium}
  From these examples, the use of this problem in finding
  amicable numbers is seen splendidly enough; but because of that
  same excess of freedom in arranging things, it is not a little
  tiresome to run through all cases according to the rules laid
  out here.  Therefore, since it suffices to have related this
  method and to have demonstrated its use, I will not linger
  longer on it, but proceed to expound the final method which I
  have used, by means of which amicable numbers can be
  extracted.  It relies, however, on remarkable properties,
  which numbers enjoy in ratio to a divisor sum, which I will
  explain when the occasion presents itself\footnote{Rudio here
    refers to \cite{euler243}. I could not see why, and suspect
    that Euler here means that he will explain the properties
    during the solution to Problem 5 (namely the deficiency of
    \(r\) and the inequalities between \(s\) and \(\int r\))
    rather than setting them up separately as lemmas first. \JE}
  lest we create tedium by giving further lemmas. But by
  explaining them it will not be difficult to resolve many more
  pertinent problems of this sort.
\end{Scholium}

\begin{Problem}[5]
  {\em To find amicable numbers of the form \(zap\) and \(zbq\),
    where the factors \(a\) and \(b\) are given, \(p\) and \(q\)
    are prime numbers and the common factor must be found.}
\end{Problem}

\begin{SolutionUnnum}
  Let \(\int a:\int b=m:n\), and since it must be that
  \(\int a\cdot (p+1)=\int b\cdot (q+1)\), we get
  \(m(p+1)=n(q+1)\). Suppose \(p+1=nx\) and \(q+1=mx\) and the
  amicable numbers will be
  \[za(nx-1)\et zb(mx-1),\] where indeed it is required that
  \(mx-1\) and \(nx-1\) are prime numbers. Now since the divisor
  sum of each number is the same
  \(=nx\int a\cdot \int z=mx\int b\cdot \int z\), it must be
  that this is equal to the sum of the numbers
  \(z\bigl((na+mb)x-a-b\bigr)\). Whence we obtain this equation:
  \[\frac{z}{\int z}=\frac{nx\int a}{(na+mb)x-a-b}.\]
  Now in order that the value of \(z\) may be found from this
  equation, the fraction \(\dfrac{nx\int a}{(na+mb)x-a-b}\) is
  reduced to lowest terms, say \(=\dfrac{r}{s}\), so that we
  have \(\dfrac{z}{\int z}=\dfrac{r}{s}\), and hence the
  following are to be noted. First \(z\) is to be equal to
  either \(r\) or a multiple of it, say \(kr\). In the first
  case, if \(z=r\), we get \(\int z=s\) and therefore \(s=\int r\). In the latter case, if \(z=kr\), we get
  \(\int z=ks=\int kr\). But whatever \(k\) is, we get
  \(\dfrac{\int kr}{\int r}>k\); for \(\int kr\) contains all
  the divisors of \(r\) each multiplied by \(k\), and on top of
  that, those divisors of\footnote{In the first edition (and
    also in \cite{eulerfuss}) this says \(r\) instead of
    \(kr\). Corrected \FR} \(kr\) which are not divisible by
  \(k\), and therefore we get \(\int kr > k\int r\). Therefore
  since \(\int z>k\int r\), we also get \(ks>k\int r\) or
  \(s>\int r\). In this way, if in the fraction \(\dfrac{r}{s}\)
  we have \(s=\int r\), we get \(z=r\); if however \(s>\int r\),
  we get that \(z\) equals some multiple of \(r\). Whence it is
  clear that if \(s<\int r\), the equation
  \(\dfrac{z}{\int z}=\dfrac{r}{s}\) is impossible, and one
  cannot find amicable numbers from this. Then,
  since\footnote{In the first edition (but not in
    \cite{eulerfuss}) here and in the next formula \(a-b\) is
    written in place of \(a+b\). Corrected. \FR}
  \[\frac{\int z}{z}=\frac{na+mb}{n\int a}-\frac{a+b}{nx\int
      a}=\frac{a}{\int a}+\frac{b}{\int b}-\frac{a+b}{nx\int
      a},\] because \(\dfrac{a}{\int a}<1\) and
  \(\dfrac{b}{\int b}<1\) we get
  \(\dfrac{\int z}{z}<2-\dfrac{a+b}{nx\int a}\) and all the more
  so \(\dfrac{z}{\int z}>\dfrac{1}{2}\), thus \(z\) is always a
  deficient number\footnote{i.e. its aliquot sum is less than
    itself. \JE}. And hence this shows that the equation
  \(\dfrac{z}{\int z}=\dfrac{r}{s}\) will always have the
  property that \(\dfrac{r}{s}>\dfrac{1}{2}\) or
  \(s<2r\). Whence if \(\int r=s\), we get \(\int r<2r\), and, if
  \(s>\int r\), we get (all the more so) \(\int r<2r\). In
  either case, therefore, \(r\) will be a deficient
  number. Wherefore if \(x\) were regarded as an unknown number,
  the value of \(x\) must be determined from the given equation
  \(\dfrac{z}{\int z}=\dfrac{nx\int a}{(na+mb)x-a-b}\) so that,
  by having reduced the fraction
  \(\dfrac{nx\int a}{(na+mb)x-a-b}\) to lowest terms
  \(\dfrac{r}{s}\), the number \(r\) becomes deficient and
  either \(s=\int r\) or \(s>\int r\).

  With these conditions observed, both \(r\) and \(s\) are
  resolved into their simplest prime factors, so that an
  equation like this holds:
  \[\frac{z}{\int z}=\frac{A^\alpha B^\beta
      C^\gamma}{E^\varepsilon F^\zeta G^\eta};\] but then,
  successively, either \(A^\alpha\) or a higher power of \(A\)
  is put as a factor of \(z\); that is, we put
  \(z=P\cdot A^{\alpha +\nu}\) and get
  \(\int z=\int A^{\alpha+\nu}\cdot \int P\) and
  \(\dfrac{z}{\int z}=\dfrac{PA^{\alpha+\nu}}{\int
    A^{\alpha+\nu}\cdot\int P}\) and so
  \[\frac{P}{\int P} = \frac{B^\beta C^\gamma\int
      A^{\alpha+\nu}}{A^\nu E^\varepsilon F^\zeta G^\eta}.\] And in
  the same way we further put \(P=B^{\beta+\mu} Q\), and proceed
  in this manner, until eventually we arrive at an equation of
  the form \(\dfrac{Z}{\int Z} = \dfrac{u}{\int u}\), from which
  we would have \(Z=u\). Indeed often this operation lacks the
  desired success, but for any given case this operation will be
  more easily taught through examples than through rules.
\end{SolutionUnnum}

\begin{Example}[1]
  Let \(a=3\), \(b=1\); we get  \(\int a=4\), \(\int b=1\) and
  \(m=4\), \(n=1\), and the amicable numbers will be
  \[3(x-1)z\et (4x-1)z,\] if \(x-1\) and\footnote{Typo: \(3x-1\)
    in the original. \JE} \(4x-1\) are prime numbers and
  \[\frac{z}{\int z}=\frac{4x}{7x-4}.\]
  But in the first place it is clear that if \(4\) were not
  cancelled from the numerator we would get \(7x-4<\int 4x\)
  because \(\int 4x=7\int x\). Therefore it is necessary that
  \(7x-4\) be an even number. Suppose \(x=4p\); we get
  \[\frac{z}{\int z}=\frac{4p}{7p-1}.\]
  Now make \(7p-1\) even by putting \(p=2q+1\); we get
  \[\frac{z}{\int z}=\frac{2(2q+1)}{7q+3}\] and \(x=8q+4\) and
  also
  \[x-1=8q+3,\quad 4x-1=32q+15.\]

  \par Whence \(q\) cannot be a multiple of three, lest \(x-1\)
  become divisible by \(3\). Therefore we get either \(q=3r+1\)
  or \(q=3r-1\); in the former case we get \(2q+1=6r+3\) and
  \(z\) would have to be\footnote{To see this, note that the
    numerator of \(z/\int z\) contains \(2q+1=3(2r+1)\) and
    there is no way to cancel this factor of \(3\) because the
    denominator is \(7q+3\) with \(q\neq 0\mod 3\), so \(3\)
    must divide the numerator in \(z/\int z\). \JE} divisible by
  \(3\), which equally cannot happen, because there is already a
  factor of \(3\) in the other number we seek
  \(3(x-1)z\). Therefore let \(q=3r-1\); we get
  \[\frac{z}{\int z}=\frac{2(6r-1)}{21r-4}\] and also
  \(x=24r-4\), \[x-1=24r-5\et 4x-1=96r-17.\] Unless two were
  cancelled from the numerator \(2(6r-1)\), \(z\) would be
  divisible by \(2\) and by putting \(z=2y\) we would get
  \[\frac{2y}{3\int y} =\frac{2(6r-1)}{21r-4}\et \frac{y}{\int
      y}=\frac{3(6r-1)}{21r-4}\] and so \(y\) (and therefore \(z\)) would have turned out to be divisible
  by \(3\), which cannot be, since \(z\) cannot have \(3\) as a
  factor\footnote{In the original, Euler starts with {\em Cum
      autem \(z\) factorem \(3\) habere nequit}, but only
    appeals to this fact at the end of the (ungainly)
    sentence. I hope I have not sacrificed any subtle meaning by
    transposing the sentence. In fact, the original has a typo
    here, and says {\em Cum autem \(z\) factorem \(4\) habere
      nequit}, which is tacitly corrected in the Opera Omnia
    edition. \JE}. Because of this, two must be cancelled from
  the numerator by putting \(r=2s\), so that
  \[x-1=48s-5,\quad 4x-1=192s-17,\] and we get
  \[\frac{z}{\int z}=\frac{12s-1}{21s-2}.\]
  Now if \(s\) were an odd number, because \(z\) is an odd
  number, also \(\int z=k(21s-2)\) would be an odd number, from
  which it follows that \(z\) would be square;\footnote{Note
    that if \(z\) is odd and not square then its divisors (all
    odd) come in pairs, which makes \(\int z\) even. \JE} but if
  however \(s\) were an even number, the common factor \(z\)
  would not be square.\footnote{Conversely, if \(z\) is odd and
    square, \(\int z\) equals a sum of even contributions (from
    pairs of odd divisors) plus the (odd) square root, so
    \(\int z\) is odd. \JE} Therefore we work out those values
  of \(s\) which make \(x-1=48s-5\) and \(4x-1=192s-17\) prime
  numbers and discern whether the equation
  \(\dfrac{z}{\int z}=\dfrac{12s-1}{21s-2}\) may be satisfied.

  Let \(s=7\); we get \(x-1=331\), \(4x-1=1327\) and
  \(\dfrac{z}{\int z}=\dfrac{83}{145}\). Now since \(z\) must be
  square, suppose \(z=83^2 A\); we get
  \(\int z=367\cdot 19\int A\) and
  \(\dfrac{A}{\int A}=\dfrac{367\cdot 19}{5\cdot 29\cdot
    83}\). But now \(19^2\) cannot have been a factor of \(A\)
  because \(\int 19^2=3\cdot 127\); indeed \(3\) would appear as
  a factor of \(A\),\footnote{To see this, set \(A=19^2B\). We
    get
    \(\tfrac{367\cdot 19}{5\cdot 29\cdot 83} = \tfrac{A}{\int
      A}= \tfrac{19^2B}{3\cdot 127\int B}\), so
    \(\tfrac{B}{\int B}=\tfrac{3\cdot 127\cdot 367}{5\cdot
      29\cdot 83\cdot 19}\) and \(3\) is now in the numerator,
    and hence a divisor of \(B\). \JE} and moreover taking
  higher powers one soon arrives at numbers so big that it is
  easy to see that the work cannot succeed.

  Let \(s=12\); we get \(x-1=571\), \(4x-1=2287\) and
  \(\dfrac{z}{\int z}=\dfrac{11\cdot 13}{2\cdot 125}\), which
  cannot be resolved by assuming that either\footnote{Typo in
    both original and Opera Omnia editions: \(11^2\). Indeed,
    \(s\) is even, so factors of \(z\) need not be square. \JE}
  \(11\) or \(13\) are factors of \(z\).

  Nor indeed has it been permitted for me to do any better from
  the larger values of \(s\).
\end{Example}

\begin{Example}[2\footnotemark]
  \footnotetext{Dickson {\cite[pp.44--45]{DicksonHistory1}}
    points out that Euler's approach in this case, restricting
    to the overlapping cases \(x=3(3q+1)\) and \(x=2(2q+1)\), is
    both incomplete and partly redundant. He gives a succinct
    summary of a more systematic approach. \JE} Let \(a=5\),
  \(b=1\); we get \(\int a=6\), \(\int b=1\), \(m=6\), \(n=1\),
  and the amicable numbers will be
  \[5(x-1)z\et (6x-1)z;\] we also have
  \[\frac{z}{\int z}=\frac{6x}{11x-6}.\] To make this equation
  possible, either two or three must be cancelled from the
  numerator \(6x\), because otherwise the numerator remains an
  abundant number.\footnote{Any number divisible by \(6\) is
    abundant. One can see this because if \(x=2^m3^ny\) with
    \(\gcd(6,y)=1\) then \(\int x=\int 2^m3^n\int y\) and one
    can prove (e.g. by double induction starting with \(m=n=1\))
    that \(\int 2^m3^n\geq 2(2^m3^n)\). \JE} Therefore we will
  have two cases to work out.

  \medskip
  
  \begin{itemize}[align=right,itemindent=2em,labelsep=2pt,labelwidth=1em,leftmargin=0pt,nosep]
  \item[I.] Let three be cancelled from the numerator by putting \(x=3p\);
    we get
    \[\frac{z}{\int z}=\frac{6p}{11p-2};\] now indeed moreover
    suppose that \(p=3q+1\), and we get
    \[\frac{z}{\int z}=\frac{2(3q+1)}{11q+3}\]
    and because \(x=9q+3\) the prime numbers must be
    \[x-1=9q+2\et 6x-1=54q+17,\] where clearly
    \(q\) must be an odd number. Therefore let \(q=2r-1\); we get
    \[x-1=18r-7,\quad 6x-1=108r-37\et \frac{z}{\int
        z}=\frac{2(6r-2)}{22r-8}=\frac{2(3r-1)}{11r-4}.\]
    We now work out the cases in which \(18r-7\) and \(108r-37\)
    become prime numbers, which are:
    \begin{itemize}
    \item[1)] \(r=1\); we get
      \[x-1=11,\quad 6x-1=71\et \frac{z}{\int z}=\frac{2\cdot
          2}{7}=\frac{4}{7}.\] Therefore since we have here
      \(7=\int 4\), we get \(z=4\) and the amicable numbers will
      be \(\amn{4\cdot 5 \cdot 11}{4\cdot 71}\) which indeed we
      have already found.\footnote{This is the pair \(220\),
        \(284\) (\S 65, I). \JE}
    \item[2)] \(r=2\); we get
      \[x-1=29,\quad 6x-1=179\et \frac{z}{\int z}=\frac{2\cdot
          5}{2\cdot 9}=\frac{5}{9}.\] But \(z\) cannot have a
      factor of \(5\).
    \item[3)] \(r=5\); we get\footnote{Typo:
        \(\frac{z}{\int a}\) in original. \JE}
      \[x-1=83,\quad 6x-1=503\et \frac{z}{\int z}=\frac{4\cdot
          7}{3\cdot 17};\] but here\footnote{Contradicting the
        fact that \(s\geq \int r\) established in \S 108. \JE}
      \(3\cdot 17<\int 4\cdot 7\).
    \item[4)] \(r=8\); we get
      \[x-1=137,\quad 6x-1=827\et \frac{z}{\int
          z}=\frac{23}{2\cdot 3\cdot 7}.\]
      Suppose\footnote{Typo: \(z=4\cdot 23P\) in original. \JE}
      \(z=23 P\); we get
      \[\int z=24\int P\et \frac{P}{\int
          P}=\frac{24}{23}\cdot\frac{z}{\int z}=\frac{4}{7};\]
      whence \(P=4\) and \(z=4\cdot 23\), which operation I will
      more succinctly represent thus
      \[\frac{z}{\int z}=\frac{23}{2\cdot 3\cdot
          7}\,\efrac{23}{24}\frac{4}{7}\efrac{4}{7};\] whence we
      get \(z=4\cdot 23\) and the amicable numbers will be
      \[\amn{4\cdot 23\cdot 5\cdot 137}{4\cdot 23\cdot 827.}\]
      The remaining values, certainly as
      far as I have examined, give no amicable numbers.
    \end{itemize}

  \item[II.] Let two be cancelled from the numerator by putting
    \(x=2p\); we get
    \[\frac{z}{\int z}=\frac{6p}{11p-3}.\]
    Now let\footnote{Otherwise \(2\) is not cancelled. \JE}
    \(p=2q+1\); we get
    \[\frac{z}{\int z}=\frac{3(2q+1)}{11q+4}\] and (because
    \(x=4q+2\)) the prime numbers must be
    \[x-1=4q+1,\quad 6x-1=24q+11;\] wherefore it is not possible
    that\footnote{Otherwise \(3\) divides \(x-1\). \JE}
    \(q=3\alpha-1\). Then since \(z\) cannot be divisible by
    \(5\), neither \(2q+1\) nor \(4q+1\) nor \(24q+11\) can be
    divisible by \(5\), whence the cases \(q=5\alpha+2\),
    \(q=5\alpha+1\) are excluded. By rejecting therefore this
    and other unusable values of \(q\) which do not give prime
    numbers for \(x-1\) and \(6x-1\), the calculation will be:
    
    \begin{center}
      \bgroup
      \def\arraystretch{2.5}
      \begin{longtable}{
          >{\raggedleft\arraybackslash}p{0.6cm}|
          >{\raggedleft\arraybackslash}p{0.8cm}|
          >{\raggedleft\arraybackslash}p{1cm}|p{9cm}}
        \(q\) & \(x-1\) & \(6x-1\) & \(\dfrac{z}{\int z}=\dfrac{3(2q+1)}{11q+4}\)\\
        \hline\endhead
        \(3\) & \(13\) & \(83\) & \(\dfrac{3\cdot 7}{37}\) which gives nothing.\\
        \(4\) & \(17\) & \(107\) & \(\dfrac{3\cdot 9}{48} = \dfrac{9}{16} \efrac{9}{13} \dfrac{13}{16} \efrac{13}{14} \dfrac{7}{8} \efrac{7}{8}\), \(z=9\cdot 7\cdot 13\);\\
        & & &or \(\dfrac{9}{16} \efrac{27}{40} \dfrac{5}{6} \efrac{5}{6}\), therefore \(z=27\cdot 5\). But this value\\
        &&& is no use because \(a=5\). The amicable numbers will\\
        &&&therefore be \(\amn{9\cdot 7\cdot 13\cdot 5\cdot 17}{9\cdot 7\cdot 13\cdot 107.}\)\\[10pt]
        \hline
        \(9\) & \(37\) & \(227\) & \(\dfrac{3\cdot 19}{103}\) which gives nothing.\\
        \(10\) & \(41\) & \(251\) & \(\dfrac{3\cdot 21}{114} = \dfrac{3\cdot 7}{2\cdot 19} \efrac{7^2}{3\cdot 19} \dfrac{3^2}{2\cdot 7} \efrac{3^2}{13} \dfrac{13}{14} \efrac{13}{14}\).\\
        & & & Therefore \(z=3^2\cdot 7^2\cdot 13\), and the amicable numbers\\
        &&&will be \(\amn{3^2\cdot 7^2\cdot 13\cdot 5\cdot 41}{3^2\cdot 7^2\cdot 13\cdot 251.}\)\\[10pt]
        \hline
        \(18\) & \(73\) & \(443\) & \(\dfrac{3\cdot 37}{202} = \dfrac{3\cdot 37}{2\cdot 101}\) which gives nothing.\\
        \(24\) & \(97\) & \(587\) & \(\dfrac{3\cdot 49}{268}=\dfrac{3\cdot 49}{4\cdot 67}\) which gives nothing.\\
        \(28\) & \(113\) & \(683\) & \(\dfrac{3\cdot 57}{312}=\dfrac{9\cdot 19}{8\cdot 39}=\dfrac{3\cdot 19}{8\cdot 13}\) which gives nothing.\\
        \(34\) & \(137\) & \(827\) & \(\dfrac{3\cdot 69}{378}=\dfrac{23}{2\cdot 21}=\dfrac{23}{2\cdot 3\cdot 7} \efrac{23}{24} \dfrac{4}{7} \efrac{4}{7}\), \(z=4\cdot 23\) as before.\\
        \(39\) & \(157\) & \(947\) & \(\dfrac{3\cdot 79}{433}\) which gives nothing.\\
        \(45\) & \(181\) & \(1091\) & \(\dfrac{3\cdot 91}{499} = \dfrac{3\cdot 7\cdot 13}{499}\)\\
        \(48\) & \(193\) & \(1163\) & \(\dfrac{3\cdot 97}{532} = \dfrac{3\cdot 97}{4\cdot 7\cdot 19} = \dfrac{3\cdot 97}{4\cdot 133} \efrac{97}{2\cdot 7^2} \dfrac{3\cdot 7}{2\cdot 19} \efrac{7^2}{3\cdot 19} \dfrac{3^2}{2\cdot 7}\)\\*
        &&& \(\efrac{3^2}{13} \dfrac{13}{14}\). Therefore
        \(z=3^2\cdot 7^2\cdot 13\cdot 97\), and the amicable
        numbers are \(\amn{3^2\cdot 7^2\cdot 13\cdot 97\cdot 5\cdot 193}{3^2\cdot 7^2\cdot 13\cdot 97\cdot 1163.}\)\\[10pt]
        \hline
        \(49\) & \(197\) & \(1187\) & \(\dfrac{3\cdot 99}{543}=\dfrac{9\cdot 11}{181}\)\\
        \(60\) & \(241\) & \(1451\) & \(\dfrac{3\cdot 121}{664}=\dfrac{3\cdot 11^2}{8\cdot 83}\)\\
        \(69\) & \(277\) & \(1667\) & \(\dfrac{3\cdot 139}{763}\) )\footnotemark\\
        \(79\) & \(317\) & \(1907\) & \(\dfrac{3\cdot 159}{873}=\dfrac{53}{97}\)\\
        \(84\) & \(337\) & \(2027\) & \(\dfrac{3\cdot 169}{928}=\dfrac{3\cdot 169}{8\cdot 116}=\dfrac{3\cdot 169}{32\cdot 29}\)\\
        \(93\) & \(373\) & \(2243\) & \(\dfrac{3\cdot 187}{1027} = \dfrac{3\cdot 11\cdot 17}{13\cdot 79}\)\\
        \(100\) & \(401\) & \(2411\) & \(\dfrac{3\cdot 201}{1104}=\dfrac{3\cdot 67}{368} = \dfrac{3\cdot 67}{16\cdot 23}\)\\
        \(244\) & \(977\) & \(5867\) & \(\dfrac{3\cdot 489}{2688}=\dfrac{3\cdot 163}{128\cdot 7}\efrac{163}{4\cdot 41} \dfrac{3\cdot 41}{32\cdot 7}\efrac{41}{2\cdot 3\cdot 7} \dfrac{3^2}{16}\efrac{3^2}{13} \dfrac{13}{16}\)\\
        &&&\(\efrac{13}{14} \dfrac{7}{8}\). Therefore \(z=3^2\cdot 7\cdot 13\cdot 41\cdot 163\) and the amicable numbers will be \(\amn{3^2\cdot 7\cdot 13\cdot 41\cdot 163\cdot 5\cdot 977}{3^2\cdot 7\cdot 13\cdot 41\cdot 163\cdot 5867.}\)
      \end{longtable}
      \footnotetext{The denominator was given as \(793\) in the original. \JE}
      \egroup
    \end{center}
    \vspace{-1cm} Hence therefore two new amicable numbers have
    appeared.\footnote{Only two because some of these appeared
      earlier (\S 65, VII and VIII). \JE}
  \end{itemize}
\end{Example}

\newpage

\begin{Example}[3]
  Let \(a=7\), \(b=1\); we get \(\int a=8\), \(\int b=1\),
  \(m=8\), \(n=1\), and the amicable numbers
  \[7(x-1)z\et (8x-1)z\]
  giving rise to
  \[\dfrac{z}{\int z}=\dfrac{8x}{15x-8}.\]
  And first indeed \(x\) must be an even number; therefore
  suppose \(x=2p\); we get
  \[x-1=2p-1,\quad 8x-1=16p-1\]
  and\runon
  \[\frac{z}{\int z}=\frac{8p}{15p-4},\]
  which equation is impossible unless the power of two in the
  numerator is brought down, because\footnote{Contradicting the
    fact that \(s\geq \int r\) established in \S 108. \JE}
  \(15p-4<\int 8p\). Therefore put \(p=4q\), so that
  \[x=8q,\quad x-1=8q-1,\quad 8x-1=64q-1\] and\runon
  \[\frac{z}{\int z}=\frac{8q}{15q-1}.\]
  Now let \(q=2r+1\); we get
  \[\frac{z}{\int z}=\frac{4(2r+1)}{15r+7}\]
  and\runon
  \[x-1=16r+7,\quad 8x-1=128r+63;\] in order that neither of
  these numbers be divisible by \(3\), we have neither
  \(r=3\alpha-1\) nor \(r=3\alpha\). Therefore let \(r=3s+1\); we get
  \[\frac{z}{\int z}=\frac{4(6s+3)}{45s+22}\tx{or}
    \frac{z}{\int z}=\frac{4\cdot 3(2s+1)}{45s+22}\]
  and\runon
  \[x-1=48s+23,\quad 8x-1=384s+191.\] Now either three or four
  must be cancelled from the numerator.\footnote{Otherwise \(z\)
    is abundant, see \S 110. \JE} But three cannot be cancelled,
  because the denominator is never divisible
  by\footnote{\(45s+22=1\mod 3\) \JE} \(3\); therefore let four be cancelled
  four, to what end I put \(s=2t\), and we get
  \[\frac{z}{\int z}=\frac{2\cdot 3(4t+1)}{45t+11};\]
  now let \(t=2u-1\); we get
  \[\frac{z}{\int z}=\frac{3(8u-3)}{45u-17};\]
  but \(s=4u-2\) and so the following numbers must be prime
  \[x-1=192u-73,\quad 8x-1=1536u-577.\]

  \begin{center}
    \bgroup
    \def\arraystretch{2.5}
    \begin{longtable}{
        >{\raggedleft\arraybackslash}p{0.6cm}|
        >{\raggedleft\arraybackslash}p{0.8cm}|
        >{\raggedleft\arraybackslash}p{1cm}|p{9cm}}
      \(u\) & \(x-1\) & \(8x-1\) & \(\dfrac{z}{\int z}\)\\
      \hline\endhead
      \(5\) & \(887\) & \(7103\) & \(\dfrac{3\cdot 37}{208} = \dfrac{3\cdot 37}{16\cdot 13} \efrac{37}{2\cdot 19} \dfrac{3\cdot 19}{8\cdot 13} \efrac{19}{4\cdot 5} \dfrac{3\cdot 5}{2\cdot 13} \efrac{5}{2\cdot 3} \dfrac{3^2}{13}\).\\
      & & & Therefore \(z=3^2\cdot 5\cdot 19\cdot 37\), and the amicable numbers\\
      &&&will be \(\amn{3^2\cdot 5\cdot 19\cdot 37\cdot 7\cdot 887}{3^2\cdot 5\cdot 19\cdot 37\cdot 7103.}\)\\[10pt]
      \hline
      \(11\) & \(2039\) & \(16319\) & \(\dfrac{3\cdot 5\cdot 17}{2\cdot 239}\))\footnotemark\\
      \(13\) & \(2423\) & \(19391\) & \(\dfrac{3\cdot 101}{8\cdot 71}\)\\
      \(26\) & \(4919\) & \(39359\) & \(\dfrac{3\cdot 205}{1153}\).
    \end{longtable}
    \footnotetext{Typo: \(\tfrac{3\cdot 5 \cdot 17}{4\cdot 107}\) in the original. \JE}
    \egroup
  \end{center}
\end{Example}

\begin{Example}[4]
  Let \(a=11\), \(b=1\); we get \(\int a=m=12\), \(\int b=n=1\);
  the numbers we seek will be
  \[11(x-1)z\et (12x-1)z\]
  and also
  \[\dfrac{z}{\int z}=\dfrac{12x}{23x-12}.\]
  Here either \(3\) or \(4\) must be cancelled from the
  numerator.

  \medskip
  
  \begin{itemize}[align=right,itemindent=2em,labelsep=2pt,labelwidth=1em,leftmargin=0pt,nosep]
  \item[I.] Let \(3\) be cancelled; suppose \(x=3p\), we get
    \[\frac{z}{\int z}=\frac{12p}{23p-4},\]
    and \(p=3q-1\); we get
    \[\frac{z}{\int z}=\frac{4(3q-1)}{23q-9}\]
    and because \(x=9q-3\), \(q\) must be odd. Let \(q=2r+1\), so that \(x=18r+6\); we get
    \[\frac{z}{\int z}=\frac{4(6r+2)}{46r+14}=\frac{4(3r+1)}{23r+7}\]
    and
    \[x-1=18r+5,\quad 12x-1=216r+71.\]

    \begin{center}
      \bgroup
      \def\arraystretch{2}
      \begin{longtable}{
          >{\raggedleft\arraybackslash}p{0.6cm}|
          >{\raggedleft\arraybackslash}p{0.8cm}|
          >{\raggedleft\arraybackslash}p{1.2cm}|p{9cm}}
        \(r\) & \(x-1\) & \(12x-1\) & \(\dfrac{z}{\int z}\)\\
        \hline\endhead
        \(0\) & \(5\) & \(71\) & \(\dfrac{4}{7}\), \(z=4\);
        amicable numbers \(\amn{4\cdot 11\cdot 5}{4\cdot
          71.}\)\\[10pt]
        \hline
        \(2\) & \(41\) & \(503\) & \(\dfrac{4\cdot 7}{53}\)\\
        \(3\) & \(59\) & \(719\) & \(\dfrac{4\cdot 10}{76}=\dfrac{2\cdot 5}{19}\) impossible.\\
        \(6\) & \(113\) & \(1367\) & \(\dfrac{4\cdot 19}{145}=\dfrac{4\cdot 19}{5\cdot 29}\) impossible.\footnotemark\\
        \(7\) & \(131\) & \(1583\) & \(\dfrac{4\cdot 22}{168}=\dfrac{11}{21}=\dfrac{11}{3\cdot 7} \efrac{11}{12}\dfrac{4}{7}\),\\
        &&&but because of the factor of \(11\) this value of \(z\) is not valid
      \end{longtable}
      \footnotetext{Typo: \(\tfrac{4\cdot 10}{5\cdot 29}\) in original. \JE}
      \egroup
    \end{center}
    \vspace{-\baselineskip}
  \item[II.] Let the factor of \(4\) be cancelled and suppose \(x=4p\), making
    \[\frac{z}{\int z}=\frac{12p}{23p-3}.\]
    Now let \(p=4q+1\); we get
    \[\frac{z}{\int z}=\frac{3(4q+1)}{23q+5}\]
    and because \(x=16q+4\) the following numbers must be prime
    \[x-1=16q+3\et 12x-1=192q+47;\]
    hence we exclude the values \(q=3\alpha\).

    \begin{center}
      \bgroup
      \def\arraystretch{2}
      \begin{longtable}{
          >{\raggedleft\arraybackslash}p{0.6cm}|
          >{\raggedleft\arraybackslash}p{0.8cm}|
          >{\raggedleft\arraybackslash}p{1.2cm}|p{9.8cm}}
        \(q\) & \(x-1\) & \(12x-1\) & \(\dfrac{z}{\int z}\)\\
        \hline\endhead
        \(0\) & \(3\) & \(47\) & \(\dfrac{3}{5}\) impossible.\\
        \(1\) & \(19\) & \(239\) & \(\dfrac{3\cdot 5}{4\cdot 7} \efrac{5}{2\cdot 3}\dfrac{3^2}{14} \efrac{3^2}{13} \dfrac{13}{14}\); \(z=3^2\cdot 5\cdot 13\) and the amicable numbers will be \(\amn{3^2\cdot 5\cdot 13\cdot 11\cdot 19}{3^2\cdot 5\cdot 13\cdot 239.}\)\\[10pt]
        \hline
        \(13\) & \(211\) & \(2543\) & \(\dfrac{3\cdot 53}{16\cdot 19} \efrac{53}{2\cdot 27} \dfrac{81}{8\cdot 19} \efrac{243}{4\cdot 7\cdot 13} \dfrac{7\cdot 13}{2\cdot 3\cdot 19} \efrac{13}{2\cdot 7} \dfrac{7^2}{3\cdot 19} \efrac{7^2}{3\cdot 19}.\) \\*
        & & & Therefore \(z=3^5\cdot 7^2\cdot 13\cdot 53\) and the amicable numbers\\
        &&&will be
        \(\amn{3^5\cdot 7^2\cdot 13\cdot 53\cdot 11\cdot 211}{3^5\cdot 7^2\cdot 13\cdot 53\cdot 2543.}\)
      \end{longtable}
      \egroup
    \end{center}
  \end{itemize}
\end{Example}

\vspace{-3\baselineskip}

\begin{Example}[5]
  Let \(a=5\), \(b=17\), and the amicable numbers
  \[5(3x-1)z\et 17(x-1)z;\]
  we get
  \[\frac{z}{\int z}=\frac{18x}{32x-22}=\frac{9x}{16x-11}.\]
  Since \(x\) must be an even number, suppose \(x=2p\); we get
  \[\frac{z}{\int z}=\frac{18p}{32p-11},\]
  and from the numerator \(18p\) either the factor \(2\) or
  \(3^2\) must be cancelled, lest the numerator be an abundant
  number. But the factor \(2\) cannot be
  cancelled;\footnote{Because the denominator is necessarily
    odd. \JE} therefore let the factor \(9\) be cancelled. To that end,
  suppose \(p=9q+4\), so that \(x=18q+8\) and
  \[x-1=18q+7\et 3x-1=54q+23;\]
  we get
  \[\frac{z}{\int z}=\frac{2(9q+4)}{32q+13}.\]

  \begin{center}
    \bgroup
    \addtocounter{footnote}{1}
    \footnotetext{Typo: \(12x-1\) in original. \JE}
    \addtocounter{footnote}{-1}
    \def\arraystretch{2}
    \begin{longtable}{
        >{\raggedleft\arraybackslash}p{0.6cm}|
        >{\raggedleft\arraybackslash}p{0.8cm}|
        >{\raggedleft\arraybackslash}p{1.3cm}|p{9cm}}
      \(q\) & \(x-1\) & \(3x-1\)\footnotemark & \(\dfrac{z}{\int z}\)\\
      \hline
      \(0\) & \(7\) & \(23\) & \(\dfrac{8}{13}\) impossible.\\
      \(2\) & \(43\) & \(131\) & \(\dfrac{4\cdot 11}{7\cdot 11}=\dfrac{4}{7}\); \(z=4\)\\
      &&&and the amicable numbers \(\amn{4\cdot 5\cdot 131}{4\cdot 17\cdot 43.}\)\\
      \hline
      \(4\) & \(79\) & \(239\) & \(\dfrac{16\cdot 5}{3\cdot 47}\)\\
      \(5\) & \(97\) & \(293\) & \(\dfrac{2\cdot 49}{173}\) \\
      \(17\) & \(313\) & \(941\) & \(\dfrac{2\cdot 157}{557}\) \\
      \(19\) & \(349\) & \(1049\) & \(\dfrac{2\cdot 5^2\cdot 7}{27\cdot 23}\) \\
      \(20\) & \(367\) & \(1103\) & \(\dfrac{16\cdot 23}{653}\) \\
      \(24\) & \(439\) & \(1319\) & \(\dfrac{8\cdot 5\cdot 11}{781}\) no use, \(=\dfrac{8\cdot 5}{71}\). 
    \end{longtable}
    \egroup
  \end{center}
\end{Example}

\vspace{-\baselineskip}

\begin{Example}[6]
  Let \(a=37\) and \(b=227\); we get \(\int a=38\),
  \(\int b=228\), and \(\dfrac{m}{n}=\dfrac{1}{6}\); whence if the
  amicable numbers are
  \[37(6x-1)z\et 227(x-1)z,\]
  which makes
  \[\dfrac{z}{\int z}=\dfrac{6\cdot 38x}{449x-264} = \dfrac{4\cdot 3\cdot 19x}{449x-264},\]
  where, since \(x\) must be an even number, we put \(x=2p\), so
  that the prime numbers must be
  \[x-1=2p-1\et 6x-1=12p-1,\] and we get
  \[\frac{z}{\int z}=\frac{4\cdot 3\cdot 19p}{449p-132}.\]
  Now from the numerator either the factor of \(4\) or the factor of \(3\) must be cancelled.

  \medskip
  
  \begin{itemize}[align=right,itemindent=2em,labelsep=2pt,labelwidth=1em,leftmargin=0pt,nosep]
  \item[I.] Let the factor of \(3\) be cancelled; to that end suppose \(p=3q\), so that
    \[\frac{z}{\int z}=\frac{4\cdot 3\cdot 19q}{449q-44};\]
    now put \(q=3r+1\), and we get
    \[\frac{z}{\int z}=\frac{4\cdot 19(3r+1)}{449q+135}\]
    and \(p=9r+3\) and\footnote{Typo: \(6-1\) instead of
      \(6x-1\) in original. \JE}
    \(\begin{cases}x-1=18r+5,\\ 6x-1=108r+35.\end{cases}\)

    \begin{center}
      \bgroup
      \def\arraystretch{2}
      \begin{longtable}{ >{\raggedleft\arraybackslash}p{0.6cm}|
          >{\raggedleft\arraybackslash}p{0.8cm}|
          >{\raggedleft\arraybackslash}p{1cm}|p{9cm}}
        \(r\) & \(x-1\) & \(6x-1\) & \(\dfrac{z}{\int z}\)\\
        \hline\endhead
        \(2\) & \(41\) & \(251\) & \(\dfrac{4\cdot 19\cdot 7}{1033}\)\\
        \(3\) & \(59\) & \(359\) & \(\dfrac{4\cdot 19\cdot 10}{1482}=\dfrac{4\cdot 5}{3\cdot 13}\)\\
        \(6\) & \(113\) & \(683\) & \(\dfrac{4\cdot 19\cdot 19}{3\cdot 23\cdot 41}\) \\
        \(13\) & \(239\) & \(1439\) & \(\dfrac{4\cdot 19\cdot 40}{4\cdot 1493}\) \\
        \(17\) & \(311\) & \(1871\) & \(\dfrac{16\cdot 13\cdot 19}{8\cdot 971}\) \\
        \(22\) & \(401\) & \(2411\) & \(\dfrac{4\cdot 19\cdot 67}{10013}=\dfrac{4\cdot 67}{17\cdot 31} \efrac{67}{4\cdot 17}\dfrac{16}{31} \efrac{16}{31}\); \(z=16\cdot 67\).\\[10pt]
        &&&Amicable numbers: \(\amn{16\cdot 67\cdot 37\cdot 2411}{16\cdot 67\cdot 227\cdot 401.}\)\\[10pt]
        \hline
        \(117\) & \(2111\) & \(12671\) & \(\dfrac{4\cdot 19\cdot 352}{52668}=\dfrac{128\cdot 11\cdot 19}{4\cdot 7\cdot 9\cdot 11\cdot 19}=\dfrac{32}{63}\); \(z=32\),\\[10pt]
        &&&and amicable numbers: \(\amn{32\cdot 37\cdot 12671}{32\cdot 227\cdot 2111.}\)
      \end{longtable}
      \egroup
    \end{center}
    \vspace{-\baselineskip}
  \item[II.] Let the factor of \(4\) be cancelled; suppose \(p=4q\); we get
    \[\frac{z}{\int z}=\frac{4\cdot 3\cdot 19q}{449q-33};\]
    now let \(q=4r+1\); we get \(p=16r+4\) and
    \[x-1=32r+7,\quad 6x-1=192r+47\]
    and also
    \[\frac{z}{\int z}=\frac{3\cdot 19(4r+1)}{449r+104}.\]

    \begin{center}
      \bgroup
      \def\arraystretch{2}
      \begin{longtable}{
          >{\raggedleft\arraybackslash}p{0.6cm}|
          >{\raggedleft\arraybackslash}p{0.8cm}|
          >{\raggedleft\arraybackslash}p{1cm}|p{9cm}}
        \(r\) & \(x-1\) & \(6x-1\) & \(\dfrac{z}{\int z}\)\\
        \hline\endhead
        \(0\) & \(7\) & \(47\) & \(\dfrac{3\cdot 19}{8\cdot 13} \efrac{19}{4\cdot 5} \dfrac{3\cdot 5}{2\cdot 13} \efrac{5}{2\cdot 3} \dfrac{3^2}{13}\); \(z=3^2\cdot 5\cdot 19\)\\[10pt]
        & & &  and amicable numbers \(\amn{3^2\cdot 5\cdot 19\cdot 37\cdot 47}{3^2\cdot 5\cdot 19\cdot 227\cdot 7.}\)\\
        \hline
        \(2\) & \(71\) & \(431\) & \(\dfrac{9\cdot 19}{2\cdot 167}\)\\
        \(8\) & \(263\) & \(1583\) & \(\dfrac{3\cdot 19\cdot 33}{16\cdot 3\cdot 7\cdot 11}=\dfrac{3\cdot 19}{16\cdot 7} \efrac{19}{4\cdot 5} \dfrac{3\cdot 5}{4\cdot 7} \efrac{5}{2\cdot 3} \dfrac{3^2}{2\cdot 7} \efrac{3^2}{13} \dfrac{13}{14}\); \\
        & & & \(z=3^2\cdot 5\cdot 13\cdot 19\),\\
        &&&and amicable numbers
        \(\amn{3^2\cdot 5\cdot 13\cdot 19\cdot 37\cdot 1583}{3^2\cdot 5\cdot 13\cdot 19\cdot 227\cdot 263.}\)\\
        \hline
        \(15\) & \(487\) & \(2927\) & \(\dfrac{3\cdot 19\cdot 61}{7\cdot 977}\) \\
        \(23\) & \(743\) & \(4463\) & \(\dfrac{9\cdot 19\cdot 31}{9\cdot 19\cdot 61}=\dfrac{31}{61}\) \\
        \(26\) & \(839\) & \(5039\) & \(\dfrac{3\cdot 19\cdot 105}{2\cdot 3\cdot 13\cdot 151}=\dfrac{3\cdot 5\cdot 7 \cdot 19}{2\cdot 13\cdot 151}\)\\
        \(30\) & \(967\) & \(5807\) & \(\dfrac{3\cdot 19\cdot 11}{2\cdot 617}\)\\
        \(41\) & \(1319\) & \(7919\) & \(\dfrac{3\cdot 19\cdot 165}{9\cdot 121\cdot 17}=\dfrac{5\cdot 19}{11\cdot 17}\). \\
      \end{longtable}
      \egroup
    \end{center}
  \end{itemize}
\end{Example}

\vspace{-2\baselineskip}

\begin{Example}[7]
  Let \(a=79\), \(b=11\cdot 19=209\), \(\int a=80\), \(\int b=240\); we get \(m=1\), \(n=3\), and the amicable numbers are
  \[79(3x-1)z\et 11\cdot 19(x-1)z;\]
  we get
  \[\frac{z}{\int z}=\frac{240x}{446x-228}=\frac{120x}{223x-144}.\]
  Let \(x=2p\); we get
  \[\frac{z}{\int z}=\frac{120p}{223p-72}\]
  and the numbers \(2p-1\) and \(6p-1\) must be prime. But now
  from the numerator \(120p\) either the factor of \(8\) or
  \(3\) must be cancelled.

  \medskip
  
  \begin{itemize}[align=right,itemindent=2em,labelsep=2pt,labelwidth=1em,leftmargin=0pt,nosep]
  \item[I]. Let the factor of \(3\) be cancelled; let \(p=9q\); we get
    \[\frac{z}{\int z}=\frac{120q}{223q-8}\]
    and put \(q=3r-1\), so that
    \[\frac{z}{\int z}=\frac{40(3r-1)}{223r-77},\]
    \[p=27r-9,\quad x-1=54r-19\tx{and} 3x-1=162r-55.\] But now,
    because \(40\) is an abundant number, either \(5\) or \(4\)
    must be cancelled.

    \begin{itemize}
    \item[\(\alpha\))] Let \(5\) be cancelled and let \(r=5s-1\); we get
      \[\frac{z}{\int z}=\frac{8(15s-4)}{223s-60}\] and the
      numbers\footnote{Typo: \(x-1=470s-73\) in original. \JE}
      \(x-1=270s-73\), \(3x-1=810s-217\) must be prime. And lest
      three enters anew into the numerator, the cases
      \(s=3\alpha-1\) are to be excluded.\footnote{It is to be
        observed indeed that the numerator is never divisible by
        \(3\). \FR} And hence nothing is found.

    \item[\(\beta\))] Since
      \(\dfrac{z}{\int z}=\dfrac{40(3r-1)}{223r-77}\), let \(4\)
      be cancelled and let \(r=4s-1\); we get
      \[\dfrac{z}{\int z}=\dfrac{10(12s-4)}{223s-75}=\dfrac{40(3s-1)}{223s-75};\]
      moreover let \(s=4t+1\); we get
      \[\dfrac{z}{\int z}=\dfrac{10(12t+2)}{223t+37}=\dfrac{20(6t+1)}{223s+37}.\]
      Moreover let \(t=2u-1\); we get
      \[\dfrac{z}{\int z}=\dfrac{10(12u-5)}{223u-93}\]
      and because \(r=16t+3=32u-13\) we get
      \(\begin{cases}x-1=1728u-721,\\
        3x-1=5184u-2161.\end{cases}\) But a smaller value than
      \(16\) for \(u\) does not render these numbers prime,
      whence we get
      \(\dfrac{z}{\int z}=\dfrac{2\cdot 11\cdot 17}{5\cdot
        139}\), which is no use because of the factor of \(11\).
    \end{itemize}
  \item[II.] Therefore from the equation
    \(\dfrac{z}{\int z}=\dfrac{120p}{223p-72}\) let the
    factor of \(8\) be cancelled. Suppose \(p=8q\); we get
    \[\frac{z}{\int z}=\frac{120q}{223q-9}\]
    and now let \(q=8r-1\); we get
    \[\frac{z}{\int z}=\frac{3\cdot 5(8r-1)}{223r-29};\]
    but because \(p=64r-8\) we get
    \[x-1=128r-17,\quad 3x-1=384r-49.\] Whence the values
    \(r=3\alpha+1\) and \(r=5\alpha\pm 1\) are excluded.

    \begin{center}
      \bgroup
      \def\arraystretch{2}
      \begin{longtable}{
          >{\raggedleft\arraybackslash}p{0.6cm}|
          >{\raggedleft\arraybackslash}p{0.8cm}|
          >{\raggedleft\arraybackslash}p{1cm}|p{9cm}}
        \(r\) & \(x-1\) & \(3x-1\) & \(\dfrac{z}{\int z}\)\\
        \hline\endhead
        \(2\) & \(239\) & \(719\) & \(\dfrac{3\cdot 5^2}{139}\)\\
        \(3\) & \(367\) & \(1103\) & \(\dfrac{3\cdot 23}{128} \efrac{23}{8\cdot 3} \dfrac{3^2}{16} \efrac{3^2}{13} \dfrac{13}{16} \efrac{13}{14} \dfrac{7}{8}\),\\
        &&&therefore \(z=3^2\cdot 7\cdot 13\cdot 23\), or\\
        & & & \(\dfrac{3\cdot 23}{128} \efrac{23}{8\cdot 3} \dfrac{3^2}{16} \efrac{3^3}{8\cdot 5} \dfrac{5}{6}\), therefore \(z=3^3\cdot 5\cdot 23\),
      \end{longtable}
      \egroup
      and the amicable numbers will be
      \[\amn{3^2\cdot 7\cdot 13\cdot 23 \cdot 79\cdot 1103}{3^2\cdot 7\cdot 13\cdot 23 \cdot 11\cdot 19\cdot 367}\tx{or}
      \amn{3^3\cdot 5\cdot 23 \cdot 79\cdot 1103}{3^3\cdot 5\cdot 23 \cdot 11\cdot 19\cdot 367.}
    \]
    \end{center}
  \end{itemize}
\end{Example}

\begin{Example}[8]
  Let \(a=17\cdot 19\), \(b=11\cdot 59\); we get
  \(\int a=18\cdot 20\), \(b=12\cdot 60\), and \(m=1\),
  \(n=2\). Therefore if we suppose the amicable numbers are
  \[\begin{array}{c}17\cdot 19(2x-1)z\\
      11\cdot 59(x-1)z,\end{array}\quad \mbox{we get}\quad
    \frac{z}{\int z}=\frac{720x}{1295x-972}.\]
  Let \(x=2p\); we get
  \[\frac{z}{\int z}=\frac{720p}{1295p-486}\et
    \begin{cases}x-1=2p-1\\ 2x-1=4p-1,\end{cases}\] neither of
  which is divisible by \(3\), it must be that \(p=3q\), so that
  \[\frac{z}{\int
      z}=\frac{720q}{1295q-162}\et\begin{cases}x-1=6q-1,\\
      2x-1=12q-1.\end{cases}\]
  Let the factor of \(16\) be cancelled from the numerator and let \(q=2r\); we get
  \[\frac{z}{\int z}=\frac{720r}{1295r-81};\]
  now let \(r=16s-1\); we get
  \[\frac{z}{\int
      z}=\frac{45(16s-1)}{1295s-86}\et \begin{cases}x-1=192s-13\\
      2x-1=384s-25.\end{cases}\]
  Let \(s=1\); we get \(x-1=179\), \(2x-1=359\) and
  \[\frac{z}{\int z}=\frac{45\cdot
      15}{1209}=\frac{225}{403}=\frac{3^2\cdot 5^2}{13\cdot
      31}\efrac{3^2}{13}\frac{5^2}{31}\efrac{5^2}{31}.\]
  Therefore \(z=3^2\cdot 5^2\) and the amicable numbers will be
  \[\amn{3^2\cdot 5^2\cdot 17\cdot 19\cdot 359}{3^2\cdot 5^2\cdot 11\cdot 59\cdot 179.}\]
\end{Example}

\begin{Scholium}
  This final method expounded in Problem 5 is completely
  different from the preceding method, which the previous four
  problems encompassed: while in this one the common factor is
  sought, in that one it is given. Each nevertheless is
  possessed of a singular kind of excellence, so that one
  without the help of the other is not apt enough to increase the
  multitude of amicable numbers. Indeed the latter method
  supplies the kind of common factors, which one could hardly
  have suspected for the use of the former; and verily the former
  suggests further factors suitable for this purpose. Moreover,
  everything which I have related here comprises a specimen of a
  highly undependable method, which, as far as possible, I have
  reduced to algebraic rules, so as to limit the vague
  uncertainty of what was to be tried. In place of the colophon,
  therefore, I will append the more than sixty pairs of amicable
  numbers which I have elicited by this method.
\end{Scholium}

\begin{Catalogue}
  \bgroup
  \def\arraystretch{2.5}
  \setlength\tabcolsep{2pt}
  \begin{center}
    \vspace{-\baselineskip}
    \begin{tabular}{p{1.2cm}p{3.5cm}p{1.2cm}p{3.5cm}p{1.2cm}p{3.5cm}}
      \nml{I.} & \(\amn{2^2\cdot 5\cdot 11}{2^2\cdot 71}\)
      & \nml{II.} & \(\amn{2^4\cdot 23\cdot 47}{2^4\cdot 1151}\)
      & \nml{III.} & \(\amn{2^7\cdot 191\cdot 383}{2^7\cdot 73727}\)
    \end{tabular}
    \vspace{-\baselineskip}
    \begin{longtable}{p{1.2cm}p{4.8cm}p{1.1cm}p{7cm}}
      \nml{IV.} & \(\amn{2^2\cdot 23\cdot 5\cdot 137}{2^2\cdot 23\cdot 827}\)&
      \nml{V.} &
      \(\amn{3^2\cdot 7\cdot 13\cdot 5 \cdot 17}{3^2\cdot 7\cdot
        13\cdot 107}\)\\
      \nml{VI.} & \(\amn{3^2\cdot 5\cdot 13\cdot 11\cdot 19}{3^2\cdot 5\cdot 13\cdot 239}\)&
      \nml{VII.} &
      \(\amn{3^2\cdot 7^2\cdot 13\cdot 5\cdot 41}{3^2\cdot
        7^2\cdot 13\cdot 251}\)\\
      \nml{VIII.} & \(\amn{3^2\cdot 5\cdot 7 \cdot 53\cdot 1889}{3^2\cdot 5\cdot 7\cdot 102059}\)&
      \nml{IX.} &
      \(\amn{2^2\cdot 13\cdot 17\cdot 389\cdot 509}{2^2\cdot
        13\cdot 17\cdot 198899}\)\\
      \nml{X.} & \(\amn{3^2\cdot 5\cdot 19\cdot 37\cdot 7\cdot 887}{3^2\cdot 5\cdot 19\cdot 37\cdot 7103}\)&
      \nml{XI.} &
      \(\amn{3^4\cdot 5\cdot 11\cdot 29\cdot 89}{3^4\cdot 5\cdot
        11\cdot 2699}\)\\
      \nml{XII.} & \(\amn{3^2\cdot 7^2\cdot 11\cdot 13\cdot 41\cdot 461}{3^2\cdot 7^2\cdot 11\cdot 13\cdot 19403}\)&
      \nml{XIII.} &
      \(\amn{3^2\cdot 5\cdot 13\cdot 19\cdot 29\cdot 569}{3^2\cdot
        5\cdot 13\cdot 19\cdot 17099}\)\\
      \nml{XIV.} & \(\amn{3^2\cdot 7^2\cdot 13\cdot 97\cdot 5\cdot 193}{3^2\cdot 7^2\cdot 13\cdot 97\cdot 1163}\)&
      \nml{XV.} &
      \(\amn{3^2\cdot 7\cdot 13\cdot 41\cdot 163\cdot 5\cdot
        977}{3^2\cdot 7\cdot 13\cdot 41\cdot 163\cdot 5867}\)\\
      \nml{XVI.} & \(\amn{2^3\cdot 17\cdot 79}{2^3\cdot 23\cdot 59}\)&
      \nml{XVII.} &
      \(\amn{2^4\cdot 23\cdot 1367}{2^4\cdot 53\cdot 607}\)\\
      \nml{XVIII.} & \(\amn{2^4\cdot 47\cdot 89}{2^4\cdot 53\cdot 79}\)&
      \nml{XIX.} &
      \(\amn{2^4\cdot 23\cdot 479}{2^4\cdot 89\cdot 127}\)\\
      \nml{XX.} & \(\amn{2^4\cdot 23\cdot 467}{2^4\cdot 103\cdot 107}\)&
      \nml{XXI.} &
      \(\amn{2^4\cdot 17\cdot 5119}{2^4\cdot 239\cdot 383}\)\\
      \nml{XXII.} & \(\amn{2^4\cdot 17\cdot 10303}{2^4\cdot 167\cdot 1103}\)&
      \nml{XXIII.} &
      \(\amn{2^4\cdot 19\cdot 1439}{2^4\cdot 149\cdot 191}\)\\
      \nml{XXIV.} & \(\amn{2^5\cdot 59\cdot 1103}{2^5\cdot 79\cdot 827}\)&
      \nml{XXV.} &
      \(\amn{2^5\cdot 37\cdot 12671}{2^5\cdot 227\cdot 2111}\)\\
      \nml{XXVI.} & \(\amn{2^5\cdot 53\cdot 10559}{2^5\cdot 79\cdot 7127}\)&
      \nml{XXVII.} &
      \(\amn{2^6\cdot 79\cdot 11087}{2^6\cdot 383\cdot 2309}\)\\
      \nml{XXVIII.} & \(\amn{2^8\cdot 383\cdot 9203}{2^8\cdot 1151\cdot 3067}\)&
      \nml{XXIX.} &
      \(\amn{2^2\cdot 11\cdot 17\cdot 263}{2^2\cdot 11\cdot
        43\cdot 107}\)\\
      \nml{XXX.} & \(\amn{3^3\cdot 5\cdot 7\cdot 71}{3^3\cdot 5\cdot 17\cdot 31}\)&
      \nml{XXXI.} &
      \(\amn{3^2\cdot 5\cdot 13\cdot 29\cdot 79}{3^2\cdot 5\cdot
        13\cdot 11\cdot 199}\)\\
      \nml{XXXII.} & \(\amn{3^2\cdot 5\cdot 13\cdot 19\cdot 47}{3^2\cdot 5\cdot 13\cdot 29\cdot 31}\)&
      \nml{XXXIII.} &
      \(\amn{3^2\cdot 5\cdot 13\cdot 19\cdot 37\cdot
        1583}{3^2\cdot 5\cdot 13\cdot 19\cdot 227\cdot 263}\) \\
      \nml{XXXIV.\footnote{Euler counted \(220499\) as a prime
          number. However even if it were prime, these would
          nevertheless fail to be amicable numbers. Indeed, we
          would have
          \(\int 11\cdot 220499=2646000=\int 89\cdot 29399\), but
          the values
          \[\int 3^2\cdot 7^2\cdot 13\cdot 19\cdot \int 11\cdot
            220499=548992080000\] and
          \[3^2\cdot 7^2\cdot 13\cdot 19(11\cdot 220499+89\cdot
            29399)=549209934000\] do not agree (\S 22). But, in
          fact, \(220499=311\cdot 709\). For this
          reason this pair XXXIV is to be deleted. \FR}} & \(\amn{3^2\cdot 7^2\cdot 13\cdot 19\cdot 11\cdot 220499}{3^2\cdot 7^2\cdot 13\cdot 19\cdot 89\cdot 29399}\)&
      \nml{XXXV.} &
      \(\amn{3^2\cdot 5\cdot 19\cdot 37\cdot 47}{3^2\cdot 5\cdot
        19\cdot 7\cdot 227}\)\\
      \nml{XXXVI.} & \(\amn{2^4\cdot 67\cdot 37\cdot 2411}{2^4\cdot 67\cdot 227\cdot 401}\)&
      \nml{XXXVII.}\footnote{In the first edition (and even in
        \cite{eulerfuss}) it reads
        \(3^2\cdot 5\cdot 7\cdot 11\cdot 29\) and
        \(3^3\cdot 5\cdot 31\cdot 89\). However, these numbers are
        not amicable. Indeed we have
        \(\int 7\cdot 11\cdot 29=2880=\int 31\cdot 89\), but the
        values
        \(\int 3^3\cdot 5\cdot \int 7\cdot 11\cdot 29=691200\) and
        \(3^3\cdot 5(7\cdot 11\cdot 29+31\cdot 89)=673920\) do not
        agree. However, from the equation
        \(z(7\cdot 11\cdot 29+31\cdot 89)=\int z\cdot \int 7\cdot
        11\cdot 29\) or
        \[\frac{z}{\int z}=\frac{2880}{4992}=\frac{3\cdot
            5}{2\cdot 13}=\frac{3\cdot 5}{2\cdot 13}
          \efrac{5}{6}\frac{3^2}{13} \efrac{3^2}{13}\] and it is
        found that \(z=3^2\cdot 5\). \FR} &
      \(\amn{3^2\cdot 5\cdot 7\cdot 11\cdot 29}{3^2\cdot 5\cdot
        31\cdot 89}\)\\
      \nml{XXXVIII.} & \(\amn{2\cdot 5\cdot 23\cdot 29\cdot 673}{2\cdot 5\cdot 7\cdot 60659}\)&
      \nml{XXXIX.} &
      \(\amn{2\cdot 5\cdot 7\cdot 19\cdot 107}{2\cdot 5\cdot
        47\cdot 359}\)\\
      \nml{XL.} & \(\amn{2^3\cdot 11\cdot 163\cdot 191}{2^3\cdot 31\cdot 11807}\)&
      \nml{XLI.} &
      \(\amn{3^2\cdot 7\cdot 13\cdot 23\cdot 11\cdot 19\cdot
        367}{3^2\cdot 7\cdot 13\cdot 23\cdot 79\cdot 1103}\)\\
      \nml{XLII.} & \(\amn{3^3\cdot 5\cdot 23\cdot 11\cdot 19\cdot 367}{3^3\cdot 5\cdot 23\cdot 79\cdot 1103}\)&
      \nml{XLIII.\footnote{In the first edition (and even in
          \cite{eulerfuss}) it reads \(57\) in place of
          \(47\). However this pair XLIII is the same as the pair
          XXVIII in the table of \cite{euler100}. It is clear that
          the incorrect number \(57\) has merely originated with a
          typographical error. \FR}} &
      \(\amn{2^3\cdot 11\cdot 59\cdot 173}{2^3\cdot 47\cdot
        2609}\)\\
      \nml{XLIV.} & \(\amn{2^3\cdot 11\cdot 23\cdot 2543}{2^3\cdot 383\cdot 1907}\)&
      \nml{XLV.} &
      \(\amn{2^3\cdot 11\cdot 23\cdot 1871}{2^3\cdot 467\cdot
        1151}\)\\
      \nml{XLVI.} & \(\amn{2^3\cdot 11\cdot 23\cdot 1619}{2^3\cdot 719\cdot 647}\)&
      \nml{XLVII.} &
      \(\amn{2^3\cdot 11\cdot 29\cdot 239}{2^3\cdot 191\cdot
        449}\)\\
      \nml{XLVIII.} & \(\amn{2^3\cdot 29\cdot 47\cdot 59}{2^3\cdot 17\cdot 4799}\)&
      \nml{XLIX.} &
      \(\amn{2^4\cdot 17\cdot 167\cdot 13679}{2^4\cdot 809\cdot
        51071}\)\\
      \nml{L.} & \(\amn{2^4\cdot 23\cdot 47\cdot 9767}{2^4\cdot 1583\cdot 7103}\) &
      \nml{LI.} &
      \(\amn{2^2\cdot 5\cdot 13\cdot 1187}{2^2\cdot 43\cdot
        2267}\)\\
      \nml{LII.} & \(\amn{3^2\cdot 7\cdot 13\cdot 5\cdot 17\cdot 1187}{3^2\cdot 7\cdot 13\cdot 131\cdot 971}\)&
      \nml{LIII.} &
      \(\amn{3^5\cdot 7^2\cdot 13\cdot 53\cdot 11\cdot
        211}{3^5\cdot 7^2\cdot 13\cdot 53\cdot 2543}\)\\
      \nml{LIV.} & \(\amn{3^2\cdot 5^2\cdot 11\cdot 59\cdot 179}{3^2\cdot 5^2\cdot 17\cdot 19\cdot 359}\)&
      \nml{LV.} &
      \(\amn{3^3\cdot 5\cdot 17\cdot 23\cdot 397}{3^3\cdot 5\cdot
        7\cdot 21491}\)\\
      \nml{LVI.} & \(\amn{3^4\cdot 7\cdot 11^2\cdot 19\cdot 47\cdot 7019}{3^4\cdot 7\cdot 11^2\cdot 19\cdot 389\cdot 863}\)&
      \nml{LVII.} &
      \(\amn{3^4\cdot 7\cdot 11^2\cdot 19\cdot 53\cdot
        6959}{3^4\cdot 7\cdot 11^2\cdot 19\cdot 179\cdot 2087}\)\\
      \nml{LVIII.} & \(\amn{3^5\cdot 7^2\cdot 13\cdot 19\cdot 47\cdot 7019}{3^5\cdot 7^2\cdot 13\cdot 19\cdot 389\cdot 863}\)&
      \nml{LIX.} &
      \(\amn{3^5\cdot 7^2\cdot 13\cdot 19\cdot 53\cdot
        6959}{3^5\cdot 7^2\cdot 13\cdot 19\cdot 179\cdot
        2087.}\)
    \end{longtable}
  \end{center}
  \egroup
  
  To this it is agreeable\footnote{However, it is also agreeable
    to add pairs VIII and IX, which are found in the table of
    \cite{euler100}, not to mention in \S 68, \S 78, \S 113 of
    this article. Indeed, those four which P. H. Fuss mentions
    in the preface to \cite{eulerfuss} (p. XXVI and LXXXI)
    reduce to these two, because pair XIII in that table is not
    valid and pair XXVIII coincides with pair XLIII of this
    table. In summary, therefore, Euler added \(59\) new pairs
    of amicable numbers to the three known before, if the error
    corrected in the footnote to pair XXXVII is to be treated as
    typographical. \FR} to add the following two pairs, which
  have a different form from the preceeding examples,
  \[\mbox{LX. }\amn{2^3\cdot 19\cdot 41}{2^5\cdot
      199}\qquad\mbox{LXI. }\amn{2^3\cdot 41\cdot 467}{2^5\cdot
      19\cdot 233.}\]
\end{Catalogue}

\input{euler}
\newpage\appendix

\begin{center}
  \medskip

  {\Huge\sc
    Appendices}

  \medskip
\end{center}

\section{Comments on \S 81}
\label{app:discussion}

In \S 81, Euler mentions in passing that, if \(f\) is prime, he
can abandon his search once he gets above a certain
threshhold. To see why, first observe that, by relabelling, we
can assume \(f\) to be the smallest of the prime factors \(p\),
\(q\), \(r\), \(f\). Now we will show that, as Euler says, if
\(f\) is taken sufficiently large, either \(p=hx-1\) or
\(q=gy-1\) will be smaller than \(f\).

We will use the notation from \S 73 and establish the claim in
general. Recall that the numbers \(a\), \(b\), and \(c\) are
fixed, and we have \(gh=\int f=f+1\),
\begin{gather*}e=bf-bgh+cgh=cf+c-b=\mathcal{O}(f),\\
\mbox{and}\quad
PQ=bbgh+be(f-1)=\mathcal{O}(f^2).
\end{gather*}
For convenience we denote the value of the product \(PQ\) by
\(R\). We will assume that \(f\) is large enough that all these
quantities (and anything we multiply or divide our inequalities
by, like \(e-b\)) is positive.

The numbers \(x\) and \(y\) are given by
\[x=\frac{P+bg}{e},\qquad y=\frac{Q+bh}{e},\] so \(p=hx-1\) and
\(q=gy-1\) are given by
\[p=\frac{Ph+b(f+1)}{e}-1,\qquad q=\frac{Qg+b(f+1)}{e}-1.\]
Since \(f<p\) and using the fact that \(gh=f+1\), we get
\[f<\frac{(f+1)(P/g+b)}{e}-1,\]
which implies
\[e-b<P/g,\tx{or} g<\frac{P}{e-b}.\]
Since \(f<q\) and using the fact that \(PQ=R\), we get
\[f<\frac{Rg/P+b(f+1)}{e}-1,\]
which implies
\[(e-b)(f+1)<Rg/P\tx{or}\frac{(e-b)(f+1)}{R}P<g.\]
Thus
\[\frac{(e-b)(f+1)}{R}P<g<\frac{P}{e-b}\]
and
\[(e-b)^2(f+1)<R.\] But the left-hand side is cubic in \(f\),
whilst \(R\) is only quadratic, so for sufficiently large \(f\)
the inequality fails.

In the case \(a=b=4\), \(c=1\), the inequality becomes
\[f^3-17f^2+35f+21<0,\] which breaks down around \(f=14.4833\)
which lies between the primes \(f=13\) and \(f=17\), just where
Euler stops.

\section{Code}
\label{app:code}

Here is the code that was used to recreate Euler's tables of
divisor sums of prime powers (this code specifically creates the
third table starting from \(193\)).

\begin{verbatim}from sage.all import sigma, latex, Primes

def cells(p,n):
    a = ''.join(['\(', str(p), '^{', str(n), '}\)'])
    b = latex(sigma(p**n,1).factor())
    b = b.replace('*', '\cdot')
    b = '\(' + b + '\)'
    return ''.join([a, ' & ', b])

def row(pns):
    if pns == ['h','h','h']:
        return '\hline'
    else:
        row_strings = [cells(*pn) for pn in pns]
    return ' & '.join(row_strings)+'\\\\'

P = Primes()[:168]

col_1 = [(p,k) for p in P if 192 < p and p < 252 for k in range(1,5)]
col_2 = [(p,k) for p in P if 256 < p and p < 314 for k in range(1,5)]
col_3 = [(p,k) for p in P if 316 < p and p < 384 for k in range(1,5)]

rows = [[a,b,c]
        for a,b,c in zip(col_1,
                         col_2,
                         col_3)]
for rw in rows:
    print(row(rw))
\end{verbatim}
\end{document}

%% file: bigtab.tex
\footnotesize
\begin{center}
\begin{tabular}{||p{0.6cm}|p{3.9cm}||p{0.6cm}|p{2.4cm}||p{0.6cm}|p{2.5cm}||}
  \hline\hline
  \scriptsize Num. & \scriptsize Divisor sum & \scriptsize Num. & \scriptsize Divisor sum & \scriptsize Num. & \scriptsize Divisor sum \\
  \hline\hline
\(2\) & \(3\) & \(3\) & \(2^{2}\) & \(11\) & \(2^{2} \cdot 3\)\\
\(2^{2}\) & \(7\) & \(3^{2}\) & \(13\) & \(11^{2}\) & \(7 \cdot 19\)\\
\(2^{3}\) & \(3 \cdot 5\) & \(3^{3}\) & \(2^{3} \cdot 5\) & \(11^{3}\) & \(2^{3} \cdot 3 \cdot 61\)\\
\(2^{4}\) & \(31\) & \(3^{4}\) & \(11^{2}\) & \(11^{4}\) & \(5 \cdot 3221\)\\
\(2^{5}\) & \(3^{2} \cdot 7\) & \(3^{5}\) & \(2^{2} \cdot 7 \cdot 13\) & \(11^{5}\) & \(2^{2} \cdot 3^{2} \cdot 7 \cdot 19 \cdot 37\)\\                                                                                                                                  
\(2^{6}\) & \(127\) & \(3^{6}\) & \(1093\) & \(11^{6}\) & \(43 \cdot 45319\)\\
\(2^{7}\) & \(3 \cdot 5 \cdot 17\) & \(3^{7}\) & \(2^{4} \cdot 5 \cdot 41\) & \(11^{7}\) & \(2^{4} \cdot 3 \cdot 61 \cdot 7321\)\\
\(2^{8}\) & \(7 \cdot 73\) & \(3^{8}\) & \(13 \cdot 757\) & \(11^{8}\) & \(7 \cdot 19 \cdot 1772893\)\\
  \(2^{9}\) & \(3 \cdot 11 \cdot 31\) & \(3^{9}\) & \(2^{2} \cdot 11^{2} \cdot 61\) & \(11^{9}\) & \(2^{2} \cdot 3 \cdot 5 \cdot 3221 \cdot 13421\)\\
  \cline{5-6}
\(2^{10}\) & \(23 \cdot 89\) & \(3^{10}\) & \(23 \cdot 3851\) & \(13\) & \(2 \cdot 7\)\\
\(2^{11}\) & \(3^{2} \cdot 5 \cdot 7 \cdot 13\) & \(3^{11}\) & \(2^{3} \cdot 5 \cdot 7 \cdot 13 \cdot 73\) & \(13^{2}\) & \(3 \cdot 61\)\\                                                                                                                               
\(2^{12}\) & \(8191\) & \(3^{12}\) & \(797161\) & \(13^{3}\) & \(2^{2} \cdot 5 \cdot 7 \cdot 17\)\\
\(2^{13}\) & \(3 \cdot 43 \cdot 127\) & \(3^{13}\) & \(2^{2} \cdot 547 \cdot 1093\) & \(13^{4}\) & \(30941\)\\
\(2^{14}\) & \(7 \cdot 31 \cdot 151\) & \(3^{14}\) & \(11^{2} \cdot 13 \cdot 4561\) & \(13^{5}\) & \(2 \cdot 3 \cdot 7 \cdot 61 \cdot 157\)\\                                                                                                                            
  \(2^{15}\) & \(3 \cdot 5 \cdot 17 \cdot 257\) & \(3^{15}\) & \(2^{5} \cdot 5 \cdot 17 \cdot 41 \cdot 193\) & \(13^{6}\) & \(5229043\)\\
  \cline{3-4}
  \(2^{16}\) & \(131071\) & & & \(13^{7}\) & \(2^{3} \cdot 5 \cdot 7 \cdot 17 \cdot 14281\)\\
  \cline{5-6}
\(2^{17}\) & \(3^{3} \cdot 7 \cdot 19 \cdot 73\) & \(5\) & \(2 \cdot 3\) &  & \\
\(2^{18}\) & \(524287\) & \(5^{2}\) & \(31\) & \(17\) & \(2 \cdot 3^{2}\)\\
\(2^{19}\) & \(3 \cdot 5^{2} \cdot 11 \cdot 31 \cdot 41\) & \(5^{3}\) & \(2^{2} \cdot 3 \cdot 13\) & \(17^{2}\) & \(307\)\\
\(2^{20}\) & \(7^{2} \cdot 127 \cdot 337\) & \(5^{4}\) & \(11 \cdot 71\) & \(17^{3}\) & \(2^{2} \cdot 3^{2} \cdot 5 \cdot 29\)\\
\(2^{21}\) & \(3 \cdot 23 \cdot 89 \cdot 683\) & \(5^{5}\) & \(2 \cdot 3^{2} \cdot 7 \cdot 31\) & \(17^{4}\) & \(88741\)\\
  \(2^{22}\) & \(47 \cdot 178481\) & \(5^{6}\) & \(19531\) & \(17^{5}\) & \(2 \cdot 3^{3} \cdot 7 \cdot 13 \cdot 307\)\\
  \cline{5-6}
\(2^{23}\) & \(3^{2} \cdot 5 \cdot 7 \cdot 13 \cdot 17 \cdot 241\) & \(5^{7}\) & \(2^{3} \cdot 3 \cdot 13 \cdot 313\) & \(19\) & 
\(2^{2} \cdot 5\)\\                                                                                                                 
\(2^{24}\) & \(31 \cdot 601 \cdot 1801\) & \(5^{8}\) & \(19 \cdot 31 \cdot 829\) & \(19^{2}\) & \(3 \cdot 127\)\\
  \(2^{25}\) & \(3 \cdot 2731 \cdot 8191\) & \(5^{9}\) & \(2 \cdot 3 \cdot 11 \cdot 71 \cdot 521\) & \(19^{3}\) & \(2^{3} \cdot 5 \cdot 181\)\\
  \cline{3-4}
\(2^{26}\) & \(7 \cdot 73 \cdot 262657\) &  &  & \(19^{4}\) & \(151 \cdot 911\)\\
  \(2^{27}\) & \(3 \cdot 5 \cdot 29 \cdot 43 \cdot 113 \cdot 127\) & \(7\) & \(2^{3}\) & \(19^{5}\) & \(2^{2} \cdot 3 \cdot 5 \cdot 7^{3} \cdot 127\)\\
  \cline{5-6}
\(2^{28}\) & \(233 \cdot 1103 \cdot 2089\) & \(7^{2}\) & \(3 \cdot 19\) &  & \\
\(2^{29}\) & \(3^{2} \cdot 7 \cdot 11 \cdot 31 \cdot 151 \cdot 331\) & \(7^{3}\) & \(2^{4} \cdot 5^{2}\) & \(23\) & \(2^{3} \cdot 3\)\\                                                                                                                              
\(2^{30}\) & \(2147483647\) & \(7^{4}\) & \(2801\) & \(23^{2}\) & \(7 \cdot 79\)\\
\(2^{31}\) & \(3 \cdot 5 \cdot 17 \cdot 257 \cdot 65537\) & \(7^{5}\) & \(2^{3} \cdot 3 \cdot 19 \cdot 43\) & \(23^{3}\) & \(2^{4} \cdot 3 \cdot 5 \cdot 53\)\\                                                                                                          
  \(2^{32}\) & \(7 \cdot 23 \cdot 89 \cdot 599479\) & \(7^{6}\) & \(29 \cdot 4733\) & \(23^{4}\) & \(292561\)\\
  \cline{5-6}
\(2^{33}\) & \(3 \cdot 43691 \cdot 131071\) & \(7^{7}\) & \(2^{5} \cdot 5^{2} \cdot 1201\) & &\\
\(2^{34}\) & \(31 \cdot 71 \cdot 127 \cdot 122921\) & \(7^{8}\) & \(3^{2} \cdot 19 \cdot 37 \cdot 1063\) & \(29\) & \(2 \cdot 3 \cdot 5\)\\                                                                                                                          
\(2^{35}\) & \(3^{3} \cdot 5 \cdot 7 \cdot 13 \cdot 19 \cdot 37 \cdot 73 \cdot 109\) & \(7^{9}\) & \(2^{3} \cdot 11 \cdot 191 \cdot 2801\) & \(29^{2}\) & \(13 \cdot 67\)\\                                                                                              
\(2^{36}\) & \(223 \cdot 616318177\) & \(7^{10}\) & \(1123 \cdot 293459\) & \(29^{3}\) & \(2^{2} \cdot 3 \cdot 5 \cdot 421\)\\
  \hline\hline
\end{tabular}

\begin{tabular}{||p{0.7cm}|p{2.5cm}||p{0.7cm}|p{2.5cm}||p{0.7cm}|p{2.5cm}||}
  \hline\hline
  Num. & Divisor sum & Num. & Divisor sum & Num. & Divisor sum \\
  \hline\hline
\(31\) & \(2^{5}\) & \(79\) & \(2^{4} \cdot 5\) & \(137\) & \(2 \cdot 3 \cdot 23\)\\
\(31^{2}\) & \(3 \cdot 331\) & \(79^{2}\) & \(3 \cdot 7^{2} \cdot 43\) & \(137^{2}\) & \(7 \cdot 37 \cdot 73\)\\
\(31^{3}\) & \(2^{6} \cdot 13 \cdot 37\) & \(79^{3}\) & \(2^{5} \cdot 5 \cdot 3121\) & \(137^{3}\) & \(2^{2} \cdot 3 \cdot 5 \cdot 23 \cdot 1877\)\\                                                                                                                     
\hline
\(37\) & \(2 \cdot 19\) & \(83\) & \(2^{2} \cdot 3 \cdot 7\) & \(139\) & \(2^{2} \cdot 5 \cdot 7\)\\
\(37^{2}\) & \(3 \cdot 7 \cdot 67\) & \(83^{2}\) & \(19 \cdot 367\) & \(139^{2}\) & \(3 \cdot 13 \cdot 499\)\\
\(37^{3}\) & \(2^{2} \cdot 5 \cdot 19 \cdot 137\) & \(83^{3}\) & \(2^{3} \cdot 3 \cdot 5 \cdot 7 \cdot 13 \cdot 53\) & \(139^{3}\) & 
\(2^{3} \cdot 5 \cdot 7 \cdot 9661\)\\                                                                                              
\hline
\(41\) & \(2 \cdot 3 \cdot 7\) & \(89\) & \(2 \cdot 3^{2} \cdot 5\) & \(149\) & \(2 \cdot 3 \cdot 5^{2}\)\\
\(41^{2}\) & \(1723\) & \(89^{2}\) & \(8011\) & \(149^{2}\) & \(7 \cdot 31 \cdot 103\)\\
\(41^{3}\) & \(2^{2} \cdot 3 \cdot 7 \cdot 29^{2}\) & \(89^{3}\) & \(2^{2} \cdot 3^{2} \cdot 5 \cdot 17 \cdot 233\) & \(149^{3}\) & \(2^{2} \cdot 3 \cdot 5^{2} \cdot 17 \cdot 653\)\\                                                                                   
\hline
\(43\) & \(2^{2} \cdot 11\) & \(97\) & \(2 \cdot 7^{2}\) & \(151\) & \(2^{3} \cdot 19\)\\
\(43^{2}\) & \(3 \cdot 631\) & \(97^{2}\) & \(3 \cdot 3169\) & \(151^{2}\) & \(3 \cdot 7 \cdot 1093\)\\
\(43^{3}\) & \(2^{3} \cdot 5^{2} \cdot 11 \cdot 37\) & \(97^{3}\) & \(2^{2} \cdot 5 \cdot 7^{2} \cdot 941\) & \(151^{3}\) & \(2^{4} \cdot 13 \cdot 19 \cdot 877\)\\                                                                                                      
\hline
\(47\) & \(2^{4} \cdot 3\) & \(101\) & \(2 \cdot 3 \cdot 17\) & \(157\) & \(2 \cdot 79\)\\
\(47^{2}\) & \(37 \cdot 61\) & \(101^{2}\) & \(10303\) & \(157^{2}\) & \(3 \cdot 8269\)\\
\(47^{3}\) & \(2^{5} \cdot 3 \cdot 5 \cdot 13 \cdot 17\) & \(101^{3}\) & \(2^{2} \cdot 3 \cdot 17 \cdot 5101\) & \(157^{3}\) & \(2^{2} \cdot 5^{2} \cdot 17 \cdot 29 \cdot 79\)\\                                                                                        
\hline
\(53\) & \(2 \cdot 3^{3}\) & \(103\) & \(2^{3} \cdot 13\) & \(163\) & \(2^{2} \cdot 41\)\\
\(53^{2}\) & \(7 \cdot 409\) & \(103^{2}\) & \(3 \cdot 3571\) & \(163^{2}\) & \(3 \cdot 7 \cdot 19 \cdot 67\)\\
\(53^{3}\) & \(2^{2} \cdot 3^{3} \cdot 5 \cdot 281\) & \(103^{3}\) & \(2^{4} \cdot 5 \cdot 13 \cdot 1061\) & \(163^{3}\) & \(2^{3} \cdot 5 \cdot 41 \cdot 2657\)\\                                                                                                       
\hline
\(59\) & \(2^{2} \cdot 3 \cdot 5\) & \(107\) & \(2^{2} \cdot 3^{3}\) & \(167\) & \(2^{3} \cdot 3 \cdot 7\)\\
\(59^{2}\) & \(3541\) & \(107^{2}\) & \(7 \cdot 13 \cdot 127\) & \(167^{2}\) & \(28057\)\\
\(59^{3}\) & \(2^{3} \cdot 3 \cdot 5 \cdot 1741\) & \(107^{3}\) & \(2^{3} \cdot 3^{3} \cdot 5^{2} \cdot 229\) & \(167^{3}\) & \(2^{4} \cdot 3 \cdot 5 \cdot 7 \cdot 2789\)\\                                                                                             
\hline
\(61\) & \(2 \cdot 31\) & \(109\) & \(2 \cdot 5 \cdot 11\) & \(173\) & \(2 \cdot 3 \cdot 29\)\\
\(61^{2}\) & \(3 \cdot 13 \cdot 97\) & \(109^{2}\) & \(3 \cdot 7 \cdot 571\) & \(173^{2}\) & \(30103\)\\
\(61^{3}\) & \(2^{2} \cdot 31 \cdot 1861\) & \(109^{3}\) & \(2^{2} \cdot 5 \cdot 11 \cdot 13 \cdot 457\) & \(173^{3}\) & \(2^{2} \cdot 3 \cdot 5 \cdot 29 \cdot 41 \cdot 73\)\\                                                                                          
\hline
\(67\) & \(2^{2} \cdot 17\) & \(113\) & \(2 \cdot 3 \cdot 19\) & \(179\) & \(2^{2} \cdot 3^{2} \cdot 5\)\\
\(67^{2}\) & \(3 \cdot 7^{2} \cdot 31\) & \(113^{2}\) & \(13 \cdot 991\) & \(179^{2}\) & \(7 \cdot 4603\)\\
\(67^{3}\) & \(2^{3} \cdot 5 \cdot 17 \cdot 449\) & \(113^{3}\) & \(2^{2} \cdot 3 \cdot 5 \cdot 19 \cdot 1277\) & \(179^{3}\) & \(2^{3} \cdot 3^{2} \cdot 5 \cdot 37 \cdot 433\)\\                                                                                       
\hline
\(71\) & \(2^{3} \cdot 3^{2}\) & \(127\) & \(2^{7}\) & \(181\) & \(2 \cdot 7 \cdot 13\)\\
\(71^{2}\) & \(5113\) & \(127^{2}\) & \(3 \cdot 5419\) & \(181^{2}\) & \(3 \cdot 79 \cdot 139\)\\
\(71^{3}\) & \(2^{4} \cdot 3^{2} \cdot 2521\) & \(127^{3}\) & \(2^{8} \cdot 5 \cdot 1613\) & \(181^{3}\) & \(2^{2} \cdot 7 \cdot 13 \cdot 16381\)\\                                                                                                                      
\hline
\(73\) & \(2 \cdot 37\) & \(131\) & \(2^{2} \cdot 3 \cdot 11\) & \(191\) & \(2^{6} \cdot 3\)\\
\(73^{2}\) & \(3 \cdot 1801\) & \(131^{2}\) & \(17293\) & \(191^{2}\) & \(7 \cdot 13^{2} \cdot 31\)\\
\(73^{3}\) & \(2^{2} \cdot 5 \cdot 13 \cdot 37 \cdot 41\) & \(131^{3}\) & \(2^{3} \cdot 3 \cdot 11 \cdot 8581\) & \(191^{3}\) & \(2^{7} \cdot 3 \cdot 17 \cdot 29 \cdot 37\)\\
\hline\hline
\end{tabular}

\begin{tabular}{||p{0.7cm}|p{2.5cm}||p{0.7cm}|p{2.5cm}||p{0.7cm}|p{2.5cm}||}
  \hline\hline
  Num. & Divisor sum & Num. & Divisor sum & Num. & Divisor sum \\
  \hline\hline
  \(193\) & \(2 \cdot 97\) & \(257\) & \(2 \cdot 3 \cdot 43\) & \(317\) & \(2 \cdot 3 \cdot 53\)\\
\(193^{2}\) & \(3 \cdot 7 \cdot 1783\) & \(257^{2}\) & \(61 \cdot 1087\) & \(317^{2}\) & \(7 \cdot 14401\)\\
\(193^{3}\) & \(2^{2} \cdot 5^{3} \cdot 97 \cdot 149\) & \(257^{3}\) & \(2^{2} \cdot 3 \cdot 5^{2} \cdot 43 \cdot 1321\) & \(317^{3}\) & \(2^{2} \cdot 3 \cdot 5 \cdot 13 \cdot 53 \cdot 773\)\\                                                                         
\hline
\(197\) & \(2 \cdot 3^{2} \cdot 11\) & \(263\) & \(2^{3} \cdot 3 \cdot 11\) & \(331\) & \(2^{2} \cdot 83\)\\
\(197^{2}\) & \(19 \cdot 2053\) & \(263^{2}\) & \(7^{2} \cdot 13 \cdot 109\) & \(331^{2}\) & \(3 \cdot 7 \cdot 5233\)\\
\(197^{3}\) & \(2^{2} \cdot 3^{2} \cdot 5 \cdot 11 \cdot 3881\) & \(263^{3}\) & \(2^{4} \cdot 3 \cdot 5 \cdot 11 \cdot 6917\) & \(331^{3}\) & \(2^{3} \cdot 29 \cdot 83 \cdot 1889\)\\                                                                                   
\hline
\(199\) & \(2^{3} \cdot 5^{2}\) & \(269\) & \(2 \cdot 3^{3} \cdot 5\) & \(337\) & \(2 \cdot 13^{2}\)\\
\(199^{2}\) & \(3 \cdot 13267\) & \(269^{2}\) & \(13 \cdot 37 \cdot 151\) & \(337^{2}\) & \(3 \cdot 43 \cdot 883\)\\
\(199^{3}\) & \(2^{4} \cdot 5^{2} \cdot 19801\) & \(269^{3}\) & \(2^{2} \cdot 3^{3} \cdot 5 \cdot 97 \cdot 373\) & \(337^{3}\) & \(2^{2} \cdot 5 \cdot 13^{2} \cdot 41 \cdot 277\)\\                                                                                     
\hline
\(211\) & \(2^{2} \cdot 53\) & \(271\) & \(2^{4} \cdot 17\) & \(347\) & \(2^{2} \cdot 3 \cdot 29\)\\
\(211^{2}\) & \(3 \cdot 13 \cdot 31 \cdot 37\) & \(271^{2}\) & \(3 \cdot 24571\) & \(347^{2}\) & \(7 \cdot 13 \cdot 1327\)\\
\(211^{3}\) & \(2^{3} \cdot 53 \cdot 113 \cdot 197\) & \(271^{3}\) & \(2^{5} \cdot 17 \cdot 36721\) & \(347^{3}\) & \(2^{3} \cdot 3 \cdot 5 \cdot 29 \cdot 12041\)\\                                                                                                     
\hline
\(223\) & \(2^{5} \cdot 7\) & \(277\) & \(2 \cdot 139\) & \(349\) & \(2 \cdot 5^{2} \cdot 7\)\\
\(223^{2}\) & \(3 \cdot 16651\) & \(277^{2}\) & \(3 \cdot 7 \cdot 19 \cdot 193\) & \(349^{2}\) & \(3 \cdot 19 \cdot 2143\)\\
\(223^{3}\) & \(2^{6} \cdot 5 \cdot 7 \cdot 4973\) & \(277^{3}\) & \(2^{2} \cdot 5 \cdot 139 \cdot 7673\) & \(349^{3}\) & \(2^{2} \cdot 5^{2} \cdot 7 \cdot 60901\)\\                                                                                                    
\hline
\(227\) & \(2^{2} \cdot 3 \cdot 19\) & \(281\) & \(2 \cdot 3 \cdot 47\) & \(353\) & \(2 \cdot 3 \cdot 59\)\\
\(227^{2}\) & \(73 \cdot 709\) & \(281^{2}\) & \(109 \cdot 727\) & \(353^{2}\) & \(19 \cdot 6577\)\\
\(227^{3}\) & \(2^{3} \cdot 3 \cdot 5 \cdot 19 \cdot 5153\) & \(281^{3}\) & \(2^{2} \cdot 3 \cdot 13 \cdot 47 \cdot 3037\) & \(353^{3}\) & \(2^{2} \cdot 3 \cdot 5 \cdot 17 \cdot 59 \cdot 733\)\\                                                                       
\hline
\(229\) & \(2 \cdot 5 \cdot 23\) & \(283\) & \(2^{2} \cdot 71\) & \(359\) & \(2^{3} \cdot 3^{2} \cdot 5\)\\
\(229^{2}\) & \(3 \cdot 97 \cdot 181\) & \(283^{2}\) & \(3 \cdot 73 \cdot 367\) & \(359^{2}\) & \(7 \cdot 37 \cdot 499\)\\
\(229^{3}\) & \(2^{2} \cdot 5 \cdot 13 \cdot 23 \cdot 2017\) & \(283^{3}\) & \(2^{3} \cdot 5 \cdot 71 \cdot 8009\) & \(359^{3}\) & \(2^{4} \cdot 3^{2} \cdot 5 \cdot 13 \cdot 4957\)\\                                                                                   
\hline
\(233\) & \(2 \cdot 3^{2} \cdot 13\) & \(293\) & \(2 \cdot 3 \cdot 7^{2}\) & \(367\) & \(2^{4} \cdot 23\)\\
\(233^{2}\) & \(7 \cdot 7789\) & \(293^{2}\) & \(86143\) & \(367^{2}\) & \(3 \cdot 13 \cdot 3463\)\\
\(233^{3}\) & \(2^{2} \cdot 3^{2} \cdot 5 \cdot 13 \cdot 61 \cdot 89\) & \(293^{3}\) & \(2^{2} \cdot 3 \cdot 5^{2} \cdot 7^{2} \cdot 17 \cdot 101\) & \(367^{3}\) & \(2^{5} \cdot 5 \cdot 23 \cdot 13469\)\\                                                             
\hline
\(239\) & \(2^{4} \cdot 3 \cdot 5\) & \(307\) & \(2^{2} \cdot 7 \cdot 11\) & \(373\) & \(2 \cdot 11 \cdot 17\)\\
\(239^{2}\) & \(19 \cdot 3019\) & \(307^{2}\) & \(3 \cdot 43 \cdot 733\) & \(373^{2}\) & \(3 \cdot 7^{2} \cdot 13 \cdot 73\)\\
\(239^{3}\) & \(2^{5} \cdot 3 \cdot 5 \cdot 13^{4}\) & \(307^{3}\) & \(2^{3} \cdot 5^{3} \cdot 7 \cdot 11 \cdot 13 \cdot 29\) & \(373^{3}\) & \(2^{2} \cdot 5 \cdot 11 \cdot 17 \cdot 13913\)\\                                                                          
\hline
\(241\) & \(2 \cdot 11^{2}\) & \(311\) & \(2^{3} \cdot 3 \cdot 13\) & \(379\) & \(2^{2} \cdot 5 \cdot 19\)\\
\(241^{2}\) & \(3 \cdot 19441\) & \(311^{2}\) & \(19 \cdot 5107\) & \(379^{2}\) & \(3 \cdot 61 \cdot 787\)\\
\(241^{3}\) & \(2^{2} \cdot 11^{2} \cdot 113 \cdot 257\) & \(311^{3}\) & \(2^{4} \cdot 3 \cdot 13 \cdot 137 \cdot 353\) & \(379^{3}\) & \(2^{3} \cdot 5 \cdot 19 \cdot 71821\)\\                                                                                         
\hline
\(251\) & \(2^{2} \cdot 3^{2} \cdot 7\) & \(313\) & \(2 \cdot 157\) & \(383\) & \(2^{7} \cdot 3\)\\
\(251^{2}\) & \(43 \cdot 1471\) & \(313^{2}\) & \(3 \cdot 181^{2}\) & \(383^{2}\) & \(147073\)\\
\(251^{3}\) & \(2^{3} \cdot 3^{2} \cdot 7 \cdot 17^{2} \cdot 109\) & \(313^{3}\) & \(2^{2} \cdot 5 \cdot 97 \cdot 101 \cdot 157\) & \(383^{3}\) & \(2^{8} \cdot 3 \cdot 5 \cdot 14669\)\\                                                                                
\hline\hline
\end{tabular}

\begin{tabular}{||p{0.7cm}|p{2.7cm}||p{0.7cm}|p{2.7cm}||p{0.7cm}|p{2.7cm}||}
  \hline\hline
  Num. & Divisor sum & Num. & Divisor sum & Num. & Divisor sum \\
  \hline\hline
\(389\) & \(2 \cdot 3 \cdot 5 \cdot 13\) & \(457\) & \(2 \cdot 229\) & \(523\) & \(2^{2} \cdot 131\)\\
\(389^{2}\) & \(7 \cdot 21673\) & \(457^{2}\) & \(3 \cdot 7 \cdot 9967\) & \(523^{2}\) & \(3 \cdot 13 \cdot 7027\)\\
\(389^{3}\) & \(2^{2} \cdot 3 \cdot 5 \cdot 13 \cdot 29 \cdot 2609\) & \(457^{3}\) & \(2^{2} \cdot 5^{2} \cdot 229 \cdot 4177\) & \(523^{3}\) & \(2^{3} \cdot 5 \cdot 17 \cdot 131 \cdot 1609\)\\                                                                        
\hline
\(397\) & \(2 \cdot 199\) & \(461\) & \(2 \cdot 3 \cdot 7 \cdot 11\) & \(541\) & \(2 \cdot 271\)\\
\(397^{2}\) & \(3 \cdot 31 \cdot 1699\) & \(461^{2}\) & \(373 \cdot 571\) & \(541^{2}\) & \(3 \cdot 7 \cdot 13963\)\\
\(397^{3}\) & \(2^{2} \cdot 5 \cdot 199 \cdot 15761\) & \(461^{3}\) & \(2^{2} \cdot 3 \cdot 7 \cdot 11 \cdot 106261\) & \(541^{3}\) & \(2^{2} \cdot 13 \cdot 271 \cdot 11257\)\\                                                                                         
\hline
\(401\) & \(2 \cdot 3 \cdot 67\) & \(463\) & \(2^{4} \cdot 29\) & \(547\) & \(2^{2} \cdot 137\)\\
\(401^{2}\) & \(7 \cdot 23029\) & \(463^{2}\) & \(3 \cdot 19 \cdot 3769\) & \(547^{2}\) & \(3 \cdot 163 \cdot 613\)\\
\(401^{3}\) & \(2^{2} \cdot 3 \cdot 37 \cdot 41 \cdot 53 \cdot 67\) & \(463^{3}\) & \(2^{5} \cdot 5 \cdot 13 \cdot 17 \cdot 29 \cdot 97\) & \(547^{3}\) & \(2^{3} \cdot 5 \cdot 137 \cdot 29921\)\\                                                                      
\hline
\(409\) & \(2 \cdot 5 \cdot 41\) & \(467\) & \(2^{2} \cdot 3^{2} \cdot 13\) & \(557\) & \(2 \cdot 3^{2} \cdot 31\)\\
\(409^{2}\) & \(3 \cdot 55897\) & \(467^{2}\) & \(19 \cdot 11503\) & \(557^{2}\) & \(7^{2} \cdot 6343\)\\
\(409^{3}\) & \(2^{2} \cdot 5 \cdot 41 \cdot 83641\) & \(467^{3}\) & \(2^{3} \cdot 3^{2} \cdot 5 \cdot 13 \cdot 113 \cdot 193\) & \(557^{3}\) & \(2^{2} \cdot 3^{2} \cdot 5^{3} \cdot 17 \cdot 31 \cdot 73\)\\                                                           
\hline
\(419\) & \(2^{2} \cdot 3 \cdot 5 \cdot 7\) & \(479\) & \(2^{5} \cdot 3 \cdot 5\) & \(563\) & \(2^{2} \cdot 3 \cdot 47\)\\                                                                                                                                   
\(419^{2}\) & \(13 \cdot 13537\) & \(479^{2}\) & \(43 \cdot 5347\) & \(563^{2}\) & \(31 \cdot 10243\)\\
\(419^{3}\) & \(2^{3} \cdot 3 \cdot 5 \cdot 7 \cdot 41 \cdot 2141\) & \(479^{3}\) & \(2^{6} \cdot 3 \cdot 5 \cdot 89 \cdot 1289\) & \(563^{3}\) & \(2^{3} \cdot 3 \cdot 5 \cdot 29 \cdot 47 \cdot 1093\)\\                                                               
\hline
\(421\) & \(2 \cdot 211\) & \(487\) & \(2^{3} \cdot 61\) & \(569\) & \(2 \cdot 3 \cdot 5 \cdot 19\)\\
\(421^{2}\) & \(3 \cdot 59221\) & \(487^{2}\) & \(3 \cdot 7 \cdot 11317\) & \(569^{2}\) & \(7^{2} \cdot 6619\)\\
\(421^{3}\) & \(2^{2} \cdot 13 \cdot 17 \cdot 211 \cdot 401\) & \(487^{3}\) & \(2^{4} \cdot 5 \cdot 37 \cdot 61 \cdot 641\) & \(569^{3}\) & \(2^{2} \cdot 3 \cdot 5 \cdot 19 \cdot 161881\)\\                                                                            
\hline
\(431\) & \(2^{4} \cdot 3^{3}\) & \(491\) & \(2^{2} \cdot 3 \cdot 41\) & \(571\) & \(2^{2} \cdot 11 \cdot 13\)\\
\(431^{2}\) & \(7 \cdot 67 \cdot 397\) & \(491^{2}\) & \(37 \cdot 6529\) & \(571^{2}\) & \(3 \cdot 7 \cdot 103 \cdot 151\)\\
\(431^{3}\) & \(2^{5} \cdot 3^{3} \cdot 293 \cdot 317\) & \(491^{3}\) & \(2^{3} \cdot 3 \cdot 41 \cdot 149 \cdot 809\) & \(571^{3}\) & \(2^{3} \cdot 11 \cdot 13 \cdot 163021\)\\                                                                                        
\hline
\(433\) & \(2 \cdot 7 \cdot 31\) & \(499\) & \(2^{2} \cdot 5^{3}\) & \(577\) & \(2 \cdot 17^{2}\)\\
\(433^{2}\) & \(3 \cdot 37 \cdot 1693\) & \(499^{2}\) & \(3 \cdot 7 \cdot 109^{2}\) & \(577^{2}\) & \(3 \cdot 19 \cdot 5851\)\\
\(433^{3}\) & \(2^{2} \cdot 5 \cdot 7 \cdot 31 \cdot 18749\) & \(499^{3}\) & \(2^{3} \cdot 5^{3} \cdot 13 \cdot 61 \cdot 157\) & \(577^{3}\) & \(2^{2} \cdot 5 \cdot 13^{2} \cdot 17^{2} \cdot 197\)\\                                                                   
\hline
\(439\) & \(2^{3} \cdot 5 \cdot 11\) & \(503\) & \(2^{3} \cdot 3^{2} \cdot 7\) & \(587\) & \(2^{2} \cdot 3 \cdot 7^{2}\)\\
\(439^{2}\) & \(3 \cdot 31^{2} \cdot 67\) & \(503^{2}\) & \(13 \cdot 19501\) & \(587^{2}\) & \(547 \cdot 631\)\\
\(439^{3}\) & \(2^{4} \cdot 5 \cdot 11 \cdot 173 \cdot 557\) & \(503^{3}\) & \(2^{4} \cdot 3^{2} \cdot 5 \cdot 7 \cdot 25301\) & \(587^{3}\) & \(2^{3} \cdot 3 \cdot 5 \cdot 7^{2} \cdot 34457\)\\                                                                       
\hline
\(443\) & \(2^{2} \cdot 3 \cdot 37\) & \(509\) & \(2 \cdot 3 \cdot 5 \cdot 17\) & \(593\) & \(2 \cdot 3^{3} \cdot 11\)\\
\(443^{2}\) & \(7 \cdot 28099\) & \(509^{2}\) & \(43 \cdot 6037\) & \(593^{2}\) & \(163 \cdot 2161\)\\
\(443^{3}\) & \(2^{3} \cdot 3 \cdot 5^{4} \cdot 37 \cdot 157\) & \(509^{3}\) & \(2^{2} \cdot 3 \cdot 5 \cdot 17 \cdot 281 \cdot 461\) & \(593^{3}\) & \(2^{2} \cdot 3^{3} \cdot 5^{2} \cdot 11 \cdot 13 \cdot 541\)\\                                                    
\hline
\(449\) & \(2 \cdot 3^{2} \cdot 5^{2}\) & \(521\) & \(2 \cdot 3^{2} \cdot 29\) & \(599\) & \(2^{3} \cdot 3 \cdot 5^{2}\)\\                                                                                                                                   
\(449^{2}\) & \(97 \cdot 2083\) & \(521^{2}\) & \(31^{2} \cdot 283\) & \(599^{2}\) & \(7 \cdot 51343\)\\
\(449^{3}\) & \(2^{2} \cdot 3^{2} \cdot 5^{2} \cdot 100801\) & \(521^{3}\) & \(2^{2} \cdot 3^{2} \cdot 29 \cdot 135721\) & \(599^{3}\) & \(2^{4} \cdot 3 \cdot 5^{2} \cdot 17 \cdot 61 \cdot 173\)\\                                                                     
\hline\hline
\end{tabular}

\begin{tabular}{||p{0.7cm}|p{2.7cm}||p{0.7cm}|p{2.7cm}||p{0.7cm}|p{2.7cm}||}
  \hline\hline
  Num. & Divisor sum & Num. & Divisor sum & Num. & Divisor sum \\
  \hline\hline
  \(601\) & \(2 \cdot 7 \cdot 43\) & \(661\) & \(2 \cdot 331\) & \(743\) & \(2^{3} \cdot 3 \cdot 31\)\\
  \(601^{2}\) & \(3 \cdot 13 \cdot 9277\) & \(661^{2}\) & \(3 \cdot 145861\) & \(743^{2}\) & \(552793\)\\
  \(601^{3}\) & \(2^{2} \cdot 7 \cdot 43 \cdot 313 \cdot 577\) & \(661^{3}\) & \(2^{2} \cdot 331 \cdot 218461\) & \(743^{3}\) & \(2^{4} \cdot 3 \cdot 5^{2} \cdot 31 \cdot 61 \cdot 181\)\\                                                                                
  \hline
  \(607\) & \(2^{5} \cdot 19\) & \(673\) & \(2 \cdot 337\) & \(751\) & \(2^{4} \cdot 47\)\\
  \(607^{2}\) & \(3 \cdot 13 \cdot 9463\) & \(673^{2}\) & \(3 \cdot 151201\) & \(751^{2}\) & \(3 \cdot 7 \cdot 26893\)\\
  \(607^{3}\) & \(2^{6} \cdot 5^{2} \cdot 19 \cdot 7369\) & \(673^{3}\) & \(2^{2} \cdot 5 \cdot 337 \cdot 45293\) & \(751^{3}\) & \(2^{5} \cdot 47 \cdot 282001\)\\                                                                                                        
  \hline
  \(613\) & \(2 \cdot 307\) & \(677\) & \(2 \cdot 3 \cdot 113\) & \(757\) & \(2 \cdot 379\)\\
  \(613^{2}\) & \(3 \cdot 7 \cdot 17923\) & \(677^{2}\) & \(459007\) & \(757^{2}\) & \(3 \cdot 13 \cdot 14713\)\\
  \(613^{3}\) & \(2^{2} \cdot 5 \cdot 53 \cdot 307 \cdot 709\) & \(677^{3}\) & \(2^{2} \cdot 3 \cdot 5 \cdot 113 \cdot 45833\) & \(757^{3}\) & \(2^{2} \cdot 5^{2} \cdot 73 \cdot 157 \cdot 379\)\\                                                                        
  \hline
  \(617\) & \(2 \cdot 3 \cdot 103\) & \(683\) & \(2^{2} \cdot 3^{2} \cdot 19\) & \(761\) & \(2 \cdot 3 \cdot 127\)\\
  \(617^{2}\) & \(97 \cdot 3931\) & \(683^{2}\) & \(7 \cdot 66739\) & \(761^{2}\) & \(579883\)\\
  \(617^{3}\) & \(2^{2} \cdot 3 \cdot 5 \cdot 103 \cdot 38069\) & \(683^{3}\) & \(2^{3} \cdot 3^{2} \cdot 5 \cdot 19 \cdot 46649\) & \(761^{3}\) & \(2^{2} \cdot 3 \cdot 17 \cdot 127 \cdot 17033\)\\                                                                      
  \hline
  \(619\) & \(2^{2} \cdot 5 \cdot 31\) & \(691\) & \(2^{2} \cdot 173\) & \(769\) & \(2 \cdot 5 \cdot 7 \cdot 11\)\\
  \(619^{2}\) & \(3 \cdot 19 \cdot 6733\) & \(691^{2}\) & \(3 \cdot 19 \cdot 8389\) & \(769^{2}\) & \(3 \cdot 31 \cdot 6367\)\\
  \(619^{3}\) & \(2^{3} \cdot 5 \cdot 13 \cdot 31 \cdot 14737\) & \(691^{3}\) & \(2^{3} \cdot 173 \cdot 193 \cdot 1237\) & \(769^{3}\) & \(2^{2} \cdot 5 \cdot 7 \cdot 11 \cdot 17 \cdot 17393\)\\                                                                         
  \hline
  \(631\) & \(2^{3} \cdot 79\) & \(701\) & \(2 \cdot 3^{3} \cdot 13\) & \(773\) & \(2 \cdot 3^{2} \cdot 43\)\\
  \(631^{2}\) & \(3 \cdot 307 \cdot 433\) & \(701^{2}\) & \(492103\) & \(773^{2}\) & \(598303\)\\
  \(631^{3}\) & \(2^{4} \cdot 79 \cdot 199081\) & \(701^{3}\) & \(2^{2} \cdot 3^{3} \cdot 13 \cdot 17 \cdot 97 \cdot 149\) & \(773^{3}\) & \(2^{2} \cdot 3^{2} \cdot 5 \cdot 43 \cdot 59753\)\\                                                                            
  \hline
  \(641\) & \(2 \cdot 3 \cdot 107\) & \(709\) & \(2 \cdot 5 \cdot 71\) & \(787\) & \(2^{2} \cdot 197\)\\
  \(641^{2}\) & \(7 \cdot 58789\) & \(709^{2}\) & \(3 \cdot 7 \cdot 23971\) & \(787^{2}\) & \(3 \cdot 37^{2} \cdot 151\)\\
  \(641^{3}\) & \(2^{2} \cdot 3 \cdot 107 \cdot 205441\) & \(709^{3}\) & \(2^{2} \cdot 5 \cdot 37 \cdot 71 \cdot 6793\) & \(787^{3}\) & \(2^{3} \cdot 5 \cdot 197 \cdot 241 \cdot 257\)\\                                                                                  
  \hline
  \(643\) & \(2^{2} \cdot 7 \cdot 23\) & \(719\) & \(2^{4} \cdot 3^{2} \cdot 5\) & \(797\) & \(2 \cdot 3 \cdot 7 \cdot 19\)
  \\                                                                                                                                  
  \(643^{2}\) & \(3 \cdot 97 \cdot 1423\) & \(719^{2}\) & \(487 \cdot 1063\) & \(797^{2}\) & \(157 \cdot 4051\)\\
  \(643^{3}\) & \(2^{3} \cdot 5^{2} \cdot 7 \cdot 23 \cdot 8269\) & \(719^{3}\) & \(2^{5} \cdot 3^{2} \cdot 5 \cdot 53 \cdot 4877\) & \(797^{3}\) & \(2^{2} \cdot 3 \cdot 5 \cdot 7 \cdot 19 \cdot 63521\)\\                                                               
  \hline
  \(647\) & \(2^{3} \cdot 3^{4}\) & \(727\) & \(2^{3} \cdot 7 \cdot 13\) & \(809\) & \(2 \cdot 3^{4} \cdot 5\)\\
  \(647^{2}\) & \(211 \cdot 1987\) & \(727^{2}\) & \(3 \cdot 176419\) & \(809^{2}\) & \(7 \cdot 13 \cdot 19 \cdot 379\)\\
  \(647^{3}\) & \(2^{4} \cdot 3^{4} \cdot 5 \cdot 41 \cdot 1021\) & \(727^{3}\) & \(2^{4} \cdot 5 \cdot 7 \cdot 13 \cdot 17 \cdot 3109\) & \(809^{3}\) & \(2^{2} \cdot 3^{4} \cdot 5 \cdot 229 \cdot 1429\)\\                                                              
  \hline
  \(653\) & \(2 \cdot 3 \cdot 109\) & \(733\) & \(2 \cdot 367\) & \(811\) & \(2^{2} \cdot 7 \cdot 29\)\\
  \(653^{2}\) & \(7 \cdot 13^{2} \cdot 19^{2}\) & \(733^{2}\) & \(3 \cdot 19 \cdot 9439\) & \(811^{2}\) & \(3 \cdot 31 \cdot 73 \cdot 97\)\\                                                                                                                               
  \(653^{3}\) & \(2^{2} \cdot 3 \cdot 5 \cdot 109 \cdot 42641\) & \(733^{3}\) & \(2^{2} \cdot 5 \cdot 13 \cdot 367 \cdot 4133\) & \(811^{3}\) & \(2^{3} \cdot 7 \cdot 13 \cdot 29 \cdot 41 \cdot 617\)\\                                                                   
  \hline
  \(659\) & \(2^{2} \cdot 3 \cdot 5 \cdot 11\) & \(739\) & \(2^{2} \cdot 5 \cdot 37\) & \(821\) & \(2 \cdot 3 \cdot 137\)\\
  \(659^{2}\) & \(13 \cdot 33457\) & \(739^{2}\) & \(3 \cdot 7 \cdot 26041\) & \(821^{2}\) & \(7 \cdot 229 \cdot 421\)\\
  \(659^{3}\) & \(2^{3} \cdot 3 \cdot 5 \cdot 11 \cdot 17 \cdot 53 \cdot 241\) & \(739^{3}\) & \(2^{3} \cdot 5 \cdot 37 \cdot 273061\) & \(821^{3}\) & \(2^{2} \cdot 3 \cdot 137 \cdot 337021\)\\
  \hline\hline
\end{tabular}

\begin{tabular}{||p{0.6cm}|p{2.9cm}||p{0.6cm}|p{2.9cm}||p{0.6cm}|p{2.9cm}||}
  \hline\hline
  Num. & Divisor sum & Num. & Divisor sum & Num. & Divisor sum \\
  \hline\hline
  \(823\) & \(2^{3} \cdot 103\) & \(881\) & \(2 \cdot 3^{2} \cdot 7^{2}\) & \(947\) & \(2^{2} \cdot 3 \cdot 79\)\\
  \(823^{2}\) & \(3 \cdot 7 \cdot 43 \cdot 751\) & \(881^{2}\) & \(19 \cdot 40897\) & \(947^{2}\) & \(7 \cdot 277 \cdot 463\)\\
  \(823^{3}\) & \(2^{4} \cdot 5 \cdot 103 \cdot 67733\) & \(881^{3}\) & \(2^{2} \cdot 3^{2} \cdot 7^{2} \cdot 388081\) & \(947^{3}\) & \(2^{3} \cdot 3 \cdot 5 \cdot 79 \cdot 89681\)\\                                                                                    
  \hline
  \(827\) & \(2^{2} \cdot 3^{2} \cdot 23\) & \(883\) & \(2^{2} \cdot 13 \cdot 17\) & \(953\) & \(2 \cdot 3^{2} \cdot 53\)\\
  \(827^{2}\) & \(684757\) & \(883^{2}\) & \(3 \cdot 260191\) & \(953^{2}\) & \(181 \cdot 5023\)\\
  \(827^{3}\) & \(2^{3} \cdot 3^{2} \cdot 5 \cdot 13 \cdot 23 \cdot 5261\) & \(883^{3}\) & \(2^{3} \cdot 5 \cdot 13 \cdot 17 \cdot 77969\) & \(953^{3}\) & \(2^{2} \cdot 3^{2} \cdot 5 \cdot 53 \cdot 90821\)\\                                                            
  \hline
  \(829\) & \(2 \cdot 5 \cdot 83\) & \(887\) & \(2^{3} \cdot 3 \cdot 37\) & \(967\) & \(2^{3} \cdot 11^{2}\)\\
  \(829^{2}\) & \(3 \cdot 211 \cdot 1087\) & \(887^{2}\) & \(13 \cdot 60589\) & \(967^{2}\) & \(3 \cdot 67 \cdot 4657\)\\
  \(829^{3}\) & \(2^{2} \cdot 5 \cdot 17^{2} \cdot 29 \cdot 41 \cdot 83\) & \(887^{3}\) & \(2^{4} \cdot 3 \cdot 5 \cdot 29 \cdot 37 \cdot 2713\) & \(967^{3}\) & \(2^{4} \cdot 5 \cdot 11^{2} \cdot 13 \cdot 7193\)\\                                                      
  \hline
  \(839\) & \(2^{3} \cdot 3 \cdot 5 \cdot 7\) & \(907\) & \(2^{2} \cdot 227\) & \(971\) & \(2^{2} \cdot 3^{5}\)\\
  \(839^{2}\) & \(704761\) & \(907^{2}\) & \(3 \cdot 7 \cdot 39217\) & \(971^{2}\) & \(13 \cdot 79 \cdot 919\)\\
  \(839^{3}\) & \(2^{4} \cdot 3 \cdot 5 \cdot 7 \cdot 109 \cdot 3229\) & \(907^{3}\) & \(2^{3} \cdot 5^{2} \cdot 227 \cdot 16453\) & \(971^{3}\) & \(2^{3} \cdot 3^{5} \cdot 197 \cdot 2393\)\\                                                                            
  \hline
  \(853\) & \(2 \cdot 7 \cdot 61\) & \(911\) & \(2^{4} \cdot 3 \cdot 19\) & \(977\) & \(2 \cdot 3 \cdot 163\)\\
  \(853^{2}\) & \(3 \cdot 43 \cdot 5647\) & \(911^{2}\) & \(830833\) & \(977^{2}\) & \(7 \cdot 136501\)\\
  \(853^{3}\) & \(2^{2} \cdot 5 \cdot 7 \cdot 13 \cdot 29 \cdot 61 \cdot 193\) & \(911^{3}\) & \(2^{5} \cdot 3 \cdot 19 \cdot 29 \cdot 41 \cdot 349\) & \(977^{3}\) & \(2^{2} \cdot 3 \cdot 5 \cdot 53 \cdot 163 \cdot 1801\)\\                                            
  \hline
  \(857\) & \(2 \cdot 3 \cdot 11 \cdot 13\) & \(919\) & \(2^{3} \cdot 5 \cdot 23\) & \(983\) & \(2^{3} \cdot 3 \cdot 41\)\\
  \(857^{2}\) & \(735307\) & \(919^{2}\) & \(3 \cdot 7 \cdot 13 \cdot 19 \cdot 163\) & \(983^{2}\) & \(103 \cdot 9391\)\\
  \(857^{3}\) & \(2^{2} \cdot 3 \cdot 5^{2} \cdot 11 \cdot 13 \cdot 37 \cdot 397\) & \(919^{3}\) & \(2^{4} \cdot 5 \cdot 23 \cdot 37 \cdot 101 \cdot 113\) & \(983^{3}\) & \(2^{4} \cdot 3 \cdot 5 \cdot 13 \cdot 41 \cdot 7433\)\\                                        
  \hline
  \(859\) & \(2^{2} \cdot 5 \cdot 43\) & \(929\) & \(2 \cdot 3 \cdot 5 \cdot 31\) & \(991\) & \(2^{5} \cdot 31\)\\
  \(859^{2}\) & \(3 \cdot 246247\) & \(929^{2}\) & \(157 \cdot 5503\) & \(991^{2}\) & \(3 \cdot 7 \cdot 13^{2} \cdot 277\)\\
  \(859^{3}\) & \(2^{3} \cdot 5 \cdot 43 \cdot 137 \cdot 2693\) & \(929^{3}\) & \(2^{2} \cdot 3 \cdot 5 \cdot 31 \cdot 431521\) & \(991^{3}\) & \(2^{6} \cdot 31 \cdot 491041\)\\
  \hline
  \(863\) & \(2^{5} \cdot 3^{3}\) & \(937\) & \(2 \cdot 7 \cdot 67\) & \(997\) & \(2 \cdot 499\)\\
  \(863^{2}\) & \(7^{2} \cdot 15217\) & \(937^{2}\) & \(3 \cdot 292969\) & \(997^{2}\) & \(3 \cdot 13 \cdot 31 \cdot 823\)\\
  \(863^{3}\) & \(2^{6} \cdot 3^{3} \cdot 5 \cdot 13 \cdot 17 \cdot 337\) & \(937^{3}\) & \(2^{2} \cdot 5 \cdot 7 \cdot 67 \cdot 87797\) & \(997^{3}\) & \(2^{2} \cdot 5 \cdot 499 \cdot 99401\)\\
  \hline
  \(877\) & \(2 \cdot 439\) & \(941\) & \(2 \cdot 3 \cdot 157\) &  & \\
  \(877^{2}\) & \(3 \cdot 7 \cdot 37 \cdot 991\) & \(941^{2}\) & \(811 \cdot 1093\) &  & \\
  \(877^{3}\) & \(2^{2} \cdot 5 \cdot 439 \cdot 76913\) & \(941^{3}\) & \(2^{2} \cdot 3 \cdot 13 \cdot 157 \cdot 34057\) & &\\                                                                                                                             
  \hline\hline
\end{tabular}
\end{center}
\normalsize